\documentclass[draftcls,onecolumn]{IEEEtran}
\usepackage{graphicx}
\usepackage{subfigure}
\usepackage{amsmath}
\usepackage{amssymb}
\usepackage{theorem}
\usepackage[usenames,dvipsnames]{pstricks}
\usepackage{pstricks-add}
\usepackage{pst-grad}
\usepackage{pst-plot} 
\usepackage{pst-eucl} 
\usepackage[eulergreek]{sansmath}

\usepackage{caption} 
\usepackage{cite}

\usepackage{algorithm}
\usepackage{algorithmicx}
\usepackage{algpseudocode}

\newcommand{\qed}{\nobreak \ifvmode \relax \else%
      \ifdim\lastskip<1.5em \hskip-\lastskip%
      \hskip1.5em plus0em minus0.5em \fi \nobreak%
      \vrule height0.75em width0.5em depth0.25em\fi}
\newcommand{\thh}{\begin{sansmath}\varphi\end{sansmath}}
\newcommand{\menge}[2]{\big\{{#1}~\big |~{#2}\big\}} 
\newcommand{\emp}{\ensuremath{\varnothing}}
\newcommand{\scal}[2]{{#2}^{\!\top}{#1}} 
  
\newcommand{\HH}{\RR^N}
\newcommand{\GG}{\RR^K}

\newcommand{\RR}{\ensuremath{\mathbb R}}
\newcommand{\CC}{\ensuremath{\mathbb C}}
\newcommand{\ZZ}{\ensuremath{\mathbb Z}}
\newcommand{\NN}{\ensuremath{\mathbb N}}
\newcommand{\RP}{\ensuremath{\left[0,+\infty\right[}}

\newcommand{\RPP}{\ensuremath{\,\left]0,+\infty\right[}}
\newcommand{\RX}{\ensuremath{\,\left]-\infty,+\infty\right]}}

\newcommand{\dom}{\ensuremath{\mathrm{dom}\,}}
\newcommand{\prox}{\ensuremath{\mathrm{prox}}}
\newcommand{\inte}{\ensuremath{\mathrm{int}\,}}

\newcommand{\Id}{\ensuremath{\mathrm{Id}}}

\newcommand{\pinf}{\ensuremath{+\infty}}

\newcommand{\xb}{\ensuremath{\boldsymbol{x}}}
\newcommand{\vb}{\ensuremath{\boldsymbol{v}}}

\newcommand{\epi}{\operatorname{epi}}
\newcommand{\argmind}[2]{\ensuremath{\underset{\substack{{#1}}}%
{\mathrm{argmin}}\;\;#2 }}
\newcommand{\Argmind}[2]{\ensuremath{\underset{\substack{{#1}}}%
{\mathrm{Argmin}}\;\;#2 }}
\newcommand{\infconv}{\ensuremath{\mbox{\footnotesize$\,\square\,$}}}
\newcommand{\minimize}[2]{\ensuremath{\underset{\substack{{#1}}}%
{\mathrm{minimize}}\;\;#2 }}
\newcommand{\maximize}[2]{\ensuremath{\underset{\substack{{#1}}}%
{\mathrm{maximize}}\;\;#2 }}
\newcommand{\rank}{\operatorname{rank}}
\newcommand{\ran}{\ensuremath{\text{\rm ran}\,}}

\theoremstyle{plain}{\theorembodyfont{\rmfamily}%
}
\theoremstyle{plain}{\theorembodyfont{\rmfamily}%
}
\theoremstyle{plain}{\theorembodyfont{\rmfamily}%
}
\theoremstyle{plain}{\theorembodyfont{\rmfamily}%
}
\theoremstyle{plain}{\theorembodyfont{\rmfamily}%
}
\theoremstyle{plain}{\theorembodyfont{\rmfamily}%
\renewcommand\theenumi{\roman{enumi})}
\renewcommand\theenumii{\alph{enumii})}

\newcounter{numperceuse}
\renewcommand\thenumperceuse{\roman{numperceuse}}

\usepackage{color}
\abovecaptionskip=-3pt
\belowcaptionskip=-15pt
\begin{document}

\title{Playing with Duality: An Overview of Recent Primal-Dual Approaches
for Solving Large-Scale Optimization Problems}
\author{Nikos Komodakis, {\em Member, IEEE}, 
and Jean-Christophe Pesquet, {\em Fellow, IEEE}
%\thanks{This work was supported by }
\thanks{N. Komodakis (corresponding author) and J.-C. Pesquet are with the Laboratoire d'Informatique Gaspard 
Monge, UMR CNRS 8049,
Universit\'e Paris-Est, 77454 Marne la Vall\'ee 
Cedex 2, France. E-mail: \texttt{nikos.komodakis@enpc.fr}, \texttt{jean-christophe.pesquet@univ-paris-est.fr}.}
%\thanks{\bfseries\em 22 January, 2014 -- version 1.0}
}
\markboth{IEEE Signal Processing Magazine}%
{Komodakis and Pesquet: Playing with Duality}

\thispagestyle{empty}
\maketitle

\vspace*{-1.5cm}
\begin{abstract}
Optimization methods are at the core of many problems in signal/image processing, computer vision, and machine learning.
For a long time, it has been recognized that looking at the dual of an optimization problem may drastically simplify its solution.
Deriving efficient strategies which jointly brings into play the primal and the dual problems is however a more recent idea which has generated many important new contributions in the last years. These novel developments are grounded on recent advances in convex analysis, discrete optimization, parallel processing, and nonsmooth optimization with emphasis on sparsity issues.
In this paper, we aim at presenting the principles of primal-dual approaches, while giving an overview of numerical methods which have been proposed in different contexts. We show the benefits which can be drawn from primal-dual algorithms both for solving large-scale convex optimization problems and discrete ones, and we provide various application examples to illustrate their usefulness.
\end{abstract}


\newpage
\setcounter{page}{1}
\section{Motivation and importance of the topic}

Optimization \cite{Bertsekas_D_2004_book_non_prog} is an extremely popular paradigm which constitutes the backbone of many branches of applied mathematics
and engineeering, such as signal processing, computer vision, machine learning, inverse problems, and network communications, to mention just a few. The popularity of optimization approaches often stems from the fact that many problems from the above fields are typically characterized by a lack of closed form solutions and by uncertainties. In signal and image processing, for instance, uncertainties can be introduced due to noise, sensor imperfectness, or ambiguities that are often inherent in the visual interpretation. As a result, perfect or exact solutions hardly exist, whereas inexact but optimal (in a statistical or an application-specific sense) solutions and their efficient computation is what one aims at. At the same time, one important characteristic that is  nowadays shared by   increasingly many optimization problems encountered in the above areas  is the fact that these problems are often of very large scale.
A good example is the field of computer vision where one often needs to solve low level problems that require  associating at least one (and typically more than one) variable to each pixel of an image (or even worse of an image sequence as in the case of video) \cite{Blake_book}. This leads to  problems that easily  can contain millions of variables, which are therefore the norm rather than the exception in this context.

Similarly,   in fields like machine learning \cite{Sra_S_2012_book_optimization_ml,Theodoridis_S_2014_book_machine_lsipp}, due to the great ease with which data can now be collected and stored,  quite often one has to cope with truly massive datasets and to train very large models, which thus naturally lead to optimization problems of very high dimensionality \cite{Bach_F_2012_j-ftml_optimization_sip}.
Of course, a similar situation arises in many other scientific domains, including application areas such as 
inverse problems (e.g., medical image reconstruction or satellite image restoration)
or telecommunications (e.g., network design, network provisioning) and industrial engineering. Due to this fact,  computational efficiency  constitutes a major issue that needs to be thoroughly addressed. This, therefore, makes mandatory the use of tractable  optimization techniques that are able to properly exploit the problem structures, but which at the same time  remain  applicable to a class of problems as wide as possible.   

A bunch of  important advances that took place in this regard over the last years concerns a particular class of optimization approaches known as \emph{primal-dual} methods.   As their name implies, these approaches proceed by concurrently solving a primal problem (corresponding to the original optimization task) as well as a dual formulation of this problem.
As it turns out, in doing so  they are able to exploit more efficiently  the problem specific properties, thus offering in many cases  important computational advantages, some of which are briefly mentioned next for two very broad classes of problems. 

\subsubsection{Convex optimization} Primal-dual methods have been primarily employed in convex optimization problems \cite{Rockafellar_R_T_1970_book_Convex_A,Boyd_S_2004_book_con_optim,Bauschke_H_2011_book_con_amo} where strong duality holds. They have been successfully applied  to various types of nonlinear and nonsmooth cost functions that are prevalent in the above-mentioned application fields.

Many such applied problems    
can essentially be expressed under the form of a minimization of a sum of terms, where each term is given by the composition of a convex  function  with a linear operator.
One first advantage of primal-dual methods pertains to the fact that they can yield very efficient splitting optimization schemes, according to which a solution to the original problem is iteratively computed through solving   a sequence of easier subproblems, each one involving only one of the  terms appearing in the objective function.

The  resulting primal-dual splitting schemes  can also handle    both differentiable and  nondifferentiable terms, the former by use of  gradient operators (i.e., through explicit steps) and the latter by use of proximity operators (i.e., through implicit steps) \cite{Moreau_J_1965_bsmf_Proximite_eddueh,Combettes_P_2010_inbook_proximal_smsp}.  Depending on the target functions, either explicit or implicit steps may be easier to implement. Therefore, the derived optimization schemes exploit the properties of the input problem, in a flexible manner, thus leading to very efficient first-order algorithms.

Even more importantly, primal-dual techniques are able to achieve what is known as \emph{full splitting} in the optimization literature, meaning that each of the operators  involved in the problem (i.e., not only the gradient or proximity
operators but also the involved linear operators) is used separately \cite{Combettes_P_2012_j-svva_pri_dsa}. As a result, no call to the inversion of a linear operator, which is  an expensive operation  for large scale problems, is required during  the optimization process. This is an important 
feature which gives these methods a significant  computational  advantage compared with all other splitting-based  approaches.

Last but not least, primal-dual methods lead to algorithms that are easily parallelizable, which is nowadays becoming increasingly important for efficiently handling high-dimensional problems.

\subsubsection{Discrete optimization}
Besides convex optimization, another important area where primal-dual methods play a prominent role is discrete optimization. This is of particular significance given that a large variety of  tasks from signal processing, computer vision, and pattern recognition
are  formulated as discrete labeling problems, where one seeks to optimize some measure related to
the quality of the labeling \cite{Li_book}. This includes, for instance, tasks such as image segmentation, optical flow estimation, image denoising, stereo matching, to mention a few examples from image analysis. The resulting discrete optimization problems not only are of very large size, but also typically exhibit highly nonconvex objective functions, which are generally 
intricate to optimize.

Similarly to the case of convex optimization, primal-dual methods  again offer many computational advantages, leading often to very fast graph-cut or message-passing-based algorithms, which are also easily parallelizable, thus providing in many cases a very efficient way   for handling    discrete optimization problems that are encountered in practice \cite{komodakis_CVIU08, komodakis_pami2011, TRW_wainwright, Wang}.
Besides being efficient, they are also successful in making little compromises regarding the quality of the estimated solutions. Techniques like the so-called \emph{primal-dual schema} are known to provide
a principled way for deriving powerful approximation algorithms to difficult combinatorial problems, thus allowing primal-dual methods to often exhibit  theoretical (i.e., worst-case) 
approximation properties. 
Furthermore, apart from the aforementioned worst-case guaranties, primal-dual  algorithms  can also provide  (for free) \emph{per-instance} approximation guaranties. This is essentially 
made possible by the fact that these methods are estimating not only  primal but also  dual solutions. 
%  
%  Such factors, which can  be continuously updated throughout an algorithm execution (see Fig.~\ref{fig:fastpd_bounds}), can be directly used for assessing the performance of a primal-dual method with respect to any particular  problem instance (and without essentially any extra computational cost).  Moreover, often in practice they turn out to take values very close to 1, which means that the corresponding estimated solutions are almost optimal.

Convex optimization and discrete optimization have different background theory originally.
Convex optimization may appear as the most
tractable topic in optimization, for which many efficient algorithms have been developed 
allowing a broad class of problems to be solved. 
By contrast, combinatorial optimization problems are generally NP-hard.
However, many convex relaxations of certain discrete problems can provide good
approximate solutions to the original ones \cite{Vazirani_book, Hochbaum_book}. The problems encountered in discrete optimization 
therefore constitute a source of inspiration for developing novel convex optimization techniques.

\vspace*{0.5em}
\noindent \textbf{Goals of this tutorial paper.}
Based on the above observations, our objectives will be the following:

\begin{enumerate}
\item To provide a thorough introduction that intuitively explains the basic principles and ideas behind primal-dual approaches.
\item To describe how these methods can be employed both in the context of continuous optimization  and in the context of discrete optimization.
\item To  explain some of the recent advances that have taken place concerning primal-dual algorithms for solving large-scale optimization problems.
\item To detail useful connections between primal-dual methods and some widely used
optimization techniques like the alternating direction method of multipliers (ADMM) \cite{Fortin_M_1983_book_augmented_lmansbvp,Boyd_S_2011_j-found-tml_distributed_osl_admm}.

\item Finally, to  provide examples of  useful applications  in the context of image analysis and signal processing.

\end{enumerate}

The remainder of the paper is structured as follows. In Section \ref{se:background}, we introduce the necessary methodological background on optimization.
Our presentation is grounded on the powerful notion of duality
known as Fenchel's duality, from which duality properties in linear programming  can be deduced. We also introduce useful tools from functional analysis and convex optimization, including the notions of subgradient and subdifferential, conjugate function, and proximity operator.
The following two sections explain and describe various primal-dual methods. Section \ref{se:optconv} is devoted to convex optimization  problems.
We discuss the merits of various algorithms and explain their connections with 
ADMM, that we show to be a special case of primal-dual proximal method. 
Section \ref{se:optdisc} deals with primal-dual methods for discrete optimization. We explain how to derive algorithms of this type based on the primal-dual schema 
which is a well-known approximation technique in combinatorial optimization, and we also present primal-dual methods based on LP relaxations and dual decomposition. 
In Section \ref{se:appli}, we present applications from the domains of signal processing and image analysis, including inverse problems and computer vision tasks related to Markov Random Field energy minimization. In Section \ref{se:conclu}, we finally conclude the tutorial with a brief summary and discussion.

\section{Optimization background} \label{se:background}
In this section, we introduce the necessary mathematical definitions and concepts used for introducing primal-dual algorithms in later sections. Although the following framework holds for general Hilbert spaces, for simplicity we will focus on the finite dimensional case.
\subsection{Notation}
\begin{figure}[htb]
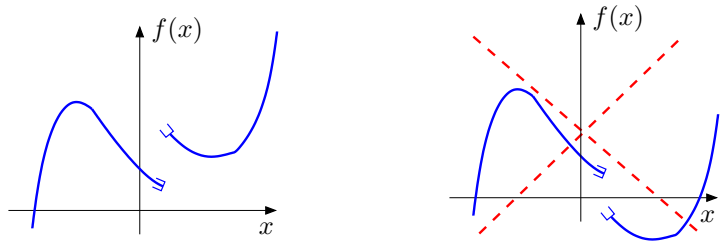

\begin{center}
\begin{tabular}{@{}c@{}c@{}}
\scalebox{1.0}{\input{sci_ex1.pstex_t}} \hspace{1cm} & \hspace{1cm}  \scalebox{1.0}{\input{sci_ex2w.pstex_t}}
\end{tabular}
\end{center}
\caption{\small Illustration of the lower-semicontinuity property.}\label{fig:lowersemi}
\end{figure}

In this paper, we will consider functions from $\HH$ to $\RX$.
The fact that we allow functions to take $+\infty$ value is useful in modern optimization to
discard some ``forbidden part" of the space when searching for an optimal solution
(for example, in image processing problems, the components of the solution often are intensity values which must be nonnegative).
The \emph{domain} of a function $f\colon \HH \to \RX$ is the subset of $\HH$ where this function takes finite values, i.e.
$\dom f = \menge{x\in \HH}{f(x) < +\infty}$. A function with a nonempty domain is said to be
\emph{proper}. A function $f$ is \emph{convex} if 
\begin{equation}
\big(\forall (x,y)\in (\HH)^2\big)
(\forall \lambda \in [0,1])\qquad f(\lambda x + (1-\lambda) y) \le \lambda f(x) + (1-\lambda) f(y).
\end{equation}
The class of functions for which most of the main results in convex analysis have been established
is $\Gamma_0(\HH)$, the class of proper, convex, lower-semicontinuous functions from
$\HH$ to $\RX$. Recall that a function  $f\colon \HH \to \RX$ is lower-semicontinuous if its \emph{epigraph}
$\epi f = \menge{(x,\zeta) \in \dom f \times \RR}{f(x) \le \zeta}$ is a closed set (see Fig. \ref{fig:lowersemi}).

If $C$ is a nonempty subset of $\HH$, the \emph{indicator function} of $C$ is defined as
\begin{equation}\label{e:defindicator}
(\forall x \in \HH) \qquad \iota_C(x) = \begin{cases}
0 & \mbox{if $x\in C$}\\
+\infty & \mbox{otherwise.}
\end{cases}
\end{equation}
This function belongs to $\Gamma_0(\HH)$ if and only if $C$ is a nonempty closed convex set.
%The convexity of $C$ means that $\big(\forall (x,y)\in (\HH)^2\big)$
%$(\forall \lambda \in [0,1])$ $\lambda x + (1-\lambda) y\in C$.

\begin{figure}[htb]
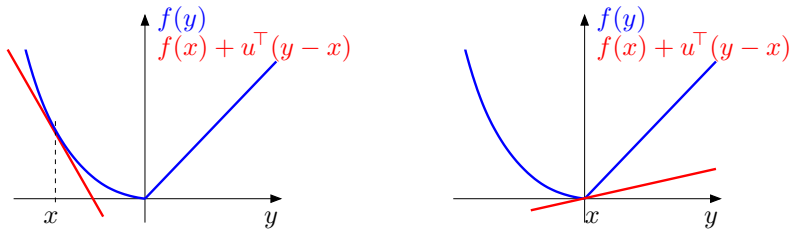

\begin{center}
\begin{tabular}{@{}c@{}c@{}}
\scalebox{1.0}{\input{sousdiff1.pstex_t}} \hspace{1cm} & \hspace{1cm}  \scalebox{1.0}{\input{sousdiff4.pstex_t}}
\end{tabular}
\end{center}
\caption{\small Examples of subgradients $u$ of a function $f$ at $x$.}\label{fig:subgrad}
\end{figure}

The Moreau \emph{subdifferential} of a function $f\colon \HH \to \RX$ at $x\in \HH$ is defined as
\begin{equation}
\partial f(x) = \menge{u\in\HH}{(\forall y \in \HH)\;\; f(y) \ge f(x)+\scal{(y-x)}{u}}.
\end{equation}
Any vector $u$ in $\partial f(x)$ is called a \emph{subgradient} of $f$ at $x$ (see Fig. \ref{fig:subgrad}).

Fermat's rule states that $0$ is a subgradient of $f$ at $x$ if and only if $x$
belongs to the set %$\Argmin f$  
of global minimizers of $f$.
If $f$ is a proper convex function which is differentiable at $x$, then its subdifferential
at $x$ reduces to the singleton consisting of its gradient, i.e. $\partial f(x) = \{\nabla f(x)\}$.
Note that, in the nonconvex case, extended definitions of the subdifferential may be useful such as
the \emph{limiting subdifferential} \cite{Mordukhovich06}, but this one reduces to the Moreau subdifferential when the function
is convex.

\subsection{Proximity operator}
\begin{figure}[htb]
\begin{center}
\includegraphics[width=7.5cm]{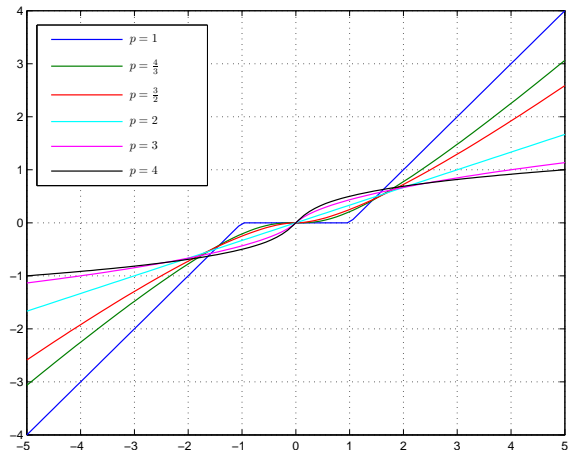}
\end{center}
\vspace{-0.2cm}
\begin{center}
\caption{\small Graph of $\prox_{|\cdot|^p}$. 
This   power $p$ function is often used to regularize inverse problems.} \label{fig:proxpowerp}
\end{center}
\end{figure}

A concept which has been of growing importance in recent developments in optimization is the concept of \emph{proximity
operator}. It must be pointed out that the proximity operator was introduced in the early work by J.~J.~Moreau (1923-2014) \cite{Moreau_J_1965_bsmf_Proximite_eddueh}. The proximity operator of a function $f\in \Gamma_0(\HH)$
is defined as
\begin{equation}
\prox_f\colon \HH \to \HH\colon x \mapsto \argmind{y\in \HH}{f(y) + \frac12 \|y-x\|^2}
\end{equation}
where $\|\cdot\|$ denotes the Euclidean norm.
For every $x\in \HH$, $\prox_f x$ can thus be interpreted as the result of a regularized minimization of $f$ 
in the neighborhood of $x$. Note that the minimization to be performed to calculate $\prox_f x$ always has a unique solution.
Fig. \ref{fig:proxpowerp} shows the variations of the $\prox_f$ function when $f\colon \RR \to \RR\colon x \mapsto |x|^p$ with $p\ge 1$. In the case when $p=1$, the classical soft-thesholding operation is obtained.

In the case when $f$ is equal to the indicator function of a nonempty closed convex set $C\subset \HH$, the proximity operator of $f$
reduces to the projection $P_C$ onto this set, i.e. $(\forall x \in \HH)$ $P_C x = \argmind{y\in C}{\|y-x\|}$.

This shows that proximity operators can be viewed as extensions of projections onto convex sets. The proximity operator enjoys many
properties of the projection, in particular it is firmly nonexpansive.
The firm nonexpansiveness can be viewed as a generalization of the strict contraction property 
which is the engine behind the Banach-Picard fixed point theorem. This property makes the proximity operator successful 
in ensuring the convergence of fixed point algorithms grounded on its use. For more details about proximity operators and their rich properties,
the reader is refered to the tutorial papers in \cite{Combettes_P_2010_inbook_proximal_smsp,Bach_F_2012_j-ftml_optimization_sip,Parikh_N_2013_j-found-tml_prox_algo}.
The definition of the proximity operator can be extended to nonconvex lower-semicontinuous functions which are lower bounded by an affine function, but $\prox_f x$ is no longer guaranteed to be uniquely defined at any given point $x$.

\subsection{Conjugate function}
\begin{figure}[htb]
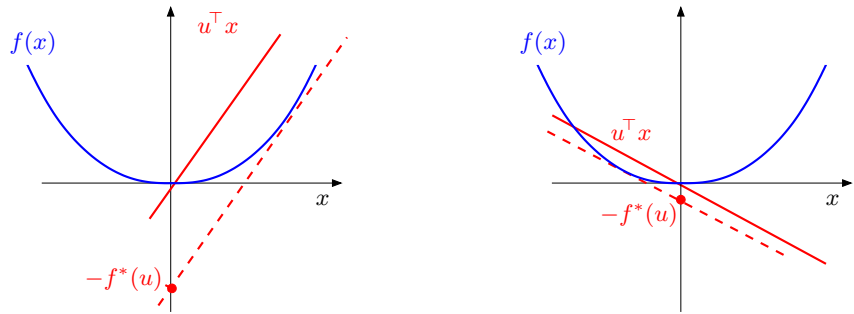

\begin{center}
\begin{tabular}{@{}c@{}c@{}}
\scalebox{1.0}{\input{conj1.pstex_t}} \hspace{1cm} & \hspace{1cm}  \scalebox{1.0}{\input{conj3.pstex_t}}
\end{tabular}
\end{center}
\caption{\small Graphical interpretation of the conjugate function.}\label{fig:conj}
\end{figure}

A fundamental notion when dealing with duality issues is the notion of \emph{conjugate function}.
The conjugate of a function $f\colon \HH \to \RX$ is the function $f^*$ defined as
\begin{equation}\label{e:deffconj}
f^*\colon \HH \to \RX\colon
u \mapsto \sup_{x\in \HH} \big(\scal{u}{x} -f(x)\big).
\end{equation}
This concept was introduced by A. M.~Legendre (1752-1833) in the one-variable case, and it was generalized by M.~W.~Fenchel (1905-1988).
A graphical illustration of the conjugate function is provided in Fig. \ref{fig:conj}.
In particular, for every vector $x \in \HH$ such that the supremum in \eqref{e:deffconj} is attained, $u$ is a subgradient of $f$ at $x$.

It must be emphasized that, even if $f$ is nonconvex, $f^*$ is a (non necessarily proper) lower-semicontinuous convex function.
In addition, when $f \in \Gamma_0(\HH)$, then $f^* \in \Gamma_0(\HH)$, and also the biconjugate of $f$ (that is the conjugate of its conjugate)
is equal to $f$. This means that we can express any function $f$ in $\Gamma_0(\HH)$ as
\begin{equation}
(\forall x \in \HH)\qquad f(x) = \sup_{u\in \HH} \big(\scal{x}{u} -f^*(u)\big).
\end{equation}
A geometrical interpretation of this result is that the epigraph of any proper lower-semicontinuous convex function always is an intersection of closed half-spaces.

As we have seen, the subdifferential plays an important role in the characterization of the minimizers of a function. A natural question is thus to enquire about the relations existing between the subdifferential of a function $f\colon \HH \to \RX$
and the subdifferential of its conjugate function. An answer is provided by the following important properties:
\begin{align}\label{e:subdifconj}
u\in \partial f(x)\quad &\Rightarrow \quad x \in \partial f^*(u) \qquad \mbox{if $f$ is proper}\nonumber\\
u\in \partial f(x)\quad &\Leftrightarrow \quad x \in \partial f^*(u) \qquad \mbox{if $f \in \Gamma_0(\HH)$.}
\end{align}
Another important property is Moreau's decomposition formula which links the proximity operator of a function
$f\in \Gamma_0(\HH)$ to the proximity operator of its conjugate:
\begin{equation}\label{e:moreaudec}
(\forall x \in \HH)(\forall \gamma \in \RPP) \qquad x = \prox_{\gamma f} x + \gamma\,\prox_{\gamma^{-1} f^*} (\gamma^{-1} x).
\end{equation} 
Other useful properties of the conjugation operation are listed in Table~\ref{t:propconj},\footnote{Throughout the paper, $\inte S$ denotes the interior of a set $S$.} where a parallel is drawn with the multidimensional Fourier transform, which is a more familiar tool in signal and image processing.
Conjugation also makes it possible to build an insightful bridge between the main two kinds of nonsmooth convex functions encountered in signal and image processing problems, namely indicator functions of feasibility constraints and sparsity measures (see framebox below.
%on page \pageref{bo:consup}).

\refstepcounter{numperceuse}
\begin{table}[hptb]
\centering
\caption{\small Parallelism between properties of the Legendre-Fenchel conjugation \cite{Combettes_P_2010_inbook_proximal_smsp} and of the Fourier transform.
$f$ is a function defined on $\HH$, $f^*$ denotes its conjugate, $\widehat{f}$  is its Fourier transform such that
$\widehat{f}(\nu) = \int_{\HH} f(x) \exp(-\jmath 2\pi \scal{\nu}{x}) dx$ where $\nu\in \HH$ and 
$\jmath$ is the imaginary unit (a similar notation is used for other functions), $h$, $g$, and $(f_m)_{1\le m \le M}$ 
are functions defined on $\HH$, $(\varphi_j)_{1\le j \le N}$ are functions
defined on $\RR$, $\psi$  is an even function defined on $\RR$, $\widetilde{\psi}$ is defined as 
$\widetilde{\psi}(\rho) = 2\pi \rho^{(2-N)/2}\int_0^{+\infty} r^{N/2} J_{(N-2)/2}(2\pi r \rho) \psi(r) dr$ where
$\rho \in \RR$ and $J_{(N-2)/2}$ is the Bessel function of order $(N-2)/2$, and $\delta$ denotes the Dirac distribution. (Some properties of the Fourier transform may require some technical assumptions.)
}\label{t:propconj}
\vspace{0.4cm}
\begin{tabular}{|p{0.4cm} l|l|l||l|l|}
\cline{3-6}
 \multicolumn{2}{c}{ } & \multicolumn{2}{|c||}{conjugation} & \multicolumn{2}{c|}{Fourier transform}\\
\hline
& Property & $h(x)$ & $h^*(u)$ & $h(x)$ & $\widehat{h}(\nu)$\\
\hline
\hline
\thenumperceuse\refstepcounter{numperceuse} &
invariant function & $\frac12 \|x\|^2$ & $\frac12 \|u\|^2$ & $\exp(-\pi \|x\|^2) $ &  $\exp(-\pi \|\nu\|^2)$\\
\hline
\thenumperceuse\refstepcounter{numperceuse} &
translation & $f(x-c)$ & $f^*(u)+\scal{u}{c}$ & $f(x-c)$ & $\exp(-\jmath 2\pi \scal{\nu}{c} ) \widehat{f}(\nu)$\\
 & $c\in \HH$           &                      &          &  &\\
\hline
\thenumperceuse\refstepcounter{numperceuse} &
dual translation & $f(x)+\scal{x}{c}$ & $f^*(u-c)$ & $\exp(\jmath 2\pi \scal{x}{c} )f(x-c)$ & $\widehat{f}(\nu-c)$\\
& $c\in \HH$           &                      &          & & \\
\hline
\thenumperceuse\refstepcounter{numperceuse} &
scalar multiplication & $\alpha f(x)$ & $\alpha f^*\left(\frac{u}{\alpha}\right)$ & $\alpha f(x)$ & $\alpha \widehat{f}(\nu)$\\
& $\alpha \in \RPP$           &                      &          & & \\
\hline
\thenumperceuse\refstepcounter{numperceuse} &
invertible linear transform & $f(Lx)$ & $f^*\big((L^{-1})^\top u\big)$ & $f(Lx)$ & $\frac{1}{|\operatorname{det}(L)|}\widehat{f}\big((L^{-1})^\top\nu\big)$\\
& $L \in \RR^{N\times N}$ invertible           &                      &          & & \\
\hline
\thenumperceuse\refstepcounter{numperceuse} &
scaling & $f\left(\frac{x}{\alpha}\right)$ & $f^*(\alpha u)$ & $f\left(\frac{x}{\alpha}\right)$ & $|\alpha| \widehat{f}(\alpha \nu)$\\
& $\alpha \in \RR^*$           &                      &          & & \\
\hline
\thenumperceuse\refstepcounter{numperceuse} &
reflection & $f(-x)$ & $f^*(-u)$ & $f(-x)$ & $\widehat{f}(-\nu)$\\
\hline
\thenumperceuse\refstepcounter{numperceuse} &
separability & $\displaystyle\sum_{j=1}^N \varphi_j(x^{(j)})$ & $\displaystyle\sum_{j=1}^N \varphi_j^*(u^{(j)})$ & $\displaystyle\prod_{j=1}^N \varphi_j(x^{(j)})$ & $\displaystyle\prod_{j=1}^N \widehat{\varphi}_j(\nu^{(j)})$\\
&             &   $x = (x^{(j)})_{1\le j \le N}$  &   $u = (u^{(j)})_{1\le j \le N}$       & $x = (x^{(j)})_{1\le j \le N}$ & 
$\nu = (\nu^{(j)})_{1\le j \le N}$\\
\hline
\thenumperceuse\refstepcounter{numperceuse} \label{e:definfconv} &
isotropy & $\psi(\|x\|)$ & $\psi^*(\|u\|)$ & $\psi(\|x\|)$ & $\widetilde{\psi}(\|\nu\|)$\\
\hline
\thenumperceuse\refstepcounter{numperceuse} \label{e:sumdual} &
inf-convolution & $(f\infconv g)(x)$ & $f^*(u)+g^*(u)$ & $(f\star g)(x)$ & $\widehat{f}(\nu)\widehat{g}(\nu)$\\
&   /convolution   &  $\displaystyle=\inf_{y\in \HH} f(y)+g(x-y)$                     &          & $\displaystyle=\int_{\HH} f(y) g(x-y) dy$ & \\
\hline
\thenumperceuse\refstepcounter{numperceuse}  &
sum/product & $f(x)+g(x)$ & $(f^* \infconv g^*)(u)$ & $f(x)g(x)$ & $(\widehat{f}\star \widehat{g})(\nu)$\\
&             &   $f \in \Gamma_0(\HH)$, $g \in \Gamma_0(\HH)$ &          & & \\
&             &   $\dom f \cap \inte(\dom g) \neq \emp$ &          & & \\
\hline
\thenumperceuse\refstepcounter{numperceuse} &
identity element   & $\iota_{\{0\}}(x)$ & $0$ & $\delta(x)$ & $1$\\
&  of convolution           &                      &          & & \\
\hline
\thenumperceuse\refstepcounter{numperceuse} &
identity element  & $0$ & $\iota_{\{0\}}(u)$ & $1$ & $\delta(\nu)$\\
&  of addition/product           &                      &          & & \\
\hline
\thenumperceuse\refstepcounter{numperceuse} &
offset & $f(x)+\alpha$ & $f^*(u)-\alpha$ & $f(x)+\alpha$ & $\widehat{f}(\nu)+\alpha \delta(\nu)$\\
& $\alpha \in \RR$           &                      &          & & \\
\hline
\thenumperceuse\refstepcounter{numperceuse} &
infinum/sum & $\displaystyle\inf_{1\le m \le M} f_m(x) $ & $\displaystyle\sup_{1\le m \le M} f_m^*(u)$ & $\displaystyle\sum_{m=1}^M f_m(x)$ & 
$\displaystyle\sum_{m=1}^M \widehat{f}_m(\nu)$\\
\hline
\thenumperceuse\refstepcounter{numperceuse} &
value at $0$ & \multicolumn{2}{c||}{$f^*(0) = -\inf f$ } & \multicolumn{2}{c|}{$\displaystyle\widehat{f}(0) = \int_{\HH} f(x) dx$ } \\
\hline
\end{tabular}
\end{table}

\begin{center}
\framebox{\label{bo:consup}
\begin{minipage}{14cm}
\footnotesize\sffamily 
{CONJUGATES OF SUPPORT FUNCTIONS}

The support function of a set $C\subset\HH$ is defined as
\begin{equation}\label{e:defsupfonc}
(\forall u \in \HH)\qquad \sigma_C(u) = \sup_{x\in C} \scal{u}{x}.
\end{equation}
In fact, a function $f$ is the support function of a nonempty closed
convex set $C$ if and only if it belongs to $\Gamma_0(\HH)$ and
it is positively homogeneous \cite{Bauschke_H_2011_book_con_amo}, i.e.
\[
(\forall x \in \HH)(\forall \alpha \in \RPP)\qquad f(\alpha x) = \alpha f(x).
\]
Examples of such functions are norms, e.g. the $\ell_1$-norm:
\[
\big(\forall x = (x^{(j)})_{1\le j \le N} \in \HH\big) \qquad f(x) = \|x\|_1 = \sum_{j=1}^N |x^{(j)}|
\]
 which is a useful convex sparsity-promoting measure in LASSO estimation \cite{Tibshirani_R_1996_j-roy-stat-soc-b_regression_ssv}
and in compressive sensing \cite{Candes_E_J_j-ieee-spm-introduction-rcs_2008}. Another famous
example is the Total Variation semi-norm \cite{Rudin_L_1992_tv_atvmaopiip} which is popular in image processing 
for retrieving constant areas with sharp contours.
An important property is that, if $C$ is a nonempty closed convex set, the conjugate of its support 
function is the indicator function of $C$. For example, the conjugate function
of the $\ell_1$-norm is the indicator function of the hypercube $[-1,1]^N$.
This shows that using sparsity measures are equivalent in the dual domain to imposing some constraints.
\end{minipage}.
}
\end{center}

\subsection{Duality results}\label{se:dualres}
A wide array of problems in signal and image processing can be expressed under the following variational form: 
\begin{equation}\label{e:primal}
\minimize{x\in \HH}{f(x)+g(Lx)}
\end{equation}
where $f\colon \HH \to \RX$, $g\colon \GG\to \RX$, and $L \in \RR^{K\times N}$.
Problem \eqref{e:primal} is usually referred to as the \emph{primal problem}
which is associated with the following \emph{dual problem} \cite{Rockafellar_R_T_1970_book_Convex_A,Bot_R_2010_book_conjugate_dco,Bauschke_H_2011_book_con_amo}:
\begin{equation}\label{e:dual}
\minimize{v\in \GG}{f^*(-L^\top v)+g^*(v)}.
\end{equation}
The latter problem may be easier to solve than the former one, especially when $K$ is much smaller than $N$.

A question however is to know whether solving the dual problem may bring some information on the solution
of the primal one. A first answer to this question is given by the Fenchel-Rockafellar duality theorem
which basically states that solving the dual problem provides a lower bound on the minimum value which can be obtained in the 
primal one. More precisely, 
if $f$ and $g$ are proper functions and if $\mu$ and $\mu^*$ denote the infima of the functions
minimized in the primal and dual problems, respectively, then \emph{weak duality} holds, which means that 
$\mu \ge - \mu^*$. If $\mu$ is finite,
$\mu+\mu^*$ is called the \emph{duality gap}. In addition, if $f\in \Gamma_0(\HH)$ and $g\in \Gamma_0(\GG)$, then, under
appropriate qualification conditions,\footnote{For example, this property is satisfied if the intersection of the 
interior of the domain of $g$ and the image of the domain of $f$ by $L$ is nonempty.} there always exists a solution to the dual
problem and the duality gap vanishes. When the duality gap is equal to zero, it is said that \emph{strong duality} holds.

\begin{center}
\framebox{\label{fb:cons_sharing}
\begin{minipage}{14cm}
\footnotesize\sffamily 
{CONSENSUS AND SHARING ARE DUAL PROBLEMS}

Suppose that our objective is to minimize a composite function $\sum_{m=1}^M g_m$
where the potential $g_m\colon \HH \to \RX$ is computed at the vertex of index $m\in \{1,\ldots,M\}$ of a graph.
A classical technique to perform this task in a distributed or parallel
manner \cite{Boyd_S_2011_j-found-tml_distributed_osl_admm} consists of reformulating this problem as a \emph{consensus problem},
where a variable is assigned to each vertex, and the defined variables $x_1,\ldots,x_M$ are updated so as to reach a consensus:
$x_1 = \ldots = x_M$. This means that, in the product space $(\HH)^M$ the original optimization problem can be rewritten as
\[
\minimize{\xb=(x_1,\ldots,x_M)\in (\HH)^M}{\iota_D(\xb)+\underbrace{\sum_{m=1}^M g_m(x_m)}_{\mbox{$g(\xb)$}}}
\]
where $D$ is the vector space defined as $D = \menge{\xb=(x_1,\ldots,x_M)\in (\HH)^M}{x_1 = \ldots = x_M}$.

By noticing that the conjugate of the indicator function of a vector space is the indicator function of its orthogonal complement,
it is easy to see that the dual of this consensus problem has the following form:
\[
\minimize{\vb=(v_1,\ldots,v_M)\in (\HH)^M}{\iota_{D^\bot}(\vb)+\underbrace{\sum_{m=1}^M g_m^*(v_m)}_{\mbox{$g^*(\vb)$}}
}
\]
where $D^\bot = \menge{\vb=(v_1,\ldots,v_M)\in (\HH)^M}{v_1+\cdots+v_M = 0}$ is the orthogonal complement of $D$.
By making the variable change $(\forall m \in \{1,\ldots,M\})$ $v_m = u_m-u/M$ where $u$ is some given vector in $\HH$,
and by setting $h_m(u_m) = -g_m^*(u_m-u/M)$, the latter minimization can be reexpressed as
\[
\maximize{\substack{u_1\in \HH,\ldots,u_M\in \HH\\ u_1+\cdots+u_M = u}}{\sum_{m=1}^M h_m(u_m)}.
\] 
This problem is known as a \emph{sharing problem} where one wants to allocate a given resource $u$
between $M$ agents while maximizing the sum of their welfares
evaluated through their individual utility functions $(h_m)_{1\le m \le M}$.

\end{minipage}.
}
\end{center}

Another useful result follows from the fact that, by using the definition of the conjugate function of $g$,
Problem~\eqref{e:primal} can be reexpressed as the following saddle-point problem:
\begin{equation}\label{e:saddleprob}
\mbox{Find}\qquad \inf_{x\in \HH}\sup_{v\in \GG}\big(f(x) + \scal{Lx}{v} - g^*(v)\big).
\end{equation}
In order to find a saddle point $(\widehat{x},\widehat{v})\in \HH \times \GG$, it thus appears natural to impose
the inclusion relations:
\begin{equation}
-L^\top\widehat{v} \in \partial f(\widehat{x}),\qquad
L \widehat{x}  \in \partial g^*(\widehat{v}).
\end{equation}
A pair $(\widehat{x},\widehat{v})$ satisfying the above conditions is called a \emph{Kuhn-Tucker point}.
Actually, under some technical assumption, by using Fermat's rule and \eqref{e:subdifconj}, it can be proved  that, if 
$(\widehat{x},\widehat{v})$ is a Kuhn-Tucker point, then $\widehat{x}$ is a solution to the primal problem
and $\widehat{v}$ is a solution to the dual one. This property especially holds when $f\in \Gamma_0(\HH)$
and $g\in \Gamma_0(\GG)$.
\subsection{Duality in linear programming}\label{se:dualLP}
In linear programming (LP) \cite{Bertsimas_D_1997_book_introduction_lo}, we are interested in convex optimization problems of the form:
\begin{equation}\label{e:primalLP}
\text{Primal-LP}:\qquad \minimize{x\in [0,+\infty[^N}{c^\top x}\quad\text{s.t.}\quad Lx \ge b,
\end{equation}
where $L = (L^{(i,j)})_{1\le i \le K,1\le j \le N}\in \RR^{K\times N}$, $b\in \RR^K$, and $c\in \RR^N$.\footnote{The vector inequality 
in \eqref{e:primalLP} means that $Lx-b \in \RP^K$.}
The above formulation can be viewed as a special case of \eqref{e:primal} where
\begin{equation}
(\forall x \in \HH) \quad f(x) = c^\top x + \iota_{[0,+\infty[^N}(x),\qquad
(\forall z \in \GG) \quad g(z) = \iota_{[0,+\infty[^K}(z-b).
\end{equation}
By using
the properties of the conjugate function and by setting $y = -v$, it is readily shown that
the dual problem \eqref{e:dual} can be reexpressed as
\begin{equation}\label{e:dualLP}
\text{Dual-LP}:\qquad \maximize{y\in [0,+\infty[^K}{b^\top y}\quad\text{s.t.}\quad L^\top y \le c.
\end{equation}
Since $f$ is a convex function, strong duality holds in LP. If $\widehat{x} = (\widehat{x}^{(j)})_{1\le j \le N}$ is a solution to
Primal-LP, a solution $\widehat{y} = (\widehat{y}^{(i)})_{1\le i \le K}$ to Dual-LP can be obtained by the \emph{primal complementary slackness condition}:
\begin{equation}\label{e:slackprimal}
(\forall j\in \{1,\ldots,N\}) \quad\text{such that}\quad \widehat{x}^{(j)} > 0, \qquad
\sum_{i=1}^K L^{(i,j)}\, \widehat{y}^{(i)} = c^{(j)}.
\end{equation}
whereas, if $\widehat{y}$ is a solution to Dual-LP, a solution $\widehat{x}$ to Primal-LP can be obtained
by the \emph{dual complementary slackness condition}:
\begin{equation}\label{e:slackdual}
(\forall i\in \{1,\ldots,K\}) \quad\text{such that}\quad \widehat{y}^{(i)} > 0, \qquad
\sum_{j=1}^N L^{(i,j)}\, \widehat{x}^{(j)} = b^{(i)}.
\end{equation}

%where $\scal{\cdot}{\cdot}$ denotes the standard Euclidean product of $\RR^N

\section{Convex optimization algorithms} \label{se:optconv}
In this section, we present several primal-dual splitting methods for solving convex optimization problems, starting from the
basic forms to the more sophisticated highly parallelized ones. 
\subsection{Problem}
A wide range of convex optimization problems can be formulated as follows:
\begin{equation}
\label{e:primalvar}
\minimize{x\in\HH}{f(x)+g(L x) +h(x)}. 
\end{equation}
where $f \in \Gamma_0(\HH)$, $g \in \Gamma_0(\GG)$, $L \in \RR^{K\times N}$,
and $h \in \Gamma_0(\HH)$ is a differentiable function having a Lipschitzian
gradient with a Lipschitz constant $\beta \in \RPP$. The latter assumption means that the gradient $\nabla h$ of $h$ is such that
\begin{equation}
\big(\forall (x,y)\in (\HH)^2\big) \qquad \|\nabla h(x) -\nabla h(y) \| \le \beta \|x-y\|.
\end{equation}
For examples, the functions $f$, $g\circ L$, and $h$ may model various data fidelity terms and regularization functions
encountered in the solution of inverse problems. In particular, the Lipschitz differentiability property is satisfied for 
%quadratic functions.
least squares criteria.

With respect to Problem~\eqref{e:primal}, we have introduced an additional smooth term $h$.
This may be useful in offering more flexibility for taking into account the structure of the problem
of interest and the properties of the involved objective function. We will however see that not all algorithms
are able to possibly take advantage of the fact that $h$ is a smooth term.

Based on the results in Section~\ref{se:dualres} and Property~\eqref{e:sumdual} in Table \ref{t:propconj}, the dual optimization problem
reads:
\begin{equation}
\label{e:dualvar}
\minimize{v\in\GG}{(f^*\infconv h^*)(-L^\top v)+g(v)}. 
\end{equation}
Note that, in the particular case when $h=0$, the inf-convolution $f^*\infconv h^*$ (see the definition in Table~\ref{t:propconj}\eqref{e:definfconv}) of the conjugate functions
of $f$ and $h$ reduces to $f^*$ and we recover the basic form \eqref{e:dual} of the dual problem.

The common trick used in the algorithms which will be presented in this section is to solve jointly Problems~\eqref{e:primalvar}
and \eqref{e:dualvar}, instead of focusing exclusively on either \eqref{e:primalvar} or \eqref{e:dualvar}. More precisely, these algorithms
aim at finding a Kuhn-Tucker point $(\widehat{x},\widehat{v}) \in \HH\times \GG$ such that
\begin{equation}\label{e:KTPD}
-L^\top \widehat{v}-\nabla h(\widehat{x}) \in \partial f(\widehat x)\quad \text{and} \quad L \widehat{x}\in \partial g^*(\widehat{v}).
\end{equation} 

It has to be mentioned that some specific forms of Problem \eqref{e:primalvar} (e.g. when $g=0$) can be solved in a quite efficient manner by simpler proximal algorithms (see \cite{Combettes_P_2010_inbook_proximal_smsp}) than those described in the following.

\subsection{ADMM} \label{se:ADMM}
The celebrated ADMM (Alternating Direction Method of Multipliers) can be viewed as a primal-dual algorithm.
This algorithm belongs to the class of \emph{augmented Lagrangian} methods since a possible way of deriving this algorithm consists of looking
for a saddle point of an augmented version of the classical Lagrange function \cite{Boyd_S_2011_j-found-tml_distributed_osl_admm}. This augmented Lagrangian is defined as
\begin{equation}
\big(\forall (x,y,z)\in \HH\times (\GG)^2\big)\qquad \widetilde{\mathcal{L}}(x,y,z) = f(x)+h(x)+g(y)+\gamma \scal{(Lx-y)}{z}+\frac{\gamma}{2} \|Lx-y\|^2
\end{equation}
where $\gamma \in \RPP$ and $\gamma z$ corresponds to a Lagrange multiplier.
ADMM simply splits the step of minimizing the
augmented Lagrangian with respect to $(x,y)$ by alternating between the two variables,
while a gradient ascent is performed with respect to the variable $z$.
The resulting iterations are given in Algorithm~\ref{algo:admm}.

{\linespread{1.5}
\begin{algorithm}
{\small 
\caption{ADMM}\label{algo:admm}
\begin{equation*}
\begin{array}{l}
\text{Set}\;y_0\in \GG\;\text{and}\;z_0\in \GG\\
\text{Set}\;\gamma\in \RPP\\
\text{For}\;n=0,1,\ldots\\
\left\lfloor
\begin{array}{l}
x_n = \argmind{x\in \HH}{\frac{1}{2} \left\|L x - y_n + z_n\right\|^2+ \frac{1}{\gamma} \big(f(x)+h(x)\big)}\\
s_n = L x_n\\
y_{n+1} = \prox_{\frac{g}{\gamma}}\left(z_n+s_n\right)\\
z_{n+1} = z_n + s_n - y_{n+1}.
\end{array}
\right.
\end{array}
\end{equation*}
}
\end{algorithm}}

This algorithm has been known for a long time \cite{Gaba76,Fortin_M_1983_book_augmented_lmansbvp} although it has attracted recently much interest in 
the signal and image processing community (see e.g. \cite{Giovanelli_JF_2005-astoastro-positive-dse,Goldstein_T_2009_j-siam-is_split_bml,Figueiredo_M_2009_ssp_Deconvolution_opiuvsaalo,Figueiredo_M_2010_t-ip_restoration_piado,AfonsoBF11,Tran-Dinh_Q_2014_Primal_Dafccm}). 
A condition for the convergence of ADMM is as follows:
\begin{center}
\framebox{
\begin{minipage}{14cm}
\footnotesize\sffamily 
{CONVERGENCE OF ADMM}

Under the assumptions that
\begin{itemize}
\item $\rank(L) = N$,
\item Problem \eqref{e:primalvar} admits a solution,
\item $\inte(\dom g) \cap L(\dom f) \neq \emp$ or $\dom g \cap \inte\big(L(\dom f)\big)\neq \emp$,\footnotemark
\end{itemize}
%it can be shown that 
$(x_n)_{n\in\NN}$ converges to a solution to the primal problem \eqref{e:primalvar}
and $(\gamma z_n)_{n\in \NN}$ converges to a solution to the dual problem \eqref{e:dualvar}.
\end{minipage}
}
\footnotetext{More general qualification conditions involving the relative interiors of the domain of $g$ and $L(\dom f)$ can be obtained \cite{Combettes_P_2010_inbook_proximal_smsp}.}
\end{center}
A convergence rate analysis is conducted in \cite{Hong_M_2013_On_tlc}.

It must be emphasized that ADMM is equivalent to the application of the Douglas-Rachford algorithm \cite{Ecks92,Jsts07},
another famous algorithm in convex optimization, to the dual problem. Other primal-dual algorithms can be deduced from
the Douglas-Rachford iteration \cite{Bot_R_2013_siam-opt_Dou_rtp} or an augmented Lagrangian approach \cite{Chen_G_1994_j-mp_pro_bdm}.

Although ADMM was observed to have a good numerical performance in many problems, its applicability 
may be limited by the computation of $x_n$ at each iteration $n\in \NN$, which may be intricate due to the presence of matrix $L$, especially when this matrix is high-dimensional and has no simple structure. 
In addition, functions $f$ and $h$ are not dealt
with separately, and so the smoothness of $h$ is not exploited here in an explicit manner.

\subsection{Methods based on a Forward-Backward approach}\label{se:FBPD}
The methods which will be presented in this subsection are based on a forward-backward approach \cite{Combettes_P_2005_j-siam-mms_signal_rpfb}: they combine a gradient descent step (forward step)
with a computation step involving a proximity operator. The latter computation corresponds to a kind of subgradient step performed in an implicit (or backward) manner \cite{Combettes_P_2010_inbook_proximal_smsp}.
A deeper justification of this terminology is provided by the theory of monotone operators
\cite{Bauschke_H_2011_book_con_amo} which allows to highlight the fact that a pair $(\widehat{x},\widehat{v}) \in \HH\times \GG$ satisfying \eqref{e:KTPD} is a zero of a sum of two maximally monotone operators.
We will not go into details which can become rather technical, but we can mention that
the algorithms presented in this section can then be viewed as offsprings of the forward-backward algorithm for finding such a zero \cite{Bauschke_H_2011_book_con_amo}. Like ADMM, this algorithm is an
instantiation of a recursion converging to a fixed point of a nonexpansive mapping.

One of the most popular primal-dual method within this class is given by Algorithm~\ref{algo:fbpd}.
In the case when $h=0$, this algorithm can be viewed as an extension of the Arrow-Hurwitz method which performs
alternating subgradient steps with respect to the primal and dual variables in order to solve the saddle point problem \eqref{e:saddleprob} \cite{Nedic_A_2009_j-ota_subgradient_msp}.
Two step-sizes $\tau$ and $\sigma$ and relaxation factors $(\lambda_n)_{n\in \NN}$ are involved in Algorithm~\ref{algo:fbpd},
which can be adjusted by the user so as to get the best convergence profile for a given application.

{\linespread{1.5}
\begin{algorithm}
{\small
\caption{FB-based primal-dual algorithm}\label{algo:fbpd}
\begin{equation*}
\begin{array}{l}
\text{Set}\;x_0\in \HH\;\text{and}\;v_0\in \GG\\
\text{Set}\;(\tau,\sigma)\in \RPP^2\\
\text{For}\;n=0,1,\ldots\\
\left\lfloor
\begin{array}{l}
p_n = \prox_{\tau f}\big(x_n - \tau \big(\nabla h(x_n) + L^\top v_n\big)\big)\\
%q_n = \prox_{\sigma g^*}\big(v_n + \sigma \big(L (2 p_n - x_n)-\nabla \ell^*(v_n)\big)\big)\\
q_n = \prox_{\sigma g^*}\big(v_n + \sigma L (2 p_n - x_n)\big)\\
\text{Set}\;\lambda_n\in \RPP\\
(x_{n+1},v_{n+1}) = (x_n,v_n) + \lambda_n \big( (p_n,q_n) - (x_n,v_n)\big).
\end{array}
\right.
\end{array}
\end{equation*}}
\end{algorithm}

Note that when $L = 0$ and $g^*=0$ the basic form of the forward-backward algorithm (also called the proximal gradient algorithm) is recovered, a popular example of which is the iterative soft-thresholding algorithm \cite{Daub04}.

A rescaled variant of the primal-dual method (see Algorithm \ref{algo:fbpdresc}) is sometimes preferred, which can be deduced from the previous one by using Moreau's decomposition
\eqref{e:moreaudec} and by making the variable changes: $q_n' \equiv q_n/\sigma$ and $v'_n \equiv v_n/\sigma$. Under this form, it can be seen that, when $N=K$, $L=\Id$, $h = 0$,
and $\tau\sigma = 1$, the algorithm  reduces to the Douglas-Rachford algorithm (see \cite{Davis_D_2014_Convergence_rafdrss} for the link existing with extensions of the Douglas-Rachford algorithm).

{\linespread{1.5}
\begin{algorithm}
{\small
\caption{Rescaled variant of Algorithm \ref{algo:fbpd}}\label{algo:fbpdresc}
\begin{equation*}
\begin{array}{l}
\text{Set}\;x_0\in \HH\;\text{and}\;v'_0\in \GG\\
\text{Set}\;(\tau,\sigma)\in \RPP^2\\
\text{For}\;n=0,1,\ldots\\
\left\lfloor
\begin{array}{l}
p_n = \prox_{\tau f}\big(x_n - \tau \big(\nabla h(x_n) + \sigma L^\top v'_n\big)\big)\\
%q_n = \prox_{\sigma g^*}\big(v_n + \sigma \big(L (2 p_n - x_n)-\nabla \ell^*(v_n)\big)\big)\\
q'_n = (\Id-\prox_{g/\sigma})\big(v'_n + L (2 p_n - x_n)\big)\\
\text{Set}\;\lambda_n\in \RPP\\
(x_{n+1},v'_{n+1}) = (x_n,v'_n) + \lambda_n \big( (p_n,q'_n) - (x_n,v'_n)\big).
\end{array}
\right.
\end{array}
\end{equation*}}
\end{algorithm}}

Also, by using the symmetry existing between the primal and the dual problems, another variant of Algorithm~\ref{algo:fbpd} can be obtained
(see Algorithm \ref{algo:fbpdsym}) which is often encountered in the literature. When $L^\top L = \mu \Id$ with $\mu\in \RPP$, $h = 0$, $\tau\sigma \mu = 1$, and $\lambda_n \equiv 1$, 
Algorithm \ref{algo:fbpdsym} reduces to ADMM by setting $\gamma = \sigma$, and $z_n \equiv v_n/\sigma$ in Algorithm \ref{algo:admm}.

{\linespread{1.5}
\begin{algorithm}
{\small
\caption{Symmetric form of Algorithm \ref{algo:fbpd}}\label{algo:fbpdsym}
\begin{equation*}
\begin{array}{l}
\text{Set}\;x_0\in \HH\;\text{and}\;v_0\in \GG\\
\text{Set}\;(\tau,\sigma)\in \RPP^2\\
\text{For}\;n=0,1,\ldots\\
\left\lfloor
\begin{array}{l}
q_n = \prox_{\sigma g^*}\big(v_n + \sigma L x_n\big)\\
p_n = \prox_{\tau f}\big(x_n - \tau \big(\nabla h(x_n)+L^\top(2 q_n-v_n)\big)\big)\\
%q_n = \prox_{\sigma g^*}\big(v_n + \sigma \big(L (2 p_n - x_n)-\nabla \ell^*(v_n)\big)\big)\\
\text{Set}\;\lambda_n\in \RPP\\
(x_{n+1},v_{n+1}) = (x_n,v_n) + \lambda_n \big( (p_n,q_n) - (x_n,v_n)\big).
\end{array}
\right.
\end{array}
\end{equation*}}
\end{algorithm}}

Convergence guarantees were established in \cite{Condat_L_2013_j-ota-primal-dsm}, as well as for a more general version of this algorithm in \cite{Vu_B_2013_j-acm_spl_adm}:
\begin{center}
\framebox{
\begin{minipage}{14cm}
\footnotesize\sffamily 
{CONVERGENCE OF ALGORITHMS \ref{algo:fbpd} and \ref{algo:fbpdsym}}

Under the following sufficient conditions:
\begin{itemize}
\item $\tau^{-1} - \sigma \|L\|^2_{\rm S} \ge \beta/2$ where $\|L\|_{\rm S}$ is the spectral norm of $L$,
\item  $(\lambda_n)_{n\in \NN}$ a sequence in $]0,\delta[$ such that $\sum_{n\in \NN}\lambda_n (\delta - \lambda_n) = +\infty$
where $\delta = 2-\beta(\tau^{-1}-\sigma \|L\|^2_{\rm S})^{-1}/2 \in [1,2[$,
\item Problem \eqref{e:primalvar} admits a solution,
\item $\inte(\dom g) \cap L(\dom f) \neq \emp$ or $\dom g \cap \inte\big(L(\dom f)\big)\neq \emp$,
\end{itemize}
the sequences $(x_n)_{n\in\NN}$ and $(v_n)_{n\in \NN}$ 
%generated by  Algorithms \ref{algo:fbpd} and \ref{algo:fbpdsym} 
are such that the former one converges to a solution to the primal problem \eqref{e:primalvar}
and the latter one converges to a solution to the dual problem \eqref{e:dualvar}.
\end{minipage}
}
\end{center}
Algorithm \ref{algo:fbpd} also constitutes a generalization  of \cite{Chambolle_A_2010_first_opdacpai,Esser_E_2010_j-siam-is_gen_fcf,He_B_2012_j-siam-is_conv_apd} 
(designated by some authors as PDHG, Primal-Dual Hybrid Gradient). 
Preconditioned or adaptive versions of this algorithm were proposed in \cite{Pock_T_2008_p-iccv_diagonal_pffo,Combettes_P_2014_j-optim_Variable_mfb,Goldstein_T_2013_adaptive_pdh,Combettes_P_2014_p-icip_forward_bvo}
which may accelerate its convergence. Convergence rate results were also recently derived in \cite{Liang_J_2014_Convergence_rwi}.

Another primal-dual method (see Algorithm \ref{algo:fbpd2}) was proposed in \cite{Loris_I_2011_generalization_ist,Chen_P_2013_j-inv-prob_prim_dfp} which also results from a forward-backward approach \cite{Combettes_P_2014_p-icip_forward_bvo}. This algorithm is restricted to the case when $f = 0$ in Problem \eqref{e:primalvar}.

{\linespread{1.5}
\begin{algorithm}
{\small
\caption{Second FB-based primal-dual algorithm}\label{algo:fbpd2}
\begin{equation*}
\begin{array}{l}
\text{Set}\;x_0\in \HH\;\text{and}\;v_0\in \GG\\
\text{Set}\;(\tau,\sigma)\in \RPP^2\\
\text{For}\;n=0,1,\ldots\\
\left\lfloor
\begin{array}{l}
s_n=x_n - \tau \nabla h(x_n)\\
y_n = s_n - \tau L^\top v_n\\
q_{n}=\prox_{\sigma g^*}\Big(v_n+\sigma L y_n\Big)\\
p_n = s_n-\tau L^\top q_n\\
\text{Set}\;\lambda_n\in \RPP\\
(x_{n+1},v_{n+1}) = (x_n,v_n) + \lambda_n \big( (p_n,q_n) - (x_n,v_n)\big).
\end{array}
\right.\\[2mm]
\end{array}
\end{equation*}}
\end{algorithm}}
As shown by the next convergence result, the conditions on the step-sizes $\tau$ and $\sigma$ are less restrictive than for Algorithm \ref{algo:fbpd}.
\begin{center}
\framebox{
\begin{minipage}{14cm}
\footnotesize\sffamily 
{CONVERGENCE OF ALGORITHM \ref{algo:fbpd2}}

Under the assumptions that
\begin{itemize}
\item $\tau \sigma \|L\|^2_{\rm S} < 1$  and $\tau < 2/\beta$,
\item  $(\lambda_n)_{n\in \NN}$ a sequence in $]0,1]$ such that $\inf_{n\in \NN} \lambda_n > 0$,
\item Problem \eqref{e:primalvar} admits a solution,
\item $\inte(\dom g) \cap \ran(L) \neq \emp$,
\end{itemize}
the sequence $(x_n)_{n\in\NN}$ converges to a solution to the primal problem \eqref{e:primalvar}
(where $f=0$) and $(v_n)_{n\in \NN}$ converges to a solution to the dual problem \eqref{e:dualvar}.
\end{minipage}
}
\end{center}

Note also that the dual forward-backward approach that was proposed in \cite{Combettes_P_2010_j-svva_dualization_srp} for solving \eqref{e:primalvar} in the specific
case when $h = \|\cdot - r\|^2/2$ with $r\in \HH$
belongs to the class of primal-dual forward-backward approaches.

It must be emphasized that Algorithms \ref{algo:fbpd}-\ref{algo:fbpd2} present two interesting features which are very useful in practice.
At first, they allow to deal with the functions involved in the optimization problem at hand either through their proximity operator
or through their gradient. Indeed, for some functions, especially non differentiable or non finite ones, the proximity operator can be a very powerful tool \cite{Chaux_C_2007_j-ip_variational_ffbip} but, for some smooth functions
(e.g. the Poisson-Gauss neg-log-likelihood \cite{Jezierska_A_2012_p-icassp_primal_dps}) the gradient may be easier to handle. Secondly, these algorithms do not require to invert any matrix, but only 
to apply $L$ and its adjoint. This advantage is of main interest when large-size problems have to be solved for
which the inverse of $L$ (or $L^\top L$) does not exist or it has a no tractable expression.

\subsection{Methods based on a Forward-Backward-Forward approach}\label{se:FBFPF}
Primal-dual methods based on a forward-backward-forward approach were among the first primal-dual proximal methods proposed in the optimization literature, inspired from the seminal work in \cite{Tseng_P_2000_j-siam-control-optim_Modified_fbs}.
They were first developed in the case when $h=0$ \cite{Briceno_L_2011_j-siam-opt_mon_ssm}, then extended to more general scenarios in 
\cite{Combettes_P_2012_j-svva_pri_dsa} 
(see also \cite{Combettes_P_2013_siam-opt_sys_smi,Bot_R_2014_jmiv_conv_primal_apd} for further refinements).

{\linespread{1.5}
\begin{algorithm}
{\small
\caption{FBF-based primal-dual algorithm}\label{algo:fbfpd}
\begin{equation*}
\begin{array}{l}
\text{Set}\;x_0\in \HH\;\text{and}\;v_0\in \GG\\
\text{For}\;n=0,1,\ldots\\
\left\lfloor
\begin{array}{l}
\text{Set}\;\gamma_n\in \RPP\\
y_{1,n} = x_n - \gamma_n \big(\nabla h(x_n) + L^\top v_n\big)\\
y_{2,n} = v_n+\gamma_n L x_n\\
p_{1,n} = \prox_{\gamma_n f} y_{1,n}\\ 
p_{2,n} = \prox_{\gamma_n g^*} y_{2,n}\\
q_{1,n} = p_{1,n} - \gamma_n \big(\nabla h(p_{1,n}) + L^\top p_{2,n}\big)\\
q_{2,n} = p_{2,n} +\gamma_n L p_{1,n}\\
(x_{n+1},v_{n+1}) = (x_n - y_{1,n}+ q_{1,n},v_n - y_{2,n}+q_{2,n}).
\end{array}
\right.\\[2mm]
\end{array}
\end{equation*}}
\end{algorithm}}

The convergence of the algorithm is guaranteed by the following result: 
\begin{center}
\framebox{
\begin{minipage}{14cm}
\footnotesize\sffamily 
{CONVERGENCE OF ALGORITHM \ref{algo:fbfpd}}

Under the following assumptions:
\begin{itemize}
\item $(\gamma_n)_{n\in \NN}$ is a sequence in $[\varepsilon,(1-\epsilon)/\mu]$ where $\varepsilon \in ]0,1/(1+\mu)[$
and $\mu = \beta+\|L\|_{\rm S}$,
\item Problem \eqref{e:primalvar} admits a solution,
\item $\inte(\dom g) \cap L(\dom f) \neq \emp$ or $\dom g \cap \inte\big(L(\dom f)\big)\neq \emp$,
\end{itemize}
the sequence $(x_n,v_n)_{n\in \NN}$ converges to to a pair of primal-dual solutions. 
\end{minipage}
}
\end{center}

Algorithm \ref{algo:fbfpd} is often refered to as the M+LFBF (Monotone+Lipschitz Forward Backward Forward) algorithm. 
It enjoys the same advantages as FB-based primal-dual algorithms we have seen before. It however makes it possible to compute the proximity
operators of scaled versions of functions $f$ and $g^*$ in parallel. In addition, the choice of its parameters in order to satisfy convergence conditions
may appear more intuitive than for Algorithms \ref{algo:fbpd}-\ref{algo:fbpdsym}.
With respect to FB-based algorithms, an extra forward step however needs to be performed. This may lead to a slower convergence if, for example, the computational
cost of the gradient is high and an iteration of a FB-based algorithm is at least as efficient as an iteration of Algorithm \ref{algo:fbfpd}.

\subsection{A projection-based primal-dual algorithm} \label{se:projPD}
Another primal-dual algorithm was recently proposed in \cite{Alotaibi_A_2014_j-siam-opt_Solving_ccm} which relies on iterative projections onto half-spaces including the set of Kuhn-Tucker points
(see Algorithm \ref{algo:projpd}). 

{\linespread{1.5}
\begin{algorithm}
{\small
\caption{Projection-based primal-dual algorithm}\label{algo:projpd}
\begin{equation*}
\begin{array}{l}
\text{Set}\;x_0\in \HH\;\text{and}\;v_0\in \GG\\
\text{For}\;n=0,1,\ldots\\
\left\lfloor
\begin{array}{l}
\text{Set}\;(\gamma_n,\mu_n)\in \RPP\\
a_n = \prox_{\gamma_n (f+h)}(x_n-\gamma_n L^\top v_n)\\
l_n = L x_n\\
b_n = \prox_{\mu_n g}(l_n + \mu_n v_n)\\
s_n = \gamma_n^{-1} (x_n - a_n) + \mu_n^{-1} L^\top (l_n - b_n)\\
t_n = b_n - L a_n\\
\tau_n = \|s_n\|^2+\|t_n\|^2\\
\text{if}\;\tau_n = 0\\
\left\lfloor
\begin{array}{l}
\widehat{x} = a_n\\
\widehat{v} = v_n + \mu_n^{-1} (l_n-b_n)\\
\text{return}
\end{array}
\right.\\[2mm]
\text{else}\\
\left\lfloor
\begin{array}{l}
\text{Set}\;\lambda_n\in \RPP\\
\theta_n = \lambda_n (\gamma_n^{-1} \|x_n - a_n\|^2 + \mu_n^{-1} \|l_n - b_n\|^2)/\tau_n\\
x_{n+1} = x_n-\theta_n s_n\\
v_{n+1} = v_n -\theta_n t_n.
\end{array}
\right.
\end{array}
\right.\\[2mm]
\end{array}
\end{equation*}}
\end{algorithm}}

We have then the following convergence result:
\begin{center}
\framebox{
\begin{minipage}{14cm}
\footnotesize\sffamily 
{CONVERGENCE OF ALGORITHM \ref{algo:projpd}}

Assume that
\begin{itemize}
\item $(\gamma_n)_{n\in \NN}$ and $(\mu_n)_{n\in \NN}$ are sequences such that $\inf_{n\in \NN} \gamma_n > 0$, $\sup_{n\in \NN} \gamma_n < \pinf$,
$\inf_{n\in \NN} \mu_n > 0$, $\sup_{n\in \NN} \mu_n < \pinf$,
\item  $(\lambda_n)_{n\in \NN}$ a sequence in $\RR$ such that $\inf_{n\in \NN} \lambda_n > 0$ and $\sup_{n\in \NN} \lambda_n < 2$,
\item Problem \eqref{e:primalvar} admits a solution,
\item $\inte(\dom g) \cap L(\dom f) \neq \emp$ or $\dom g \cap \inte\big(L(\dom f)\big)\neq \emp$,
\end{itemize}
then, either the algorithm terminates in a finite number of iterations at a pair of primal-dual solutions $(\widehat{x},\widehat{v})$, or it generates
a sequence $(x_n,v_n)_{n\in \NN}$ converging to such a point. \end{minipage}
}
\end{center}

Although few numerical experiments have been performed with this algorithm, 
one of its potential advantages is that it introduces few constraints
on the choice of the parameters $\gamma_n$, $\mu_n$ and $\lambda_n$ at iteration $n$ and that it does not require
any knowledge on the norm of the matrix $L$. Nonetheless, the use of this algorithm does not allow us to exploit the fact
that $h$ is a differentiable function.

\subsection{Extensions} \label{se:extsplit}
More generally, one may be interested in more challenging convex optimization problems of the form:
\begin{equation}
\label{e:primalvare}
\minimize{x\in\HH}{f(x)+\sum_{m=1}^M\,(g_m\infconv\ell_m)(L_m x)+h(x)}, 
\end{equation}
where $f\in \Gamma_0(\HH)$, $h\in \Gamma_0(\HH)$, and, for every $m\in \{1,\ldots,M\}$, $g_m \in \Gamma_0(\RR^{K_m})$, 
$\ell_m \in \Gamma_0(\RR^{K_m})$, and $L_m \in \RR^{K_m \times N}$. 
The dual problem then reads
\begin{align}
\label{e:dualvare}
\minimize{v_1\in\RR^{K_1},\ldots,v_m\in\RR^{K_M}}{}
&\big(f^*\infconv h^*\big)\bigg(-\sum_{m=1}^M L_m^\top v_m\bigg)+\sum_{m=1}^M \big(g_m^*(v_m)+\ell_m^*(v_m)\big).
\end{align}
Some comments can be made on this general formulation. At first, one of its benefits is to split an original objective function in
a sum of a number of simpler terms. Such splitting strategy is often the key of an efficient resolution of difficult optimization problems.
For example, the proximity operator of the global objective function may be quite involved, while the proximity operators of the individual
functions may have an explicit form. A second point is that we have now introduced in the formulation, additional functions $(\ell_m)_{1\le m \le M}$.
These functions may be useful in some models \cite{Bec_S_2014_j-nonlinear-conv-anal_alg_sps}, but they present also the conceptual advantage to make the primal problem and its dual form
quite symmetric. For instance, this fact accounts for the symmetric roles played by Algorithms \ref{algo:fbpd} and \ref{algo:fbpdsym}. 
An assumption which is commonly adopted is to assume that, whereas $h$ is Lipschitz differentiable, the functions $(\ell_m)_{1\le m \le M}$ are strongly convex,
i.e. their conjugates are Lipschitz differentiable.
A last point to be emphasized is that, such split
forms are amenable to efficient parallel implementations. 
Using parallelized versions of primal-dual algorithms on multi-core architectures may 
render these methods even more successful for dealing with large-scale problems.

\begin{center}
\framebox{
\begin{minipage}{14cm}
\footnotesize\sffamily 
{HOW TO PARALLELIZE PRIMAL-DUAL METHODS ?}

Two main ideas can be used in order to put a primal-dual method under a parallel form.
%The main idea for putting a primal-dual method under a parallel form consists of embedding
%the optimization problem to be solved in a higher-dimensional space.

Let us first consider the following simplified form of  Problem \eqref{e:primalvare}:
\begin{equation}
\minimize{x\in\HH}{\sum_{m=1}^M\,g_m(L_m x)}. 
\end{equation}
A possibility consists of reformulating this problem in a higher-dimensional space as
\begin{equation}
\minimize{y_1\in \RR^{K_1},\ldots,y_M \in \RR^{K_M} }{f(\boldsymbol{y})+\sum_{m=1}^M\,g_m(y_m)},
\end{equation}
where $\boldsymbol{y} = [y_1^\top,\ldots,y_M^\top]^\top\in \RR^K$ with $K=K_1+\cdots+K_M$, and
$f$ is the indicator function of $\ran(\boldsymbol{L})$,
where $\boldsymbol{L} = [L_1^\top,\ldots,L_M^\top]^\top \in \RR^{K\times N}$.
Function $f$ serves to enforce the constraint: $(\forall m \in \{1,\ldots,M\})$ $y_m = L_m x$. By defining
the separable function $g\colon \boldsymbol{y} \mapsto \sum_{m=1}^M\,g_m(y_m)$, we are thus led
to the minimization of $f+g$ in the space $\RR^{K}$.
This optimization can be performed by the algorithms described in Sections \ref{se:ADMM}-\ref{se:projPD}. The proximity operator
of $f$ reduces to the linear projection onto $\ran(\boldsymbol{L})$, whereas the separability of $g$ ensures that its proximity operator
can be obtained by computing in parallel the proximity operators of the function $(g_m)_{1\le m \le M}$. Note that,
when $L_1 = \ldots = L_M = \Id$, we recover a consensus-based approach that we have already discussed. This technique can be used
to derive parallel forms of the Douglas-Rachford algorithm, namely the Parallel ProXimal Algorithm (PPXA) \cite{Combettes_PL_2008_j-ip_proximal_apdmfscvip} 
and PPXA+ \cite{Pesquet_J_2012_j-pjpjoo_par_ipo}, as well as 
parallel versions of ADMM (Simultaneous Direction Method of Multipliers or SDMM) \cite{Setzer_S_2009_j-jvcir_deblurring_pibsbt}.

The second approach is even more direct since it requires no projection onto $\ran(\boldsymbol{L})$. 
For simplicity, let us consider the following instance of Problem \eqref{e:primalvare}:
\begin{equation}
\minimize{x\in\HH}{f(x)+\sum_{m=1}^M\,g_m(L_m x)+h(x)}.
\end{equation}
By defining the function $g$ and the matrix $\boldsymbol{L}$ as in the previous approach, the problem can be recast as
\begin{equation}
\minimize{x\in\HH}{f(x)+g(\boldsymbol{L} x)+h(x)}.
\end{equation}
Once again, under appropriate assumptions on the involved functions, this formulation allows us to employ
the algorithms proposed in Sections \ref{se:FBPD}-\ref{se:projPD} and we still have the ability to compute the proximity operator of $g$ in a parallel manner.

\end{minipage}.
}
\end{center}

\section{Discrete optimization algorithms} \label{se:optdisc}

\subsection{Background on discrete optimization}
As already mentioned in the introduction, another common class of problems   in  signal processing and image analysis  are discrete optimization problems, for which primal-dual algorithms 
also  play an important role. Problems of this type are often  stated as \emph{integer linear programs} (ILPs), which can be expressed under the following form:
\begin{center}\vspace{3pt}\begin{tabular}{@{\ }l@{\ \ }
l@{\hspace{0.5cm}} l@{\ \ } l}
  $\mspace{-13mu}${Primal-ILP}: & $\minimize{x\in \HH}{c^\top x}$                 \\
          & s.t. $\ \,Lx\ge b$, $\quad x\in \mathcal{N}\subset\NN^N$,\\
%& $\mspace{37mu} x\in \mathcal{N}\subset\NN^N$,
\end{tabular}\end{center}
where $L = (L^{(i,j)})_{1\le i \le K,1\le j \le N}$ represents a 
 matrix of size $K\times N$, and $b = (b^{(i)})_{1\le i \le K}$, $c=(c^{(j)})_{1\le j \le N}$ are column vectors of size $K$ and $N$, respectively. Note that integer linear programming provides a very general formulation suitable for modeling a very broad range of problems, and will thus  form  the setting that we will consider hereafter.
Among the problems encountered in practice, many of them lead to a Primal-ILP that is NP-hard to solve. In such cases, a principled approach
for finding an approximate solution  is through the use of convex \emph{relaxations} (see framebox), where the original NP-hard problem is  approximated
with a surrogate one  (the so-called relaxed problem), which is convex and thus much easier to solve. The premise is the following: to the
extent that the surrogate problem provides a reasonably good approximation to the original optimization task,  one can expect to obtain an approximately optimal
solution for the latter   by essentially making use of or solving the former. 

\begin{center}
\framebox{
\begin{minipage}{14cm}
\footnotesize\sffamily 
{RELAXATIONS AND DISCRETE OPTIMIZATION}

Relaxations are  very useful for solving approximately discrete optimization problems.
Formally, given   %a minimization 
a problem 
\[
(\mathcal{P}): \minimize{x\in C}{f(x)}
\]
where $C$ is a subset of $\HH$, we say that  
\[
(\mathcal{P}'): \minimize{x\in C'}{f'(x)}
\]
with $C'\subset \HH$ is
a relaxation of $(\mathcal{P})$ if and only if (\textit{i}) $C\subset C'$, and (\textit{ii}) 
$(\forall x \in C')$ $f(x)\ge f'(x)$.

For instance, let us consider the integer linear program 
defined by $(\forall x \in \HH)$ $f(x) = c^\top x$ and
$C = S\cap \ZZ^N$, where $c\in \HH\setminus\{0\}$ and $S$ is a nonempty closed polyhedron defined as
\[
S = \menge{x\in \HH}{Lx \ge b}
\]
 with $L\in\RR^{K\times N}$ and $b\in \GG$. One possible linear programming relaxation 
of $(\mathcal{P})$ is obtained by setting $f'=f$ and $C'= S$,
which is  typically   much easier than $(\mathcal{P})$  (which is generally NP-hard). 
The quality of  $(\mathcal{P}')$ is quantified by its so-called integrality gap defined as 
$\frac{\inf f(C)}{\inf f'(C')} \ge 1$ (provided that $-\infty < \inf f'(C')\neq 0$).
%$\frac{\mathrm{opt(P)}}{\mathrm{opt(P')}}$, where $\mathrm{opt}(\cdot)$ is a problem's optimum.

Hence, for approximation purposes,   LP relaxations are  not all of equal value. If  
\[
(\mathcal{P}''): \minimize{x\in C''}{c^\top x}
\]
is another relaxation of $(\mathcal{P})$ with $C'' \subset C'$, then relaxation  $(\mathcal{P}'')$   
is tighter.
Interestingly, $(\mathcal{P})$ always has a tight  LP relaxation (with integrality gap  1)  given by
$C''=  \operatorname{conv}(S\cap\ZZ^N)$, where $\operatorname{conv}(C)$ is the convex hull polyhedron of $C$. 
Note, however, that if $(\mathcal{P})$ is NP-hard, polyhedron $\operatorname{conv}(S\cap\ZZ^N)$ will involve exponentially many inequalities.

The relaxations in all of the previous examples involve expanding the original feasible set. But, as mentioned, we can also derive relaxations by modifying the original objective function. For instance, in  so-called submodular relaxations \cite{Kolmogorov10,kahl-strandmark-dam-12}, one uses as  new objective a maximum submodular function that lower bounds the original objective.
{More generally, convex relaxations allow us to make use of the  well-developed duality theory of convex programming for dealing with discrete nonconvex problems.}\end{minipage}
}
\end{center}

The type of relaxations that are typically preferred  in large scale discrete optimization are  based on linear programming, involving the minimization of a linear function subject to linear inequality constraints. These can be naturally obtained by simply relaxing the integrality constraints of { Primal-ILP, thus leading to the relaxed primal problem  \eqref{e:primalLP} as well as its dual \eqref{e:dualLP}.}
It should be noted  that the use of LP-relaxations is often dictated by the need
of maintaining a reasonable computational cost. Although more powerful convex relaxations 
do exist in many cases, these may become intractable as the
number of variables grows larger, especially for  Semidefinite Programming (SDP) or Second-Order Cone Programming (SOCP) relaxations.

Based on the above observations, in the following we aim to present some  very general
primal-dual optimization strategies that can be used in this context, focusing a lot on  their underlying principles, which are based on two powerful techniques, the so-called \emph{primal-dual schema} and \emph{dual decomposition}. As we shall see, in order to estimate an approximate solution to Primal-ILP, both approaches make heavy use of the dual of the underlying
LP relaxation, i.e., Problem \eqref{e:dualLP}. But their strategies for doing so is quite different: the second  one essentially aims at solving this dual LP (and then converting the fractional solution into an integral one, trying not to increase the cost too much in the process), whereas the {first} one simply uses it in the design of the algorithm.

\subsection{The primal-dual schema {for integer linear programming}}\label{sec:primal_dual}
The primal-dual schema is a  well-known technique  in the combinatorial optimization community that has its origins in LP duality theory.
It is worth noting that it started as an exact method for solving linear programs.  As such, it had initially been used in  deriving exact polynomial-time algorithms for many  cornerstone problems in combinatorial optimization that   have a tight LP relaxation.  Its first use probably goes  back to Edmond's famous Blossom algorithm for constructing maximum matchings on graphs, but it had been also applied to many other combinatorial problems including max-flow (e.g., Ford and Fulkerson's  augmenting path-based techniques for max-flow can  essentially be understood in terms  of this schema), shortest path, minimum branching, and minimum spanning tree \cite{PapadimitriouS82}. In all of these cases, the primal-dual schema is driven by the fact that optimal LP solutions should satisfy 
%some conditions, called 
the \emph{complementary slackness conditions} (see \eqref{e:slackprimal} and \eqref{e:slackdual}). 
% JCP : THESE CONDITIONS WERE NOT INTRODUCED BEFORE
Starting with an initial  primal-dual pair of feasible solutions, it  therefore iteratively steers them towards satisfying these complementary slackness conditions (by trying at each step to minimize their total violation). Once this is achieved, both solutions (the primal and the dual) are guaranteed to be optimal. Moreover, since the primal is always chosen to be updated integrally during the iterations, it is ensured  that an integral optimal solution is obtained at the end. A notable feature of the primal-dual method %also turns out  to be the fact that  
is that it often reduces the original LP, which is a weighted optimization problem, to a series of purely combinatorial unweighted ones (related to minimizing  the violation of complementary slackness conditions  at each step).

Interestingly,  today the primal-dual schema is no longer used for providing exact algorithms. Instead, its main use concerns deriving approximation algorithms to NP-hard discrete problems that admit an ILP formulation, for which it has proved to be a very powerful and widely applicable tool. As such, it has been applied  to  many   NP-hard combinatorial problems up to now, including set-cover, Steiner-network, scheduling, Steiner tree, feedback vertex set,  facility location, to mention only a few \cite{Vazirani_book, Hochbaum_book}. With regard to problems from the domains of computer vision and image analysis, the primal-dual schema has been introduced recently in \cite{komodakis_TPAMI07,komodakis_CVIU08}, and  has been used for modeling a broad class of  tasks from these fields. %related to the discrete  optimization of Markov Random Fields.

It should be noted that for NP-hard ILPs    an integral solution  is no longer guaranteed to satisfy the complementary slackness conditions (since the LP-relaxation is not exact). How could it  then be possible to apply this schema to such problems? It turns out that the  answer to  this question consists of using an appropriate relaxation of the above conditions.
 To  understand exactly how   we need to proceed in this case, let us consider the problem Primal-ILP above.
As already explained, the goal is to compute an optimal solution to it, but, due to the integrality constraints $x \in \mathcal{N}$, this is assumed to be  NP-hard, 
and so we can only %hope to 
estimate an approximate solution. To achieve that, we will first need to relax the integrality constraints, thus giving  rise to the relaxed primal problem in \eqref{e:primalLP} as well as its dual \eqref{e:dualLP}.
A primal-dual 
algorithm attempts to compute an approximate solution to Primal-ILP by relying on the following principle (see framebox for an explanation):\\ 
\textbf{Primal-dual principle in the discrete case}: 
\emph{
Let $x\in \HH$ and $y\in \GG$ be integral-primal and dual feasible solutions (i.e.
$x \in \mathcal{N}$ and $Lx\ge b$, and $y\in [0,+\infty[^K$ and $L^\top y \le c$). Assume that there
exists $\nu \in [1,+\infty[$ such that
\begin{equation} \label{primad_dual_eq} c^\top x  \le \nu\, b^\top y.
\end{equation}
Then, $x$ {can be shown to be} a $\nu$-approximation to an unknown optimal integral solution $\widehat{x}$,
 i.e.
 \begin{equation}\label{eq:primal_dual_principle_ineq}
 c^\top \widehat{x}\!\leq c^\top x\leq \nu\, c^\top \widehat{x}.
 \end{equation}}

\begin{center}
\framebox{
\begin{minipage}{14cm}
\footnotesize\sffamily 
{PRIMAL-DUAL PRINCIPLE IN THE DISCRETE CASE}

Essentially, {the proof of }this principle relies on  the fact that the sequence of optimal costs of problems   Dual-LP, Primal-LP, and Primal-ILP is increasing.

\begin{center}
\vspace{5pt}
\includegraphics[width=0.7\linewidth]{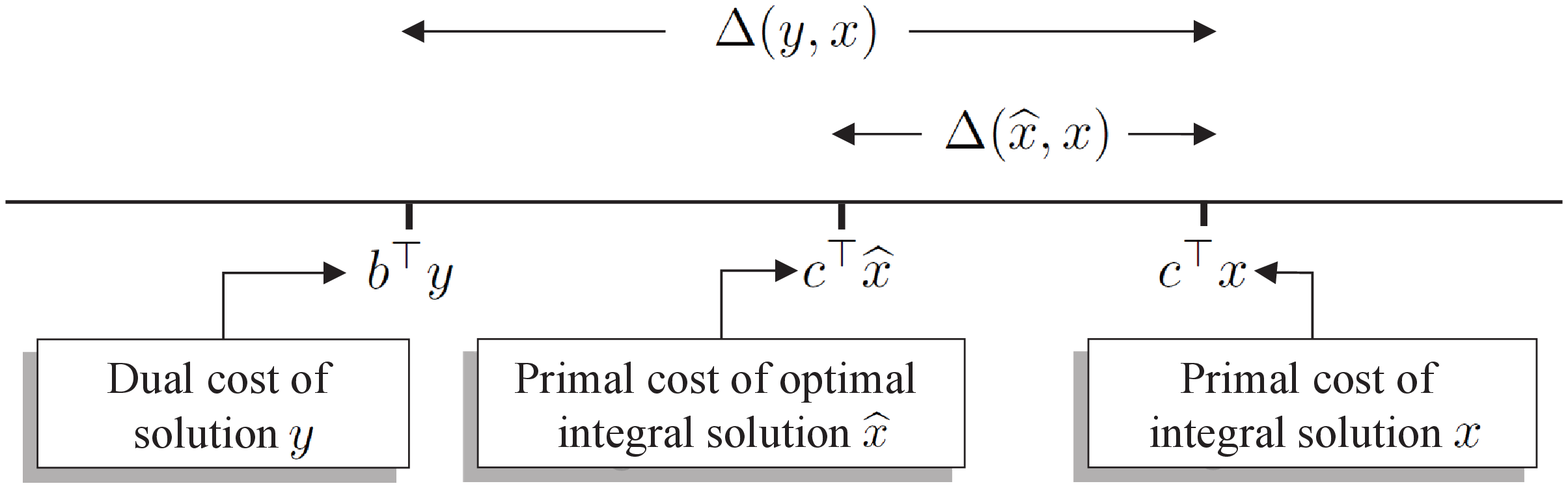}
\vspace{-17pt}
\end{center}

Specifically, by weak LP  duality, {the optimal cost of  Dual-LP  is known to not exceed the optimal cost of  Primal-LP.
% , i.e.,  if $x$ and $y$  satisfy respectively the constraints of Primal-LP and Dual-LP, then it holds $c^\top{x}\le b^\top y$. 
As a result of this fact,   the
 cost $c^\top\widehat{x}$ (of an unknown optimal integral solution $\widehat{x}$)  is guaranteed to  be at least as large as  the cost $b^\top y$ of any dual feasible solution $y$. On the other hand, by definition,  $c^\top\widehat{x}$ cannot exceed the cost $c^\top x$ of an integral-primal feasible solution $x$. } Therefore, if  the gap $\Delta(y,x)$ between the costs of $y$ and $x$ is small (e.g., it holds $c^\top x \le \nu\, b^\top y$), the same will be true for  the gap $\Delta(\widehat{x},x)$ between the costs of $\widehat{x}$ and $x$  (i.e., $c^\top x \le \nu\,c^\top\widehat{x}$),  thus proving that $x$ is a $\nu$-approximation to optimal solution  $\widehat{x}$.

\end{minipage}
}
\end{center}

% \begin{figure}[t]
% \center
% \renewcommand{\captionlabelfont}{\bf}
% \renewcommand{\captionfont}{\small}
%  \includegraphics[width=0.7\linewidth]{images/pd_principle.eps}
%  \label{fig:pd_schema1}
% \caption{By weak LP  duality,  the
%  cost $c^\top\widehat{x}$ (of an unknown optimal integral solution $\widehat{x}$) is guaranteed to lie between  costs $b^\top y$ and $c^\top x$ of any pair  of integral-primal and dual feasible solutions $(x,y)$. 
% Therefore, if  the gap $\Delta(y,x)$ between the costs of $y$ and $x$ is small (e.g., it holds $c^\top x \le \nu\, b^\top y$), the same will be true for  the gap $\Delta(\widehat{x},x)$ between the costs of $\widehat{x}$ and $x$  (i.e., $c^\top x \le \nu\,c^\top\widehat{x}$),  thus proving that $x$ is a $\nu$-approximation to optimal solution  $\widehat{x}$.
% }
% \end{figure}

Although the above principle lies at the heart of many
primal-dual techniques (i.e.,   in  one way or another, primal-dual methods often try to fulfill the assumptions imposed by this principle), it does not directly specify how to 
estimate a primal-dual pair of solutions $(x,y)$ that satisfies these assumptions. This is where the so-called \emph{relaxed complementary slackness conditions} come into play, 
%in handy, 
as they typically provide an alternative and more convenient (from an algorithmic viewpoint) way for generating such a pair of solutions.
These conditions generalize the complementary slackness conditions associated with an arbitrary pair of primal-dual linear programs
(see Section~\ref{se:dualLP}). The latter conditions    apply only in cases when there is no duality gap,  like between Primal-LP and Dual-LP, but they are not  applicable to cases like Primal-ILP and Dual-LP, when a duality gap exists as a result of the integrality constraint imposed on variable $x$.
 As in the exact case,  two  types of relaxed complementary slackness conditions exist, depending on whether the primal or dual variables are checked for being zero. 

\vspace{3pt}
\textbf{Relaxed Primal Complementary Slackness Conditions} with relaxation factor $\nu_{\mathrm{primal}}\le 1$.
For a given $x=(x^{(j)})_{1\le j\le N}\in \HH$, $y = (y^{(i)})_{1\le i \le K}\in \GG$, {the following conditions are assumed to hold:
}\begin{flalign}\label{eq:primal_csc}
\mspace{50mu}&(\forall j\in J_x)\qquad {\nu_{\mathrm{primal}}}\, c^{(j)}\le \sum_{i=1}^K L^{(i,j)}y^{(i)} \le c^{(j)}&
\end{flalign}
where $J_x = \menge{j\in\{1,\ldots,N\}}{x^{(j)}>0}$.
%The set of vectors $y$ satisfying this property will be denoted by $D_{x,\nu_{\mathrm{primal}}}$

\vspace{3pt}
\textbf{Relaxed Dual Complementary Slackness Conditions} with relaxation factor $\nu_{\mathrm{dual}}\ge 1$.
For a given $y = (y^{(i)})_{1\le i \le K}\in \GG$, $x=(x^{(j)})_{1\le j\le N}\in \HH$, {the following conditions are assumed to hold}:
\begin{flalign}\label{eq:dual_csc}
\mspace{50mu}&(\forall i\in I_y) \qquad b^{(i)}\le \sum_{j=1}^N L^{(i,j)}x^{(j)} \le \nu_{\mathrm{dual}}\, b^{(i)}&
\end{flalign}
where $I_y = \menge{i\in\{1,\ldots,K\}}{y^{(i)}>0}$.
%The set of vectors $x$ satisfying this property will be denoted by $C_{y,\nu_{\mathrm{dual}}}$.

When both $\nu_{\mathrm{primal}}=1$ and $\nu_{\mathrm{dual}}=1$, we recover the exact complementary slackness conditions in
\eqref{e:slackprimal} and \eqref{e:slackdual}. The use of the above conditions in the context of a primal-dual approximation algorithm becomes clear by the following result: 
 
%\begin{PrimalDual1}\label{thm:relaxed_comp_slack} 
\emph{
If $x = (x^{(j)})_{1\le j \le N}$ and $y = (y^{(i)})_{1\le i \le K}$ are feasible with respect to \emph{Primal-ILP} and \emph{Dual-LP} respectively, and  satisfy the  relaxed complementary slackness conditions \eqref{eq:primal_csc} and \eqref{eq:dual_csc},
then  the pair $(x,y)$
satisfies the  primal-dual principle in the discrete case with $\nu = \frac{\nu_{\mathrm{dual}}}{\nu_{\mathrm{primal}}}$. Therefore, $x$ is a $\nu$-approximate solution to \emph{Primal-ILP}.
%\end{PrimalDual1}
}

This result simply follows from the inequalities
%\begin{proof}
\begin{align}
c^\top x = \sum_{j=1}^N c^{(j)} x^{(j)}\overset{\eqref{eq:primal_csc}}{\le} \sum_{j=1}^N \Big(\frac{1}{\nu_{\mathrm{primal}}}\sum_{i=1}^K L^{(i,j)}y^{(i)}\Big)x^{(j)}&=\frac{1}{\nu_{\mathrm{primal}}}\sum_{i=1}^K\Big(\sum_{j=1}^N L^{(i,j)}x^{(j)}\Big)y^{(i)}\nonumber\\
&\overset{\eqref{eq:dual_csc}}{\le}\frac{\nu_{\mathrm{dual}}}{\nu_{\mathrm{primal}}}\sum_{i=1}^K  b^{(i)} y^{(i)} = \frac{\nu_{\mathrm{dual}}}{\nu_{\mathrm{primal}}}b^\top y.
\end{align}
%\qed
%\end{proof}

%\looseness=-1 
Based on  the above result, iterative schemes can be devised yielding  a primal-dual $\nu$-approximation algorithm. For example, we can employ the following algorithm:
%\begin{PrimalDual2}[Primal-dual schema]

{\linespread{1.5}
\begin{algorithm} 
{\small
\caption{Primal-dual schema} \label{e:primduaschema} 
\emph{Generate a sequence $(x_n,y_n)_{n\in \NN}$ of elements of $\RR^N\times \RR^K$
as follows:}
\begin{equation}
\begin{array}{l}
\text{Set}\;\nu_{\mathrm{primal}}\le 1\;\text{and}\;\nu_{\mathrm{dual}}\ge 1\\
\text{Set}\;y_0\in [0,+\infty[^K\;\text{such that}\;L^\top y_0 \le c\\
\text{For}\;n=0,1,\ldots\\
\left\lfloor
\begin{array}{l}
\text{Find}\;x_n\in \menge{x\in \mathcal{N}}{Lx\ge b}\;\text{minimizing}\\
\qquad \sum_{i\in I_{y_n}} q^{(i)}\;\;\text{s.t.}\;\; (\forall i \in I_{y_n})\quad \sum_{j=1}^N L^{(i,j)}x^{(j)} \le \nu_{\mathrm{dual}}\, b^{(i)}+q^{(i)},\; q^{(i)} \ge 0\\
\text{Find}\;y_{n+1}\in \menge{y\in[0,+\infty[^K}{L^\top y \le c}\;\text{minimizing}\\
\qquad \sum_{j\in J_{x_n}} r^{(j)}\;\;\text{s.t.}\;\; (\forall j \in J_{x_n})\quad 
\sum_{i=1}^K L^{(i,j)}y^{(i)}+r^{(j)}\ge {\nu_{\mathrm{primal}}}\, c^{(j)},\; r^{(j)} \ge 0.\\
\end{array}
\right.
\end{array}
\end{equation}
}
\end{algorithm}
}
%series of pairs of integral-primal and dual solutions
%$\{(\mathbf{x}^k,\mathbf{y}^k)\}_{k=1}^t$, such that the  elements $\mathbf{x}^t$, $\mathbf{y}^t$ of the last
%pair are both feasible and satisfy the relaxed primal complementary
%slackness conditions \eqref{eq:primal_csc} and \eqref{eq:dual_csc}.
%\end{PrimalDual2}

Note that, in this scheme,  primal solutions are always updated integrally. Also, note that, when applying the primal-dual schema, different implementation strategies are possible. 
The  strategy described in Algorithm \ref{e:primduaschema} is to maintain feasible primal-dual solutions $(x_n,y_n)$ at iteration $n$, and  iteratively improve how tightly the  
(primal or dual) complementary slackness conditions get satisfied.
This is performed through the introduction of slackness variables $(q^{(i)})_{i\in I_{y_n}}$ and $(r^{(j)})_{j\in J_{x_n}}$
the sums of which measure the degrees of violation of each relaxed slackness condition and have thus to be minimized.
 Alternatively, for example, we can opt to maintain solutions  $(x_n,y_n)$   that  satisfy the relaxed complementary slackness conditions, but may be infeasible, and 
 iteratively improve the feasibility of the generated solutions. For instance, if we  start with a  feasible dual solution but with an infeasible primal  solution, such a scheme 
 would result into improving  the feasibility of the primal solution,  as well as the optimality of  the dual solution at each iteration, ensuring that  a feasible primal solution 
 is obtained at the end. No matter which one of the above two strategies we choose to follow, the end result will be to gradually bring the
 primal and dual costs $c^\top x_n$ and $b^\top y_n$  closer and closer together  so that asymptotically the primal-dual principle gets satisfied with the desired approximation factor. Essentially, at each iteration,  through the coupling by the complementary slackness conditions the current primal solution is used to improve the dual, and vice versa. 
%(i.e.,  the current dual solution  is used to improve the primal).

Three remarks are  worth making at this point: the first one relates to the fact that the
two processes, i.e. the primal and the dual, make only local improvements
to each other. Yet, in the end they manage to yield a result that is almost  globally optimal. The second point to emphasize is that, for computing this approximately optimal result, the algorithm requires no solution to the Primal-LP or Dual-LP  to be computed, which are replaced by simpler optimization problems. This is an important advantage from a computational standpoint since, for large scale problems, solving these relaxations can often be quite costly.
In fact, in most cases where we apply the primal-dual schema, purely combinatorial algorithms can be obtained that contain no sign of linear programming in the end. A last point to be noticed is that these algorithms require appropriate choices of the relaxation factors $\nu_{\mathrm{primal}}$ and $\nu_{\mathrm{dual}}$, which are often application-guided.

\noindent\textbf{{Application to the set cover problem}}:
For a simple illustration   of the primal-dual schema, let us consider the problem of set-cover, which is known to be NP-hard. In this problem, we are given as input a 
finite set $\mathcal{V}$ of $K$ elements $(\upsilon^{(i)})_{1\le i \le K}$, a collection of (non disjoint) subsets $\mathcal{S}=\{S_j\}_{1\le j\le N}$ where, 
for every $j\in \{1,\ldots,N\}$, $S_j\subset \mathcal{V}$, and $\bigcup_{j=1}^N S_j = \mathcal{V}$. Let $\varphi\colon\mathcal{S}\rightarrow \mathbb{R}$ be a function that assigns a cost {$c_j= \varphi(S_j)$} for each subset $S_j$. The goal is to find a set cover (i.e. a subcollection of $\mathcal{S}$ that covers all elements of $\mathcal{V}$)  that has minimum cost {(see Fig.~\ref{fig:setcover}). }
\begin{figure}[htb]
\begin{center}
\includegraphics[height=2cm]{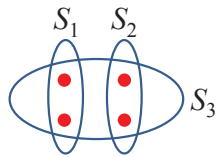}
\end{center}
\caption{\small {A toy set-cover instance with $K=4$ and $N=3$, where $\varphi(S_1)=\frac{1}{2}$, $\varphi(S_2)=1$, $\varphi(S_3)=2$. In this case, the optimal set-cover is  $\{S_1, S_2\}$ and has a cost of $\frac{3}{2}$.}}\label{fig:setcover}
\end{figure}

The above problem can  be expressed as the following ILP:
\begin{align}
&\minimize{x = (x^{(j)})_{1\le j \le N} }{ \sum_{j=1}^N \varphi(S_j)\,x^{(j)}}\ \label{eq:sc1}\\
&\quad\mathrm{s.t.} \;\;(\forall i \in \{1,\ldots,K\}) \quad \sum_{\substack{j\in \{1,\ldots,N\}\\ \upsilon^{(i)} \in S_j}} x^{(j)} \ge 1,  \quad\; x \in \{0,1\}^N,\label{eq:sc2}
%& \qquad\quad\; x \in \{0,1\}^N, \label{eq:sc3}
\end{align}
where   indicator variables $(x^{(j)})_{1\le j \le N}$ are used for determining if a set in $\mathcal{S}$ has been included in the set cover or not, and \eqref{eq:sc2}     %-\eqref{eq:sc3} 
 ensures that each one of the elements of $\mathcal{V}$ is contained in at least one of the sets that were chosen for participating to  the set cover.

An LP-relaxation for this problem is obtained by simply  replacing the Boolean constraint %\eqref{eq:sc3} 
with the constraint $x \in [0,+\infty[^N$. The  dual of this LP relaxation is given by the following linear program:
\begin{align}
&\maximize{y = (y^{(i)})_{1\le i \le K} \in [0,+\infty[^K }{ \sum_{i=1}^K y^{(i)}}\ \\
&\quad\mathrm{s.t.} \;\;(\forall j \in \{1,\ldots,N\}) \quad \sum_{\substack{i\in \{1,\ldots,K\}\\ \upsilon^{(i)} \in S_j } } y^{(i)} \le \varphi(S_j).  \label{eq:scd2}
\end{align}
%\begin{align}
%\max & \sum_{p\in \mathcal{V}}y_p & \mspace{-150mu}\label{eq:scd1}\\
%\mathrm{s.t.} \quad & \sum_{p\in S} y_p\le f(S), &\mspace{-250mu}\forall S\in \mathcal{S}  \label{eq:scd2}\\
%& y_p\ge 0, &\mspace{-250mu}\forall p\in \mathcal{V}  \label{eq:scd3}
%\end{align}

Let us  denote by $F_{\max}$ the maximum frequency of an element in $\mathcal{V}$, where by the term  \emph{frequency} %of an element 
we mean the number of sets this element belongs to. In this case, we will use the primal-dual schema to derive an $F_{\max}$-approximation algorithm by choosing $\nu_{\mathrm{primal}}=1$, $\nu_{\mathrm{dual}}=F_{\max}$.  This results into the  following  complementary slackness conditions, which we will need to satisfy:

\noindent{Primal Complementary Slackness Conditions}
\begin{flalign}\label{eq:primal_csc_sc}
%\mspace{50mu}&\forall S\in\mathcal{S}, \text{ if } x_S>0 \text{ then  } \sum_{p\in S} y_p = f(S)&
\mspace{50mu}&(\forall j\in \{1,\ldots,N\}) \text{ if } x^{(j)}>0 \text{ then  } \sum_{\substack{i\in \{1,\ldots,K\}\\ \upsilon^{(i)} \in S_j}} y^{(i)} = \varphi(S_j)&
\end{flalign}

\noindent{Relaxed Dual Complementary Slackness Conditions} (with relaxation factor $F_{\max}$)\begin{flalign}\label{eq:dual_csc_sc}
\mspace{50mu}&(\forall i\in\{1,\ldots,K\}) \text{ if } y^{(i)}>0 \text{ then } \sum_{\substack{j\in \{1,\ldots,N\}\\ \upsilon^{(i)} \in S_j } }x^{(j)}\le F_{\max}.&
%\mspace{50mu}&\forall p\in\mathcal{V}, \text{ if } y_p>0 \text{ then } \sum_{S:\, p\in S}x_S\le F_{\max}&
\end{flalign}
A set $S_j$ with $j\in \{1,\ldots,N\}$ for which $\sum_{\substack{i\in \{1,\ldots,K\}\\ \upsilon^{(i)} \in S_j}} y^{(i)} = \varphi(S_j)$ will be called \emph{packed}. Based on this definition, and given that the primal variables $(x^{(j)})_{1\le j \le N}$ are always  kept integral (i.e., either $0$ or $1$) during the primal-dual schema, Conditions \eqref{eq:primal_csc_sc} basically say that only packed sets can be included in the set cover (note that  overpacked sets are already forbidden by feasibility constraints \eqref{eq:scd2}).
Similarly,
Conditions \eqref{eq:dual_csc_sc}  require that an element $\upsilon^{(i)}$ with $i\in \{1,\ldots,K\}$ associated with a nonzero dual variable $y^{(i)}$ should not be covered more than $F_{\max}$ times, which is, of course, trivially satisfied given  that $F_{\max}$ represents the maximum frequency of any element in $\mathcal{V}$.

%
%
%\incmargin{1em}
%\linesnumbered 
%\setalcapskip{2ex}
%\begin{algorithm}
        %\SetLine
        %Set $\mathbf{x}\leftarrow \mathbf{0}$, $\mathbf{y}\leftarrow \mathbf{0}$\;\\
        %Declare all elements in $\mathcal{V}$ as uncovered\;\\
        %\While{$\mathcal{V}$  contains uncovered elements }
        %{
                %Select an uncovered element $p$ and increase $y_p$ until some\ set becomes packed\;% (i.e., \eqref{eq:primal_csc_sc} is satisfied)\;\\
                %Let $\mathcal{S}_0$ denote all sets from $\mathcal{S}$ that are packed\; %satisfy condition \eqref{eq:primal_csc_sc}\;\\
                %Set $x_S\leftarrow 1$, $\forall \in \mathcal{S}_0$   \tcc*{include all  sets that are packed in the cover}\\
                %Declare all elements belonging to at least one set $S$ with $x_S=1$ as covered\;
        %}
        %\caption{Primal-dual schema for set-cover}
        %\label{algo:sc}
%\end{algorithm}
%
%

{\linespread{1.5}
\begin{algorithm}
{\small
\caption{Primal-dual schema for set-cover.}\label{algo:sc}
\begin{equation*}
\begin{array}{l}
\text{Set}\;x_0 \leftarrow  0, y_0 \leftarrow  0\\
\text{Declare all elements in $\mathcal{V}$ as uncovered}\\
\text{While $\mathcal{V}$  contains uncovered elements}\\
\left\lfloor
\begin{array}{l}
\text{Select an uncovered element $\upsilon^{(i)}$ with $i\in \{1,\ldots,K\}$ and increase $y^{(i)}$ until some set becomes packed}\\% (i.e., \eqref{eq:primal_csc_sc} is satisfied)
\text{For every packed set $S_j$ with $j\in \{1,\ldots,N\}$, set $x^{(j)}\leftarrow 1$}\\
\qquad  \text{(include all the sets that are packed in the cover)}\\
\text{Declare all the elements belonging to at least one set $S_j$ with $x^{(j)} = 1$ as covered.}
\end{array}
\right.
\end{array}
\end{equation*}}
\end{algorithm}
}

Based on the above observations, the  iterative method whose pseudocode is shown in Algorithm~\ref{algo:sc} emerges naturally as a simple variant of Algorithm \ref{e:primduaschema}.
Upon its termination, both $x$ and $y$ will be feasible given that there will be no uncovered element and no set that violates \eqref{eq:scd2}. Furthermore, given that the final pair $(x,y)$  satisfies the relaxed complementary slackness conditions with $\nu_{\mathrm{primal}}=1$, $\nu_{\mathrm{dual}}=F_{\max}$, the set cover defined by  $x$ will provide an $F_{\max}$-approximate solution.

\subsection{Dual decomposition}\label{sec:dd_sec_}
We will next examine a different  approach for discrete optimization, which is based on the principle of dual decomposition \cite{komodakis_pami2011, Komodakis_iccv07, Bertsekas_D_2004_book_non_prog}. The core idea behind this principle   essentially  follows a divide and conquer strategy: that is, given a difficult or high-dimensional optimization problem,
we decompose it into smaller easy-to-handle subproblems  and then extract an overall solution by cleverly combining the solutions from these subproblems.

To explain this technique, 
we will consider the general problem of minimizing the energy of a discrete Markov Random Field (MRF), which is a ubiquitous  problem in the fields of computer vision and image analysis (applied with great success on a wide variety of tasks from these domains {such as stereo-matching, image segmentation, optical flow estimation, image restoration and inpainting, or object detection}) \cite{Blake_book}. 
This problem involves a graph $G$  with vertex set $\mathcal{V}$ and edge set $\mathcal{E}$ (i.e., $G=(\mathcal{V},\mathcal{E})$) plus a finite label set $\mathcal{L}$. The goal is to find a labeling 
$z= (z^{(p)})_{p \in \mathcal{V}} \in \mathcal{L}^{|\mathcal{V}|}$
for the graph vertices that has minimum cost, that is  
\begin{equation}\label{eq:mrf_energy_}
\minimize{z\in \mathcal{L}^{|\mathcal{V}|}}{\sum_{p\in\mathcal{V}}\varphi_p(z^{(p)})+\sum_{e\in \mathcal{E}} \thh_e(\mathsf{z}^{(e)}})
\end{equation}
where, for every $p\in \mathcal{V}$ and $e\in \mathcal{E}$, $\varphi_p\colon \mathcal{L}  \to \,\left]-\infty,+\infty\right[$ and $\thh_e\colon \mathcal{L}^2  \to \,\left]-\infty,+\infty\right[$
represent the  unary and pairwise costs
(also known connectively as MRF potentials $\varphi=\left\{\{\varphi_p\}_{p\in\mathcal{V}},\{\thh_e\}_{e\in \mathcal{E}}\right\}$), and $\mathsf{z}^{(e)}$ denotes the pair of components of
$z$ defined   by the variables 
corresponding to  vertices connected by $e$ (i.e.,  $\mathsf{z}^{(e)}=(z^{(p)}, z^{(q)})$ for $e=(p,q)\in \mathcal{E}$).

The above problem is NP-hard, and much of the recent work on MRF optimization  revolves around the following equivalent ILP formulation of \eqref{eq:mrf_energy_}
 \cite{chekuri}, which is the one that we  will also use here:
 \begin{equation}\label{eq:orig_MRF_}
 \minimize{x\in C_G}{f(x;\varphi)= \sum_{p\in \mathcal{V},\,z^{(p)}\in \mathcal{L}} \varphi_p(z^{(p)})\, x_p(z^{(p)})+
  \sum_{e\in \mathcal{E},\, \mathsf{z}^{(e)}\in \mathcal{L}^2}\thh_e(\mathsf{z}^{(e)})\, \mathsf{x}_e(\mathsf{z}^{(e)})},
 \end{equation}
 where the set $C_G$ is defined for any graph $G=(\mathcal{V},\mathcal{E})$ as 
\begin{equation}\label{eq:mrf_ilp_}
C_G =
\left\{ {
        x = \big\{\{x_p\}_{p\in\mathcal{V},z\in\mathcal{L}},\{\mathsf{x}_e\}_{e\in\mathcal{V},\mathsf{z}\in\mathcal{L}^2}\big\}
        }
\left|
\begin{array}{l}
   (\forall p\in \mathcal{V}) \qquad\qquad\qquad\qquad\quad\;  \sum_{z^{(p)}\in \mathcal{L}} x_p(z^{(p)})=1  \\
   (\forall e=(p,q)\in \mathcal{E})(\forall z^{(q)}\in \mathcal{L})\quad 
\sum_{\mathsf{z}^{(e)}\in \mathcal{L}\times \{z^{(q)}\}} \mathsf{x}_e(\mathsf{z}^{(e)})=x_q(z^{(q)})   \\
   (\forall e=(p,q)\in \mathcal{E})(\forall z^{(p)}\in \mathcal{L})\quad 
\sum_{\mathsf{z}^{(e)}\in \{z^{(p)}\}\times \mathcal{L}} \mathsf{x}_e(\mathsf{z}^{(e)})=x_p(z^{(p)})  \\
{(\forall p\in \mathcal{V})\qquad\qquad\qquad\qquad\qquad x_p(\cdot)\colon\mathcal{L}\mapsto \{0,1\}}\\
{(\forall e\in \mathcal{E})\qquad\qquad\qquad\qquad\qquad \mathsf{x}_e(\cdot)\colon \mathcal{L}^2\to \{0,1\}}
\end{array}\right. \right\}.
\end{equation}
In the above formulation, for every $p\in \mathcal{V}$ and $e\in \mathcal{E}$,
the unary binary {function} $x_p(\cdot)$ and the pairwise binary {function} $\mathsf{x}_{e}(\cdot)$  {indicate} the labels  assigned to vertex $p$ and 
to the pair of vertices connected by edge $e=(p',q')$ respectively, i.e., 
\begin{align}
&(\forall z^{(p)}\in \mathcal{L}) \qquad\qquad\qquad\qquad\;\;x_p(z^{(p)})=1 \quad\Leftrightarrow \quad p \text{ is assigned label } z^{(p)}\\
&(\forall \mathsf{z}^{(e)}=(z^{(p')},z^{(q')})\in \mathcal{L}^2)\qquad \mathsf{x}_{e}(\mathsf{z}^{(e)})=1 \quad \Leftrightarrow\quad p', q' \text{ are assigned labels } z^{(p')}, z^{(q')}.
\end{align}
{Minimizing with respect to the vector $x$ regrouping all these binary {functions} is equivalent
to searching for an optimal binary vector of dimension
$N = |\mathcal{V}| |\mathcal{L}|+ |\mathcal{E}| |\mathcal{L}|^2$.}
The first constraints in \eqref{eq:mrf_ilp_} simply encode the fact that each vertex must be assigned exactly one label, whereas the rest of the  constraints enforces consistency between unary {functions} $x_p(\cdot)$, $x_q(\cdot)$ and the pairwise {function} $\mathsf{x}_e(\cdot)$ for edge $e=(p,q)$,  
ensuring essentially  that if $x_p(z^{(p)})=x_q(z^{(q)})=1$, then $\mathsf{x}_e(z^{(p)},z^{(q)})=1$.

As mentioned above, our goal will be to decompose  the MRF problem \eqref{eq:orig_MRF_} into easier subproblems (called \emph{slaves}), which, in this case, involve optimizing MRFs defined on subgraphs of $G$. More specifically, let $\{G_m=(\mathcal{V}_m,\mathcal{E}_m)\}_{1\le m \le M}$ be a set of subgraphs that form a decomposition of $G=(\mathcal{V},\mathcal{E}) $ (i.e., $\displaystyle\cup_{m=1}^M \mathcal{V}_m = \mathcal{V}$, $\displaystyle\cup_{m=1}^M \mathcal{E}_m = \mathcal{E}$). On each of these subgraphs, %$G_i$, 
we  define a local MRF with corresponding (unary and pairwise) potentials $\varphi^m=\left\{\{\varphi^m_p\}_{p\in \mathcal{V}_m},\{\thh^m_e\}_{e\in\mathcal{E}_m}\right\}$, whose cost function $f^m(x;\varphi^m)$ is thus given by  
\begin{equation}
f^m(x;\varphi^m) = \sum_{p\in \mathcal{V}_m,\,z^{(p)}\in \mathcal{L}} \varphi^m_p(z^{(p)})\, x_p(z^{(p)})+
  \sum_{e\in \mathcal{E}_m,\, \mathsf{z}^{(e)}\in \mathcal{L}^2}\thh^m_e(\mathsf{z}^{(e)})\, \mathsf{x}_e(\mathsf{z}^{(e)}).
\end{equation}
Moreover, the sum (over  $m$) of the potential functions  $\varphi^{m}$ is ensured to give  back  the potentials $\varphi$ of the original MRF on $G$, i.e.,\footnote{For instance, to ensure \eqref{eq:theta_decomposition} we can simply set: $(\forall m\in\{1,\ldots,M\})$
$\varphi^m_p = \frac{\varphi_p}{|\{|m' \mid p\in \mathcal{V}_{m'}\}|}$ and $\thh^m_e = \frac{\thh_e}{|\{m'\mid e\in \mathcal{E}_{m'}\}|}$.}
\begin{equation}
\label{eq:theta_decomposition}
(\forall p\in \mathcal{V})(\forall e\in \mathcal{E})\quad
\sum_{\substack{m\in \{1,\ldots,M\}:p\in \mathcal{V}_m}}\mspace{-10mu} \varphi^m_p=\varphi_p,\quad
\sum_{\substack{m\in \{1,\ldots,M\}:e\in \mathcal{E}_m}} \mspace{-10mu} \thh^m_{e}=\thh_{e}.
\end{equation}
This  guarantees that  $f=\sum_{m=1}^M f^m$, thus allowing  us to re-express problem \eqref{eq:orig_MRF_}  as follows
\begin{equation}\label{eq:orig_MRF2_}
\minimize{x\in C_G}{\sum_{m=1}^M f^m(x;\varphi^m)}.
\end{equation}

An  assumption that often holds in practice is that minimizing separately  each of the $f^m$ (over  $x$) is easy, but minimizing their sum %the sum
is hard. Therefore, to leverage this fact, we introduce, for every $m\in\{1,\ldots,M\}$, an \emph{auxiliary} \emph{copy} $x^m\in C_{G_m}$ for the variables  of  the local MRF defined on $G_m$, which are thus constrained to coincide with the corresponding  variables in vector $x$, i.e., it holds $x^m=x_{|G_m}$, where  $x_{|G_m}$ is used to denote the subvector of $x$ containing  only those variables associated with vertices and edges of subgraph $G_m$.
In this way,  Problem \eqref{eq:orig_MRF2_}  can be transformed into 
\begin{align}\label{eq:dp2_}
\minimize{x\in C_G,\{x^m\in C_{G_m}\}_{1\le m \le M}}{} &  \sum_{m=1}^M f^m(x^m;\varphi^m)\nonumber\\
 \mathrm{s.t.} &\ \ \big(\forall m\in\{1,\ldots,M\}\big)\qquad x^m=x_{|G_m}.
\end{align}
By considering the dual of   \eqref{eq:dp2_}, using a technique similar to the one described in
framebox on page~\pageref{fb:cons_sharing}, and noticing that
\begin{equation}
x\in C_G \quad \Leftrightarrow \quad (\forall m \in \{1,\ldots,M\})\quad x^m\in C_{G_m},
\end{equation}
we finally end up with the following  problem: % \begin{equation}
\begin{equation}\label{eq:drel_}
\maximize{ (v^m)_{1\le m \le M} \in \Lambda}{\sum_{m=1}^M h^m(v^m)},
\end{equation}
where, for every $m\in \{1,\ldots,M\}$, the dual variable $v^m$ consists of
$\left\{\{v^m_p\}_{p\in \mathcal{V}_m},\{\mathsf{v}^m_e\}_{e\in\mathcal{E}_m}\right\}$
similarly to $\varphi^m$, and function $h^m$ is related to the following optimization of a slave MRF on $G_m$:
\begin{equation}
h^m(v^m)=\min_{x^m\in C_{G_m}} f^m(x^m;\varphi^m+ v^m).
\end{equation}
The feasible set $\Lambda$ {is given by} 
\begin{equation}
\Lambda  = \left\{
   {v=\left\{\{v^m_p\}_{p\in \mathcal{V}_m},\{\mathsf{v}^m_e\}_{e\in\mathcal{E}_m}\right\}_{1\le m \le M}}
        \left|
        \begin{array}{l}
   (\forall p \in \mathcal{V})(\forall z^{(p)}\in \mathcal{L})\quad\qquad\; \sum\limits_{\substack{m\in \{1,\ldots,M\}:p\in \mathcal{V}_m}} v^m_p(z^{(p)})=0,  \\
   (\forall e \in \mathcal{E})(\forall \mathsf{z}^{(e)}\in \mathcal{L}^2)
\quad\qquad\sum\limits_{\substack{m\in \{1,\ldots,M\}:e\in \mathcal{E}_m}} \mathsf{v}^m_e(\mathsf{z}^{(e)}) = 0   \\
{(\forall m \in \{1,\ldots,M\})(\forall p\in \mathcal{V})\;\; v^m_p(\cdot)\colon\mathcal{L}\mapsto\mathbb{R}}\\
{(\forall m \in \{1,\ldots,M\})(\forall e\in \mathcal{E})\;\; \mathsf{v}^m_e(\cdot)\colon\mathcal{L}^2\mapsto \mathbb{R}}
\end{array} 
\right. 
\right\}.
\end{equation}

The above dual problem provides a relaxation to the original problem \eqref{eq:orig_MRF_}-\eqref{eq:mrf_ilp_}. Furthermore, note that this relaxation leads to a convex optimization problem,\footnote{{In order to see this, notice that $h^m(v^m)$ is equal to a pointwise minimum of a set of linear functions of $v^m$, and thus it is a concave function.}} although the original one is not. 
As such, it can always be solved in an optimal manner. 
A possible way of doing this consists of using a projected subgradient method. 
Exploiting the form of the projection onto  the vector space $\Lambda$ yields Algorithm \ref{dd_pseudocode} 
where $(\gamma_n)_{n\in\NN}$ is a summable sequence of positive step-sizes and 
$\left\{\{\widehat{x}^m_{p,n}\}_{p\in \mathcal{V}_m},\{\widehat{\mathsf{x}}^m_{e,n}\}_{e\in\mathcal{E}_m}\right\}$
corresponds to a subgradient of function $h^m$ with $m\in \{1,\ldots,M\}$
computed at iteration $n$ \cite{komodakis_pami2011}.
Note that this algorithm requires \emph{only  solutions to  local subproblems} to be computed, which is, of course, a task much easier that furthermore can be executed in
a parallel manner. The solution to the master MRF is filled in from local  solutions 
$\left\{\{\widehat{x}^m_{p,n}\}_{p\in \mathcal{V}_m},\{\widehat{\mathsf{x}}^m_{e,n}\}_{e\in\mathcal{E}_m}\right\}_{1\le m \le M}$
 after convergence of the algorithm.
{\linespread{1.5}
\begin{algorithm}
{\small
\caption{Dual decomposition for MRF optimization.}\label{dd_pseudocode}
\begin{equation*}
\begin{array}{l}
\text{Choose a decomposition $\{G_m=(\mathcal{V}_m,\mathcal{E}_m)\}_{1 \le m \le M}$ of  $G$}\\
\text{Initialize potentials of slave MRFs:}\\
(\forall m \in \{1,\ldots,M\})(\forall p\in \mathcal{V}_m)\;
\varphi^m_{p,0} = \frac{\varphi_p}{|\{m'\mid p\in \mathcal{V}_{m'}\}|}, 
(\forall e\in \mathcal{E}_m)\;
\thh^m_{e,0} = \frac{\thh_e}{|\{m'\mid e\in \mathcal{E}_{m'}\}|} \\
\text{for}\;n=0,\ldots\\
\left\lfloor
\begin{array}{l}
\text{Compute minimizers of slave MRF problems:}\;
(\forall m \in \{1,\ldots,M\})\;
%\widehat{x}^m_n \in\argmind{x^m\in C_{G_m}}{f^m(x^m;\varphi^m_n)}\qquad\qquad\\
\left\{\{\widehat{x}^m_{p,n}\}_{p\in \mathcal{V}_m},\{\widehat{\mathsf{x}}^m_{e,n}\}_{e\in\mathcal{E}_m}\right\}\in\Argmind{x^m\in C_{G_m}}{f^m(x^m;\varphi^m_n)}\qquad\qquad\\
\text{Update potentials of slave MRFs:}\\
(\forall m \in \{1,\ldots,M\})(\forall p\in \mathcal{V}_m)\;
\varphi^m_{p,n+1} = \varphi^m_{p,n+1}+\gamma_n \left( \widehat{x}^m_{p,n} - \frac{\sum_{m:\,p\in \mathcal{V}_m}{\widehat{x}^m_{p,n}}}{|\{m'\mid p\in \mathcal{V}_{m'}\}|} \right)\\
(\forall m \in \{1,\ldots,M\})(\forall e\in \mathcal{E}_m)\;
 \thh^m_{e,n+1} = \thh^m_{e,n}+\gamma_n \left( \widehat{\mathsf{x}}^m_{e,n} - \frac{\sum_{m:\,e\in \mathcal{E}_m}{\widehat{\mathsf{x}}^m_{e,n}}}{|\{m'\mid p\in \mathcal{V}_{m'}\}|} \right).\\
%\text{$x^{(t)}$ = fill-in solution for master MRF from local  solutions  $\{\hat{x}^i\}$}\\
\end{array}
\right.
\end{array}
\end{equation*}}
\end{algorithm}}

For a better intuition for the  updates of variables 
$\left\{\{\varphi^m_{p,n}\}_{p\in \mathcal{V}_m},\{\thh^m_{e,n}\}_{e\in\mathcal{E}_m}\right\}_{1\le m \le M,n\in\NN}$
 in Algorithm \ref{dd_pseudocode},  we should note that  their aim is essentially to bring a consensus among the solutions of the local subproblems. In other words, they try to adjust the potentials of the slave MRFs so that  in the end the corresponding local solutions are consistent with each other, i.e., all variables corresponding to a common vertex or edge are assigned the same value by the different  subproblems.
If this condition is satisfied (i.e.,  there is a full consensus) then the overall solution that results from combining the consistent local solutions is guaranteed to be optimal. In general,  though, this  might not  always be true given that the above procedure is solving only a  \emph{relaxation} of the original NP-hard problem.

\begin{center}
\framebox{
\begin{minipage}{14cm}
\footnotesize\sffamily 
{MASTER-SLAVE COMMUNICATION}

During dual decomposition a communication between a master process and the slaves (local subproblems)  can be thought of as taking place, which can  also be interpreted as a resource allocation/pricing stage.
\begin{center}
\includegraphics[width=0.6\linewidth]{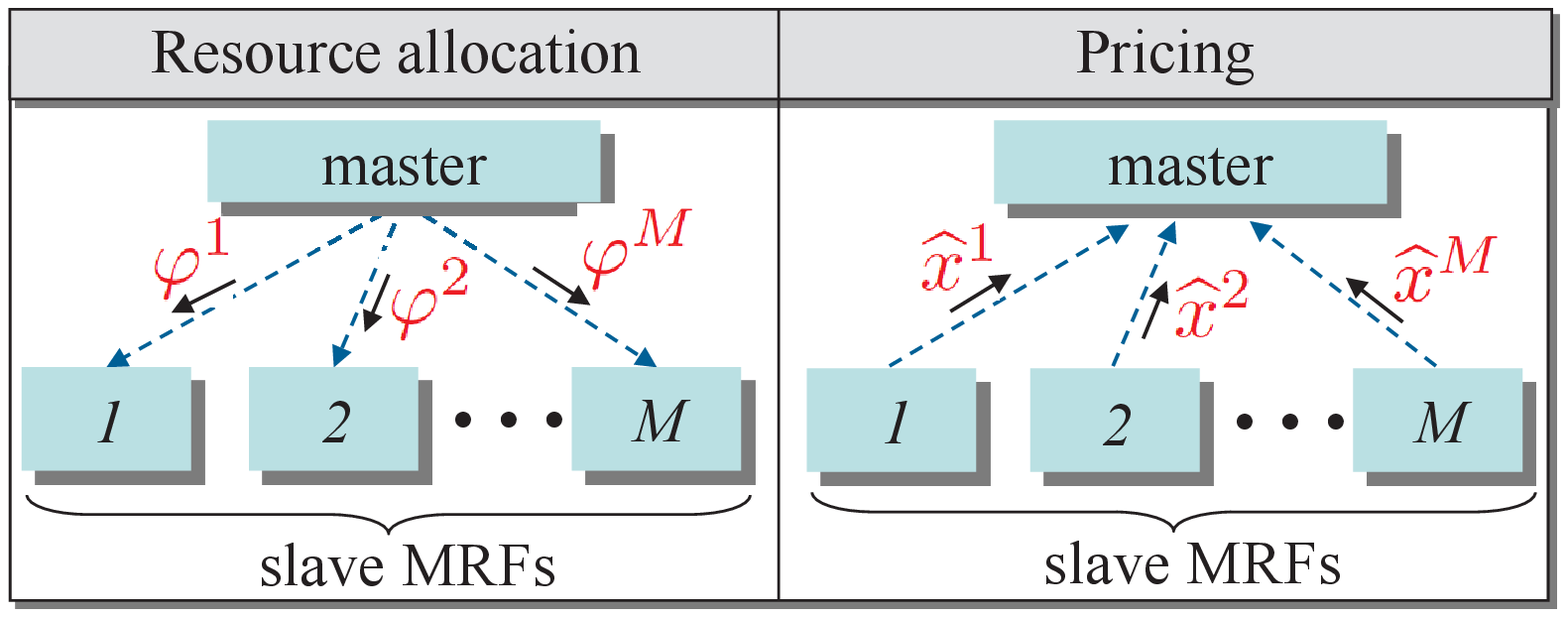}
\vspace{-10pt}
\end{center}
 \textbf{Resource allocation:} At each iteration, the master assigns new MRF potentials (i.e., resources) $(\varphi^m)_{1\le m \le M}$ to the slaves based on the current local solutions $(\widehat{x}^m)_{1\le m \le M}$. 

\textbf{Pricing:} The slaves respond by adjusting their local solutions $(\widehat{x}^m)_{1\le m \le M}$ (i.e., the prices) so as to maximize their welfares based on the newly assigned resources $(\widehat{x}^m)_{1\le m \le M}$.
\end{minipage}
}
\end{center}

\begin{center}
\framebox{\label{fb:dd_rel}
\begin{minipage}{14cm}
\footnotesize\sffamily 
{DECOMPOSITIONS AND RELAXATIONS}

Different decompositions  can lead to different relaxations and/or can affect the speed of convergence. We show below, for instance, 3 possible decompositions for an MRF assumed to be defined on a $5\times 5$ image grid.

\begin{center}
\includegraphics[width=0.5\linewidth]{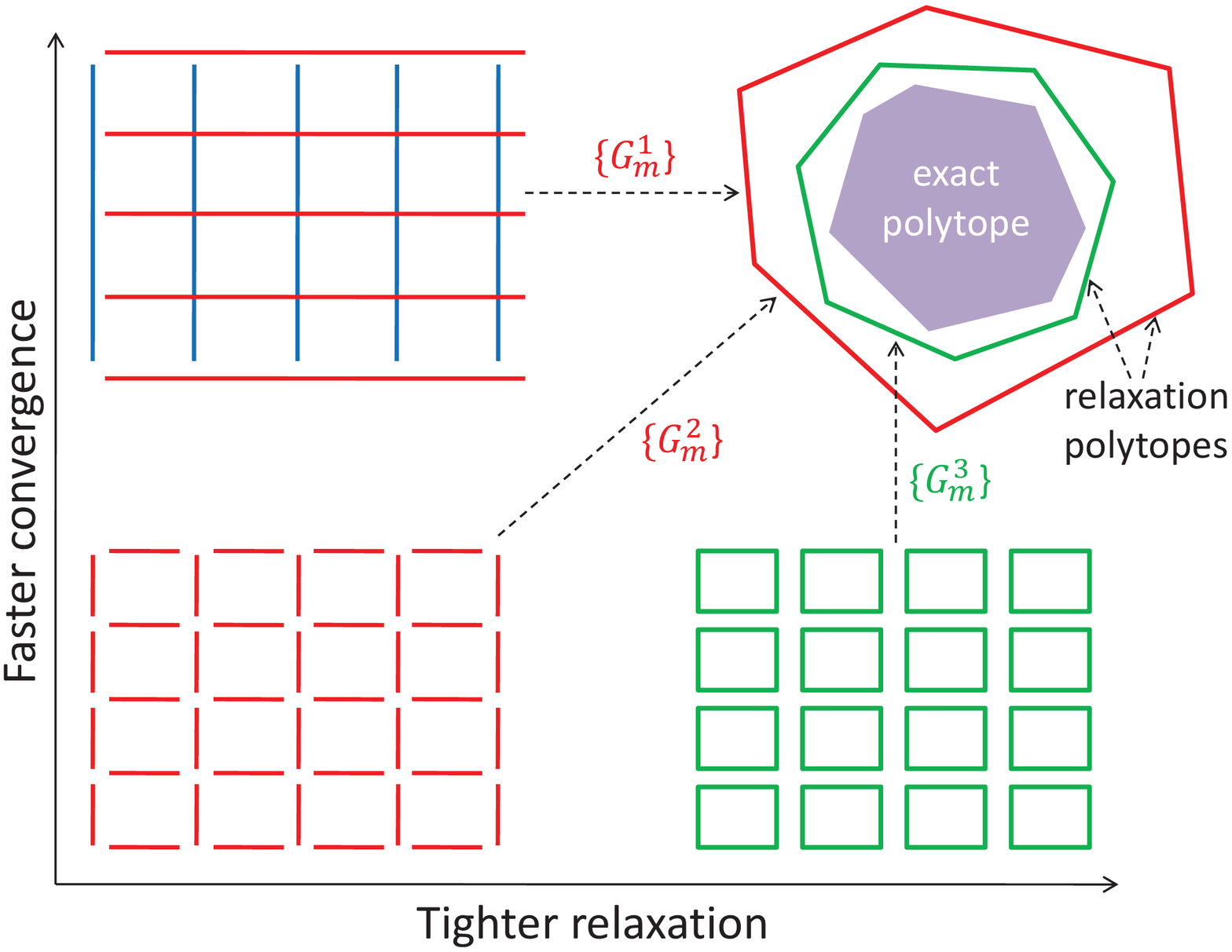}
\vspace{-10pt}
\end{center}
Decompositions $\{G^1_m\}, \{G^2_m\}, \{G^3_m\}$ consist respectively of one suproblem per row and column, one subproblem per edge, and one subproblem per $2\times 2$ subgrid of the original $5\times 5$ grid. Both  $\{G^1_m\}$ and $\{G^2_m\}$ (due to using solely   subgraphs that are trees) lead to the same LP relaxation of  \eqref{eq:orig_MRF_}, whereas   $\{G^3_m\}$ leads to a  relaxation that is tighter (due to containing loopy subgraphs).

 On the other hand, decomposition $\{G^1_m\}$ leads to faster convergence compared with $\{G^2_m\}$ due to using larger subgraphs that allow a faster propagation of information during message-passing.
\end{minipage}
}
\end{center}

Interestingly, if we choose to use a decomposition consisting only of %tree-structured
 subgraphs that are trees,
then the resulting  relaxation can be shown to actually coincide with the standard LP-relaxation of linear integer program \eqref{eq:orig_MRF_} (generated by replacing the integrality constraints with non-negativity constraints on the variables).
This also means that when this LP-relaxation is tight, an optimal MRF solution is computed. This, for instance, leads to the  result that \emph{dual decomposition approaches can estimate a globally  optimal solution for binary submodular MRFs} (although it should be noted  that much faster graph-cut based techniques exist for  submodular problems of this type - see framebox on page~\pageref{fb:gc}). Furthermore, when using subgraphs that are trees, a minimizer to each slave problem can be  computed efficiently by applying the Belief Propagation algorithm \cite{pearl}, which is a message-passing method. Therefore, in this case, Algorithm~\ref{dd_pseudocode} essentially reduces to a continuous exchange of messages between the nodes of  graph $G$. Such an algorithm   relates to or generalizes various other message-passing approaches \cite{TRW_wainwright, TRW_S, werner_pami, globerson,  yanover, DBLP:journals/tit/HazanS10}. In general, besides tree-structured subgraphs, other types of decompositions or subproblems can be used as well (such as binary planar problems, or problems on loopy subgraphs with small tree-width, for which MRF optimization can still be solved efficiently), which can lead to even tighter relaxations  (see framebox on page~\pageref{fb:dd_rel}) \cite{JegelkaBS13, Schraudolph10, OsokinVK11, DBLP:conf/uai/YarkonyMIF11, sontag2, komodakis_eccv08}.

\begin{center}
\framebox{\label{fb:gc}
\begin{minipage}{15cm}
\footnotesize\sffamily 
{GRAPH-CUTS AND MRF OPTIMIZATION}

For certain  MRFs, optimizing their cost is known to be equivalent to  solving a polynomial  mincut problem \cite{BorosH02, KolmogorovZ04}. These are exactly all the binary MRFs ($|\mathcal{L}|=2$)
 with submodular pairwise potentials such that, for every $e\in \mathcal{E}$,
\begin{equation}
\label{eq:subm}
\thh_e(0,0)+\thh_e(1,1)\le \thh_e(0,1)+\thh_e(1,0).
\end{equation}
Due to \eqref{eq:subm}, the cost  $f(x)$ of a binary labeling $x= (x^{(p)})_{1\le p \le |\mathcal{V}|}\in \{0,1\}^{|\mathcal{V}|}$  
for such MRFs can always be written (up to an additive constant) as 
\begin{equation}
f(x)=\sum_{p\in \mathcal{V}_P}a_{p}x^{(p)}+ \sum_{p\in \mathcal{V}_N}a^{(p)}(1-x^{(p)})+\sum_{(p,q)\in \mathcal{E}}a_{p,q}x^{(p)}(1-x^{(q)}),
\end{equation}    
 where all coefficients $(a_p)_{p\in \mathcal{V}}$ and $(a_{p,q})_{(p,q)\in \mathcal{E}}$ are nonnegative 
($\mathcal{V}_P \subset \mathcal{V}$, $\mathcal{V}_N\subset\mathcal{V}$).

In this case, we can  associate to $f$ a capacitated network that has vertex set  $\mathcal{V}_f=\mathcal{V}\cup\{s, t\}$.
A source vertex $s$ and a sink one $t$ have thus been added. The new edge set $\mathcal{E}_f$ is deduced from the one used
to express $f$:
\[
\mathcal{E}_f=\{(p,t)\mid p\in \mathcal{V}_P\}\cup \{(s,p)\mid p\in \mathcal{V}_N\} \cup \mathcal{E},
\]
and its edge capacities are defined as $(\forall p \in \mathcal{V}_P \cup \mathcal{V}_N)$ $c_{p,t}=c_{s,p}=a_p$, and
$(\forall (p,q)\in \mathcal{E})$ $c_{p,q}=a_{p,q}$.

A one-to-one correspondence between $s$-$t$ cuts and MRF labelings then exists: 
\[x\in \{0,1\}^{|\mathcal{V}|}\leftrightarrow \mathrm{cut}(x)=\{s\}\cup \{p \mid x^{(p)}=1\}\]
for which it is easy to see that 
\[f(x)=\sum_{u\in \mathrm{cut}(x),\upsilon\notin \mathrm{cut}(x)}c_{u,\upsilon}=\text{cost of }\mathrm{cut}(x)\enspace.\]
Computing a mincut, in this case,  solves the LP relaxation of \eqref{eq:orig_MRF_}, which is  tight, whereas  computing a max-flow
 solves  the dual LP. % a dual relaxation \eqref{eq:dual_}.
\end{minipage}
}
\end{center}

Furthermore, besides the projected subgradient method, one can alternatively apply an ADMM scheme for solving relaxation \eqref{eq:drel_} (see Section \ref{se:ADMM}). The main difference, in this case, is that the optimization of a slave MRF problem is performed by solving a (usually simple) local quadratic problem, which can again be solved efficiently for an appropriate choice of the decomposition (see Section \ref{se:extsplit}). This method again penalizes disagreements among slaves, but it does so even more aggressively than the subgradient method 
since there is no longer a requirement for step-sizes $(\gamma_n)_{n\in \NN}$ converging to zero.
Furthermore, alternative smoothed accelerated schemes exist and can be applied as well \cite{Jojic+al:ICML10, SavchynskyyKSS11, SavchynskyySKS12}.

\section{Applications} \label{se:appli}
Although the presented primal-dual algorithms can be applied virtually  to any area
where optimization problems have to be solved, we now mention a few common applications of these
techniques.

\subsection{Inverse problems}
For a long time, convex optimization approaches have been successfully used for solving
inverse problems such as signal restoration, signal reconstruction, or interpolation of missing data.
Most of the time, these problems are ill-posed and, in order to recover the signal
of interest in a satisfactory manner, some prior information needs to be introduced. To do this, an objective function
can be minimized which includes a data fidelity term modelling knowledge about the noise statistics and 
possibly involves a linear observation matrix (e.g. a convolutive blur), and a regularization (or penalization) term which corresponds to the additional 
prior information.  This formulation can also often be justified statistically as the determination of a Maximum A Posteriori (MAP) estimate.
In early developed methods, in particular in Tikhonov regularization, a quadratic penalty function is employed.
Alternatively, hard constraints can be imposed on the solution (for example, bounds on the signal values), leading to
signal feasibility problems. Nowadays, a hybrid regularization \cite{Pustelnik_N_2011_j-ieee-tip_par_pai}
may be prefered so as to combine various kinds of regularity measures, possibly computed for different representations
of the signal (Fourier, wavelets,...), some of them like total variation \cite{Rudin_L_1992_tv_atvmaopiip} and its nonlocal extensions \cite{Zhang_X_2010_j-siam-is_Bregmanized_nrfd} being taylored for 
preserving discontinuities such as image edges. In this context, constraint sets can be translated into penalization terms
being equal to the indicator functions of these sets (see \eqref{e:defindicator}). Altogether, these lead to global cost functions which can be quite involved,
often with many variables, for which the splitting techniques described in Section \ref{se:extsplit} are very useful. An extensive literature
exists on the use of ADMM methods for solving inverse problems
(e.g., see \cite{Giovanelli_JF_2005-astoastro-positive-dse,Goldstein_T_2009_j-siam-is_split_bml,Figueiredo_M_2009_ssp_Deconvolution_opiuvsaalo,Figueiredo_M_2010_t-ip_restoration_piado,AfonsoBF11}).
With the advent of more recent primal-dual algorithms, many works have been mainly focused on image recovery applications 
\cite{Chambolle_A_2010_first_opdacpai,Esser_E_2010_j-siam-is_gen_fcf,He_B_2012_j-siam-is_conv_apd,Pock_T_2008_p-iccv_diagonal_pffo,Goldstein_T_2013_adaptive_pdh,Bot_R_2014_jmiv_conv_primal_apd,Loris_I_2011_generalization_ist,Jezierska_A_2012_p-icassp_primal_dps,Chen_P_2013_j-inv-prob_prim_dfp,Bec_S_2014_j-nonlinear-conv-anal_alg_sps,Bonettini_S_2012_j-math-imaging-vis_convergence_pdhgatvir,Repetti_A_2012_p-eusipco_penalized_wlsardcsdn,Harizanov_S_2013_Epigraphical_psls,Teuber_T_2014_j-inv-prob_Minimization_pes,Burger_M_2014_First_oavip}. Two illustrations are now provided.

In \cite{Couprie_C_j-siam-is_dual_ctvb}, a generalization of the total variation is defined for an arbitrary graph in order to address a variety of inverse problems. {For denoising applications, 
the optimization problem to be solved is of the form \eqref{e:primalvar}
where
\begin{equation}
f = 0,\qquad g = \sigma_C,\qquad h\colon x \mapsto \frac12 \|x-y\|^2,
%the cost function
%which is minimized reads
%\begin{equation}
%(\forall x \in \HH)\qquad f(x) = \sigma_C(Lx)+\frac12 \|x-y\|^2
\end{equation}
%where 
$x$ is a vector of variables associated with each vertex of a weighted graph, and $y\in \HH$ is a vector of data observed at each vertex.}
The matrix $L\in \RR^{K\times N}$ is equal to $\operatorname{Diag}(\sqrt{\omega}_1,\ldots,\sqrt{\omega}_K)\,A$ where $(\omega_1,\ldots\omega_K) \in \RP^K$ is the vector of edge weights and $A\in \RR^{K\times N}$ is the graph incidence matrix playing a role similar to a gradient operator on the graph. The set $C$ is defined as an intersection of closed semi-balls in such a way that
its support function $\sigma_C$ (see \eqref{e:defsupfonc}) allows us to define a class of functions extending the total variation semi-norm (see \cite{Couprie_C_j-siam-is_dual_ctvb} for more details). 
Good image denoising results can be obtained by building the graph in a nonlocal manner following the strategy in \cite{Zhang_X_2010_j-siam-is_Bregmanized_nrfd}. Results obtained for Barbara image are
displayed in Fig. \ref{fig:nonlocal}. Interestingly, the ability of methods such as those presented in Section \ref{se:FBFPF} to circumvent matrix inversions leads to a significant
decrease of the convergence time for irregular graphs in comparison with algorithms based on the Douglas-Rachford iteration or ADMM (see Fig. \ref{fig:Compa_FBF_PPXA}).

\begin{figure}[ht]
\centering
\subfigure[Original image]{\includegraphics[scale=0.29]{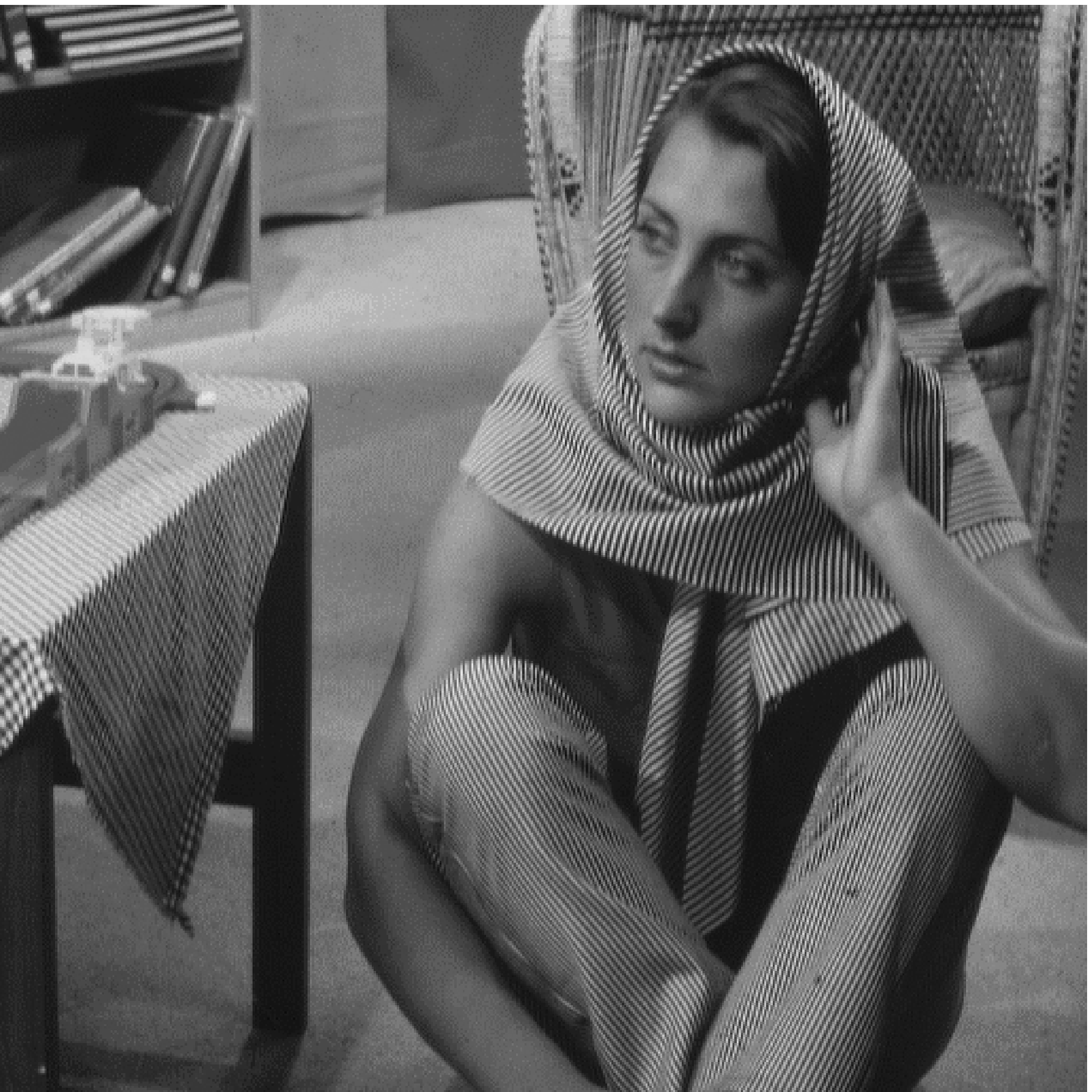}}
\subfigure[Noisy~{\scriptsize  SNR~=~$14.47$ dB}]{\includegraphics[scale=0.29]{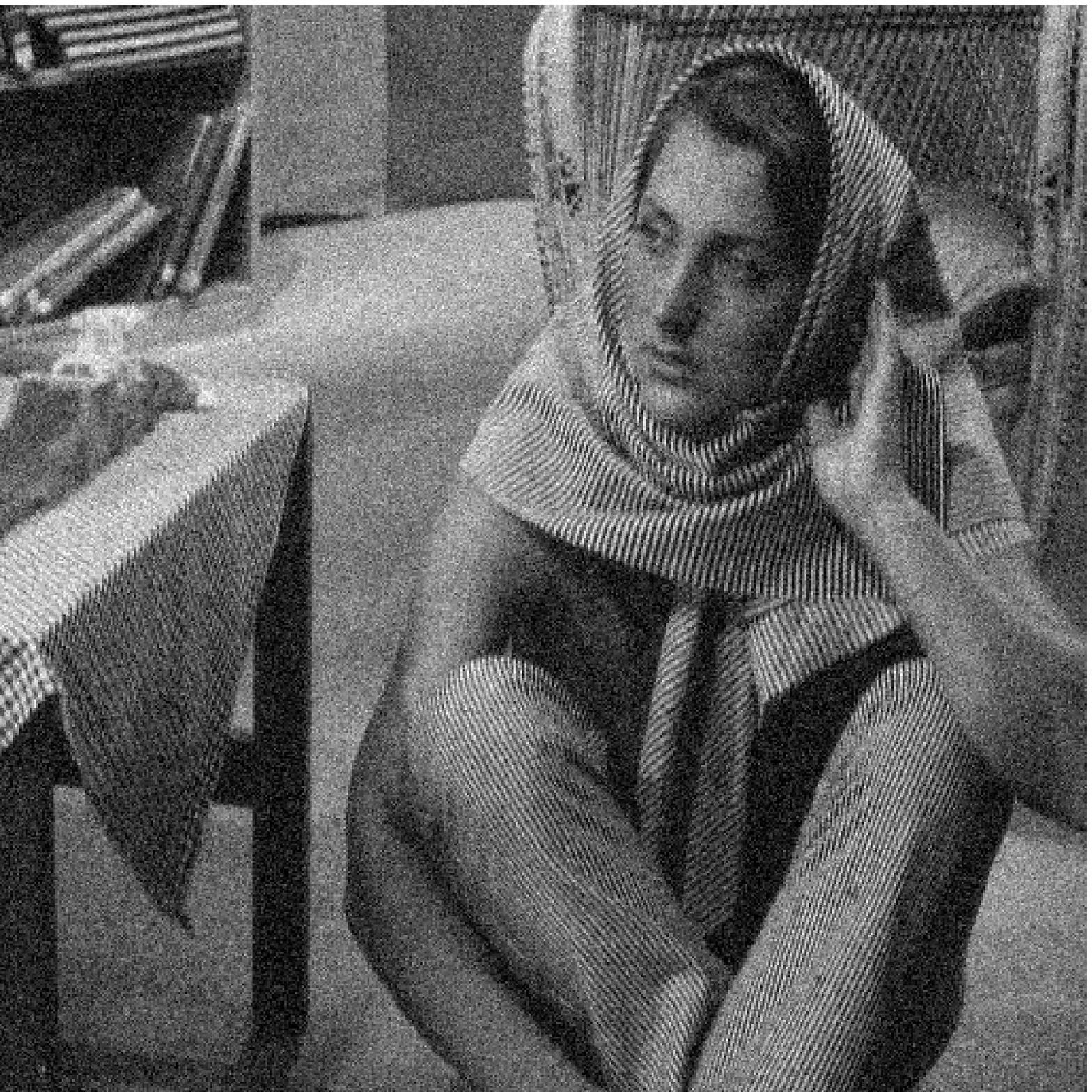}}
\subfigure[{\scriptsize Nonlocal TV~{\scriptsize SNR~=~$20.78$~dB}}]{\includegraphics[scale=0.29]{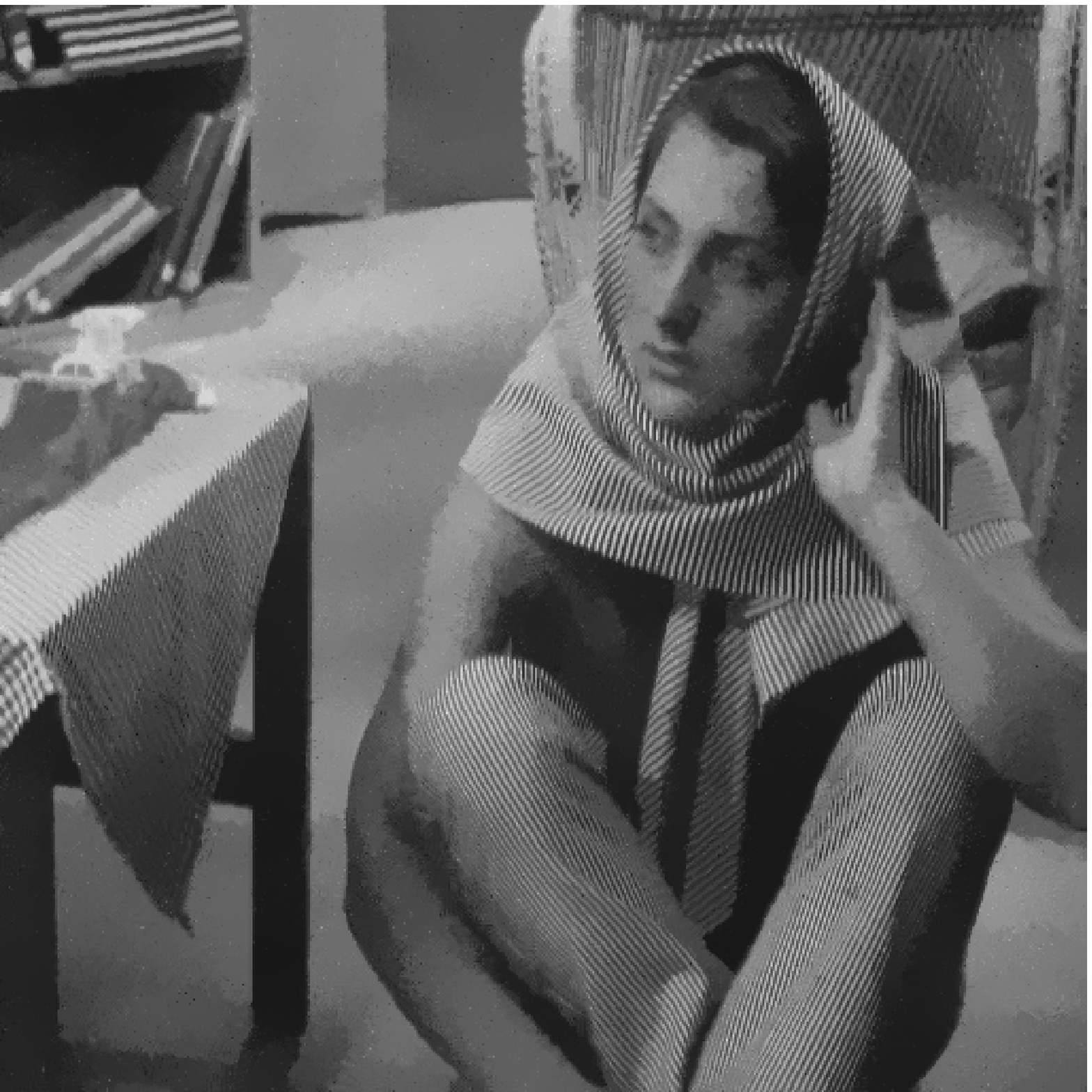}}
%\vspace{-10pt}
\caption{{\small Nonlocal denoising (additive white zero-mean Gaussian noise with variance $\sigma^2=20$).}}
\label{fig:nonlocal}
\end{figure}

\begin{figure}[h]
\begin{center}
%\begin{tabular}{cc}
% \includegraphics[scale=0.37]{Comp_FBF_PPXA_nb_iter_regular.png} &
%\begin{minipage}{0.48\textwidth}
%\subfigure[regular
%  graph]{\includegraphics[scale=0.38]{images/logtime_reg.pdf}}
%\end{minipage}&
%\includegraphics[scale=0.37]{iteration.png} &
%\begin{minipage}{0.48\textwidth}
%\subfigure[non-regular
%  graph]{
\includegraphics[scale=0.48]{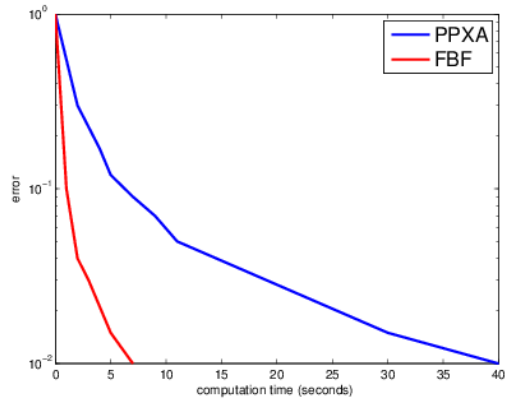}
%}
%\end{minipage}
\end{center}
%\end{tabular}
\caption{\small Comparison of the convergence speed of a Douglas-Rachford based algorithm (PPXA \cite{Combettes_PL_2008_j-ip_proximal_apdmfscvip}) (blue) and 
an FBF-based primal-dual algorithm (red) for image denoising using a non-regular graph, Matlab implementation on an
Intel Xeon 2.5GHz 8-core system.}
\label{fig:Compa_FBF_PPXA}
\end{figure}

Another application example of primal-dual proximal algorithms is Parallel Magnetic Resonance Imaging (PMRI) reconstruction. A set of measurement vectors $(z_j)_{1 \leq j\leq J}$ is acquired
from $J$ coils. These observations are related to the original full FOV (Field Of View) image $\overline{x}\in\CC^N$ corresponding to a spin density. An estimate of $\overline{x}$
is obtained by solving the following problem: 
{
\begin{equation}\label{e:probMRI}
\minimize{x\in \CC^N}{f(x)+g(L x)+\underbrace{\sum_{j=1}^J \| \Sigma F S_j x  - z_j\|^2_{\Lambda_j^{-1}}}}_{\mbox{$h(x)$}}
\end{equation}
}
where $(\forall j \in \{1,\ldots,J\})$ $\|\cdot\|^2_{\Lambda_j^{-1}} = (\cdot)^{\rm H}\Lambda_j^{-1} (\cdot)$,
$\Lambda_j$ is the noise covariance matrix for the $j$-the channel,
$S_{j}\in \CC^{N\times N}$ is a diagonal matrix modelling the sensitivity of the coil, $F \in \CC^{N\times N}$  is a 2D discrete Fourier transform, 
$\Sigma \in \{0,1\}^{\lfloor\frac{N}{R}\rfloor\times N}$ is a subsampling matrix,
$g \in \Gamma_0(\CC^K)$ is a sparsity measure (e.g. a weighted $\ell_1$-norm), $L\in \CC^{K\times N}$ is a (possibly redundant) frame analysis operator, and {$f$ is the indicator function
of a vector} subspace of $\CC^N$ serving to set to zero the 
image areas corresponding to the background.\footnote{$(\cdot)^{\rm H}$ denotes the transconjugate operation and $\lfloor \cdot \rfloor$ designates the 
lower rounding operation.} Combining suitable subsampling strategies in the k-space with the use of an array of coils allows us to reduce the acquisition time 
while maintaining a good image quality. The subsampling factor $R > 1$ thus corresponds to an \emph{acceleration factor}. For a more detailed account on the considered approach,
the reader is refered to \cite{Chaari_L_2011_j-media_wav_brr,Florescu_A_2014_j-sp_majorize_mmg} and the references therein. Reconstruction results are shown in Fig.~\ref{Fig:Reconst}. Fig.~\ref{fig:convergencePMRI} also allows us to evaluate the convergence time for various algorithms. It can be observed that smaller differences between the implemented primal-dual strategies are apparent in this example.
Due to the form of the subsampling matrix, the matrix inversion involved at each iteration of ADMM however requires to make use of a few subiterations of a linear conjugate gradient method.

\begin{figure}
\subfigure[]{
\includegraphics[width=4cm]{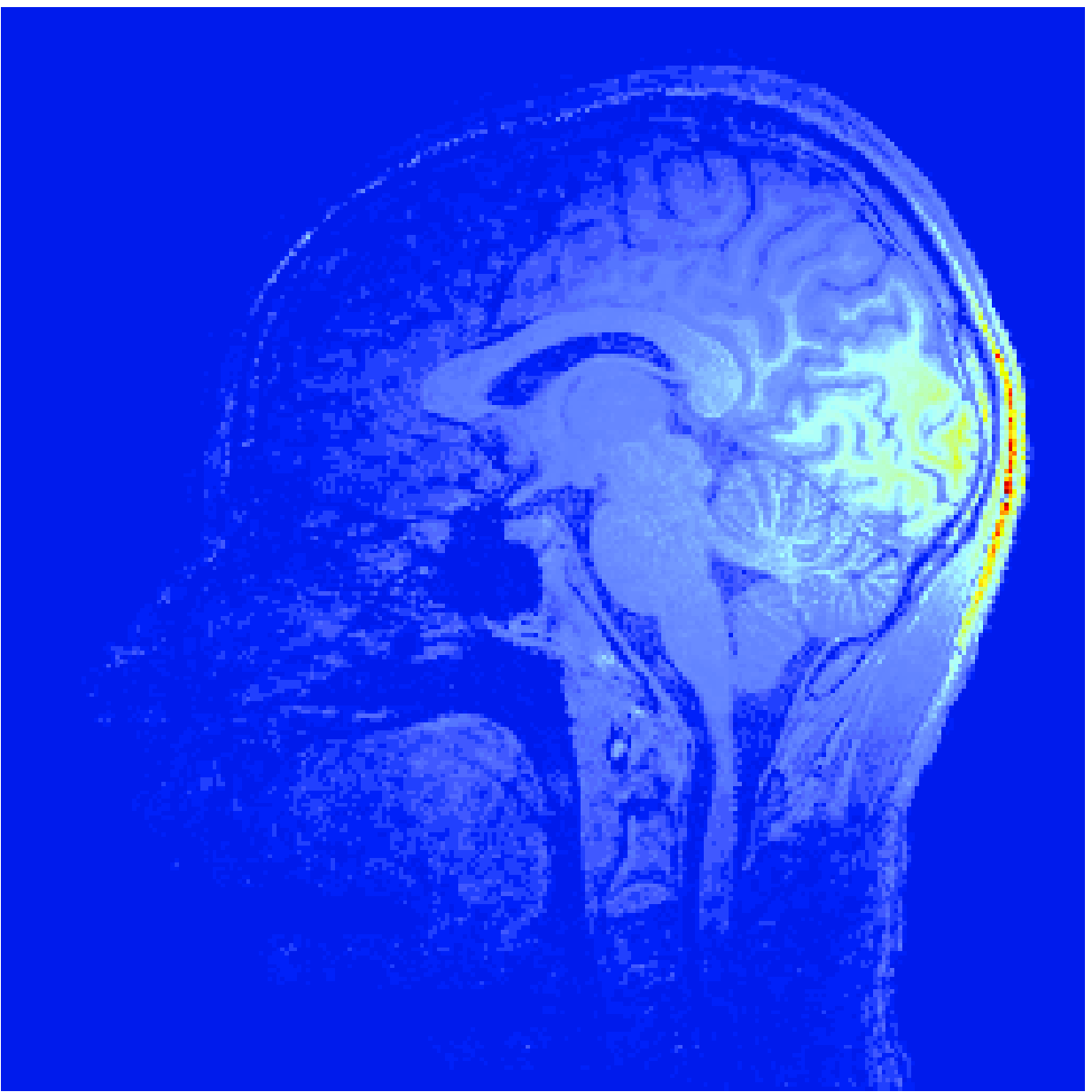}
\includegraphics[width=4cm]{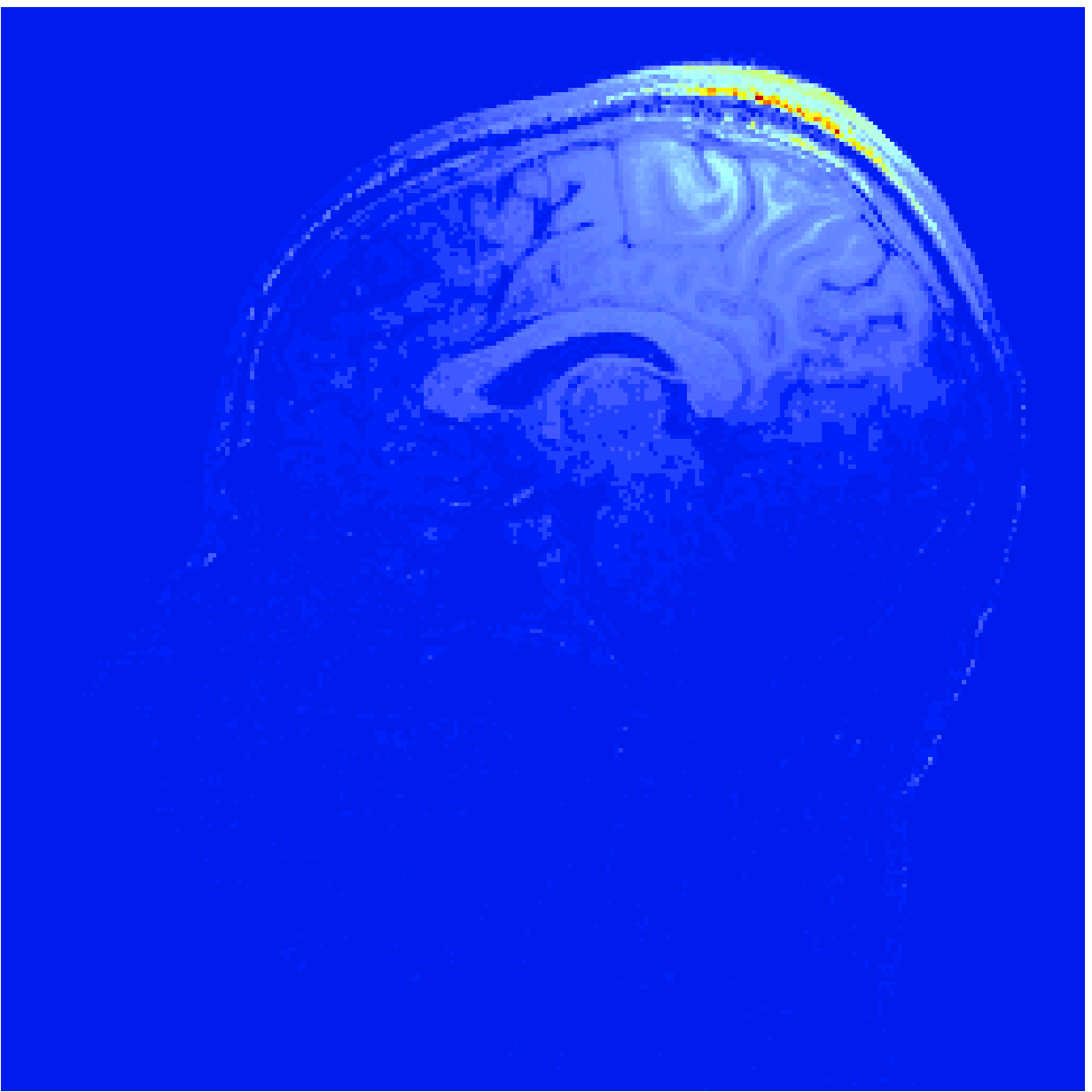}}
\hspace{10pt}
\subfigure[]{
\includegraphics[width=4cm]{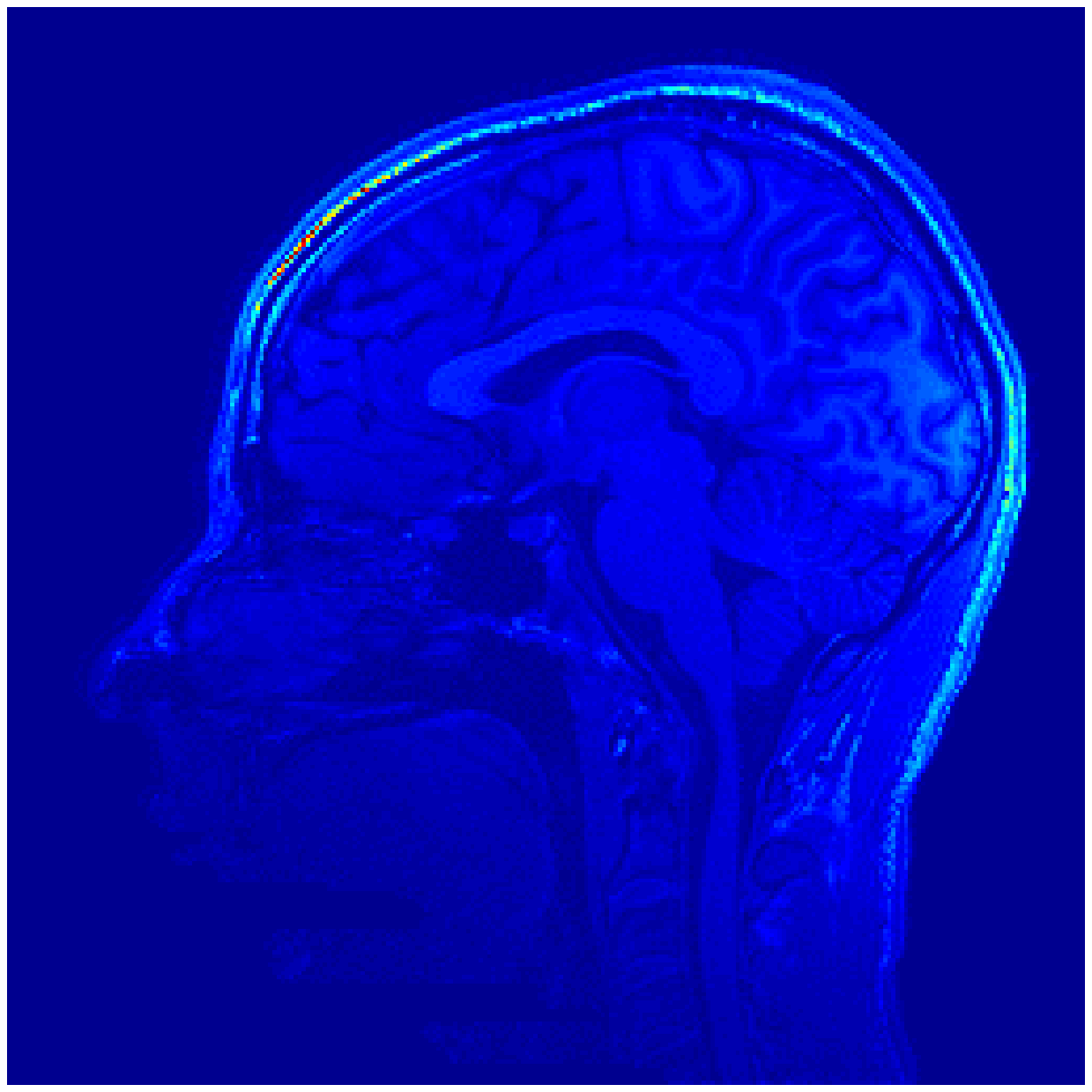}
\includegraphics[width=4cm]{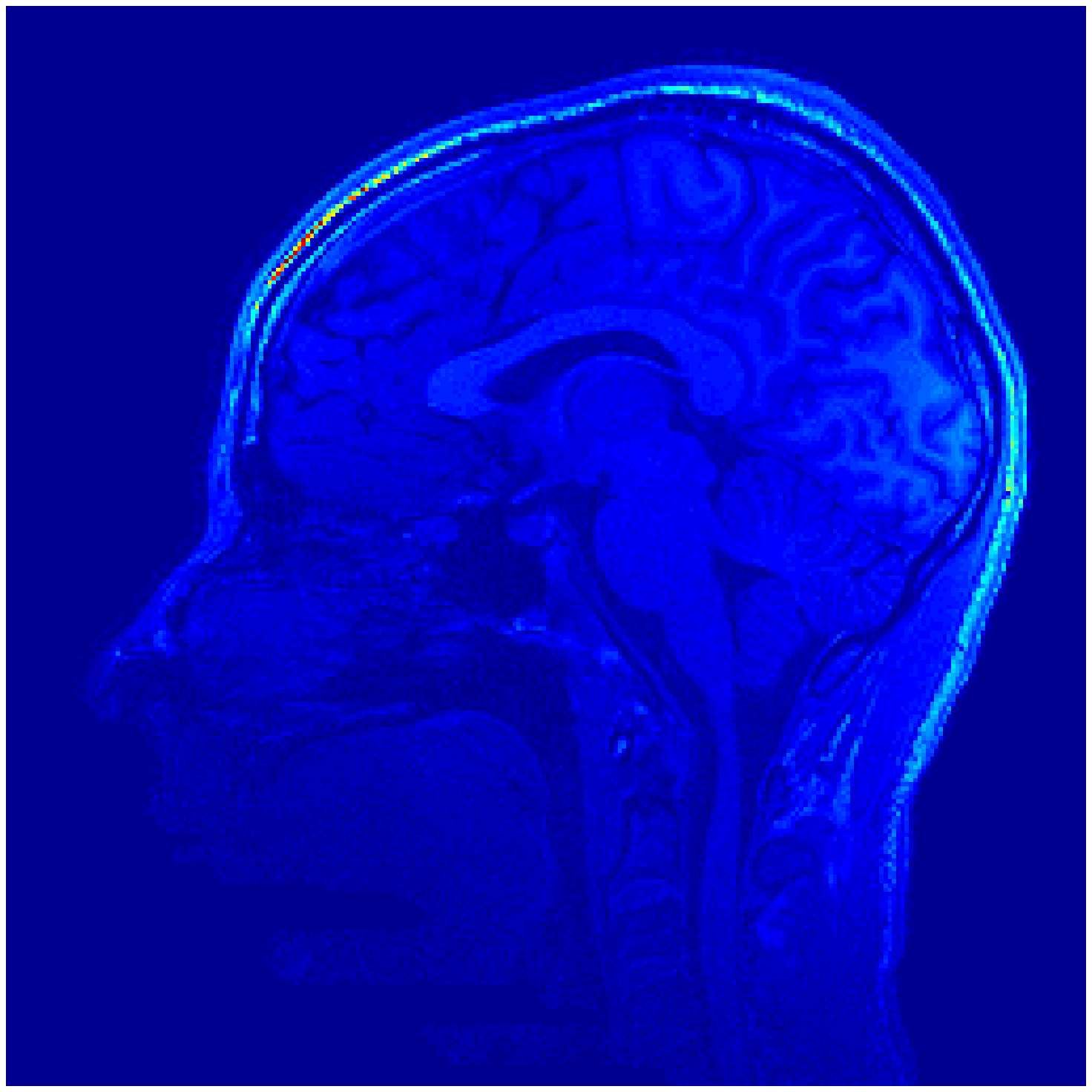}}
\caption{\small \textbf{(a)} Effects of the sensitivity matrices in the spatial domain in the absence of subsampling: the moduli of the images corresponding
to $(S_j \overline{x})_{2\leq j \leq 3}$ are displayed for 2 channels out of 32.
\textbf{(b)} Reconstruction quality:  moduli of the original slice $\overline{x}$ and the reconstructed one  with SNR~=~$20.03$ dB (from left to right)
 using polynomial sampling of order 1 with $R=5$, a wavelet frame, and an $\ell_1$ regularization.}
\label{Fig:Reconst}
\end{figure}

% \begin{figure}
% \begin{center}
% \begin{tabular}{c@{}c@{}}
% \includegraphics[width=4cm]{images/S12.eps}&
% %\includegraphics[width=3cm]{S16}
% \includegraphics[width=4cm]{images/S20.eps}\\
% %\includegraphics[width=4cm]{images/ReferenceImage_S.eps}&
% \includegraphics[width=4cm]{images/Orig_PMRI.eps}&
% \includegraphics[width=4cm]{images/im_poly1_cvx_wav_82.eps}
% \end{tabular} 
% \caption{Effects of the sensitivity matrices in the spatial domain in the absence of subsampling (top line) : the moduli of the images corresponding
% to $(S_j \overline{x})_{2\leq j \leq 3}$ are displayed for 2 channels out of 32.
% Reconstruction quality (bottom line) :  moduli of the original slice $\overline{x}$ and the reconstructed one  SNR~=~$20.57$ dB (from left to right)
%  using polynomial sampling of order 1 with $R=5$, a wavelet frame and an $\ell_1$ regularization.}
% \label{Fig:Reconst}
% \end{center}
% \end{figure}

\begin{figure}
\begin{center}
\includegraphics[height=6cm]{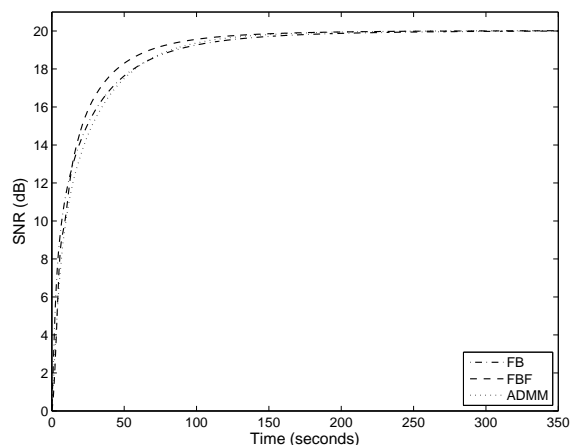}
\caption{\small Signal-to-Noise Ratio as a function of computation time using ADMM, and FB or FBF-based primal-dual methods for a given slice,
Matlab implementation on an Intel i7-3520M CPU@2.9 GHz system.}
\label{fig:convergencePMRI}
\end{center}
\vspace{-15pt}
\end{figure}

Note that convex primal-dual proximal optimization algorithms have been applied to other fields than image recovery, in particular to 
machine learning \cite{Bach_F_2012_j-ftml_optimization_sip,Mahadevan_S_2014_Proximal_rlntsdmpds}, system identification \cite{Ono_S_2013_p-icassp_sparse_siu}, audio processing \cite{Bayram_I_2014_sivp_prmal_daad}, optimal transport \cite{Papadkis_N_2014_j-siam-is_optimal_tps}, empirical mode decomposition \cite{Pustelnik_N_2014_j-sp_empirical_mdr}, seimics~\cite{Pham_M_2014_sparse_tbaf}, database management \cite{Moerkotte_G_2014_Proximal_oqf}, and data streaming over networks 
\cite{Towfic_Z_2014_Stability_plapdn}.

\subsection{Computer vision and image analysis}

The great majority of problems in computer vision involve image observation data that are of very high dimensionality, inherently ambiguous, noisy, incomplete, and often only provide a partial view of the desired space. Hence, any successful model that aims to explain such data usually requires a reasonable regularization, a robust data measure, and a compact structure between the variables of interest to efficiently characterize their relationships. Probabilistic graphical models, and in particular discrete Markov Random Fields, have led to a suitable methodology for solving such visual perception
problems \cite{Li_book, Wang}. This type of models offer great representational power, and are able to take into account dependencies in the data, encode prior knowledge, and model (soft or hard) contextual constraints in a very efficient and modular manner.  Furthermore, they offer the important ability to make  use of very powerful data likelihood terms consisting of arbitrary nonconvex and non-continuous functions that are often crucial for accurately representing the problem at hand. As a result,  MAP-inference for these models leads to discrete optimization problems that are (in most cases) highly nonconvex (NP-hard) and also of very large scale \cite{Kappes_2013, Szeliski}. {These discrete problems take the form \eqref{eq:mrf_energy_}, where typically the unary terms $\varphi_p(\cdot)$  encode the data likelihood and the higher-order terms $\thh_e(\cdot)$ encode problem specific priors.}

Primal-dual approaches can offer important computational advantages when dealing with such problems. One such characteristic example is the FastPD algorithm \cite{komodakis_CVIU08}, which currently provides one of the most efficient methods for solving generic MRF\ optimization problems of this type, also guaranteeing at the same time the  convergence to solutions that are approximately optimal.  The theoretical derivation of this method relies on the use of the primal-dual schema described in Section~\ref{se:optdisc}, which results, in this case,  into a very fast graph-cut based inference scheme that generalizes previous state-of-the-art approaches such as the $\alpha$-expansion algorithm \cite{boykov} (see Fig.~\ref{fig:fastpd}). More generally, due to the very wide applicability of MRF models to computer vision or image analysis problems, primal-dual approaches can and have been applied to a broad class of both low-level and high-level problems from these domains, including image segmentation \cite{StrandmarkKS11, DBLP:conf/iccv/VicenteKR09, DBLP:conf/cvpr/PockCCB09, DBLP:conf/iccv/WoodfordRK09}, stereo matching and 3D multi-view reconstruction \cite{Cremers_2011_book_conv_rtf, DBLP:conf/cvpr/HaneZCAP13}, graph-matching \cite{DBLP:journals/pami/TorresaniKR13}, 3D surface tracking \cite{DBLP:conf/cvpr/ZengWWGSP11}, optical flow estimation \cite{DBLP:conf/cvpr/GlockerPKTN08}, scene understanding \cite{DBLP:conf/cvpr/KumarK10}, image deblurring \cite{DBLP:conf/accv/KomodakisP12}, panoramic image stitching \cite{DBLP:journals/disopt/KolmogorovS09}, category-level segmentation \cite{Batra:2012:DMS:2403138.2403140}, and motion tracking \cite{Tsai:2012:MCT:2385520.2385534}. In the following we mention very briefly just a few examples.

\begin{figure}[t]
    \center
    \renewcommand{\figurename}{Fig.}
    \renewcommand{\captionlabelfont}{\bf}
    \renewcommand{\captionfont}{\small}
    \subfigure[`Penguin' image denoising (from left to right: noisy input image, FastPD output, time comparison plot)]
    {
    \includegraphics[height=4cm]{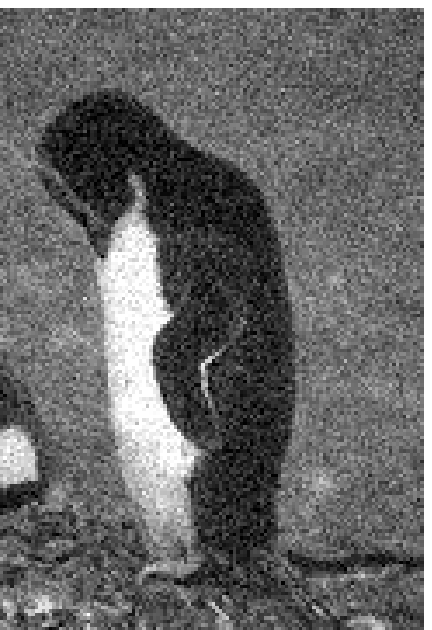} \hspace{1cm}
    \includegraphics[height=4cm]{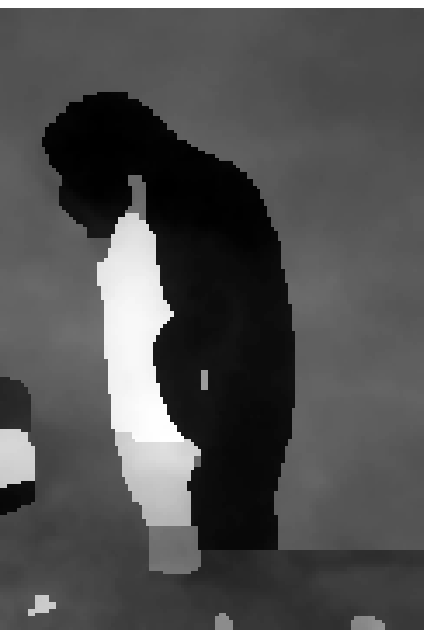} \hspace{1cm}
    \includegraphics[width = .25\linewidth]{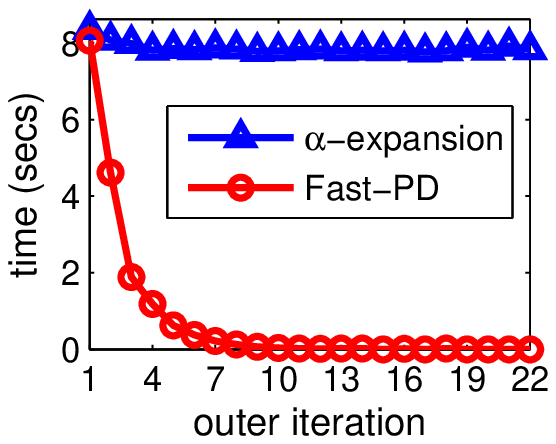}
    }
    \subfigure[`Tsukuba' stereo matching (from left to right: left image, FastPD output, time comparison plot)]
    {
    \includegraphics[width = .25\linewidth]{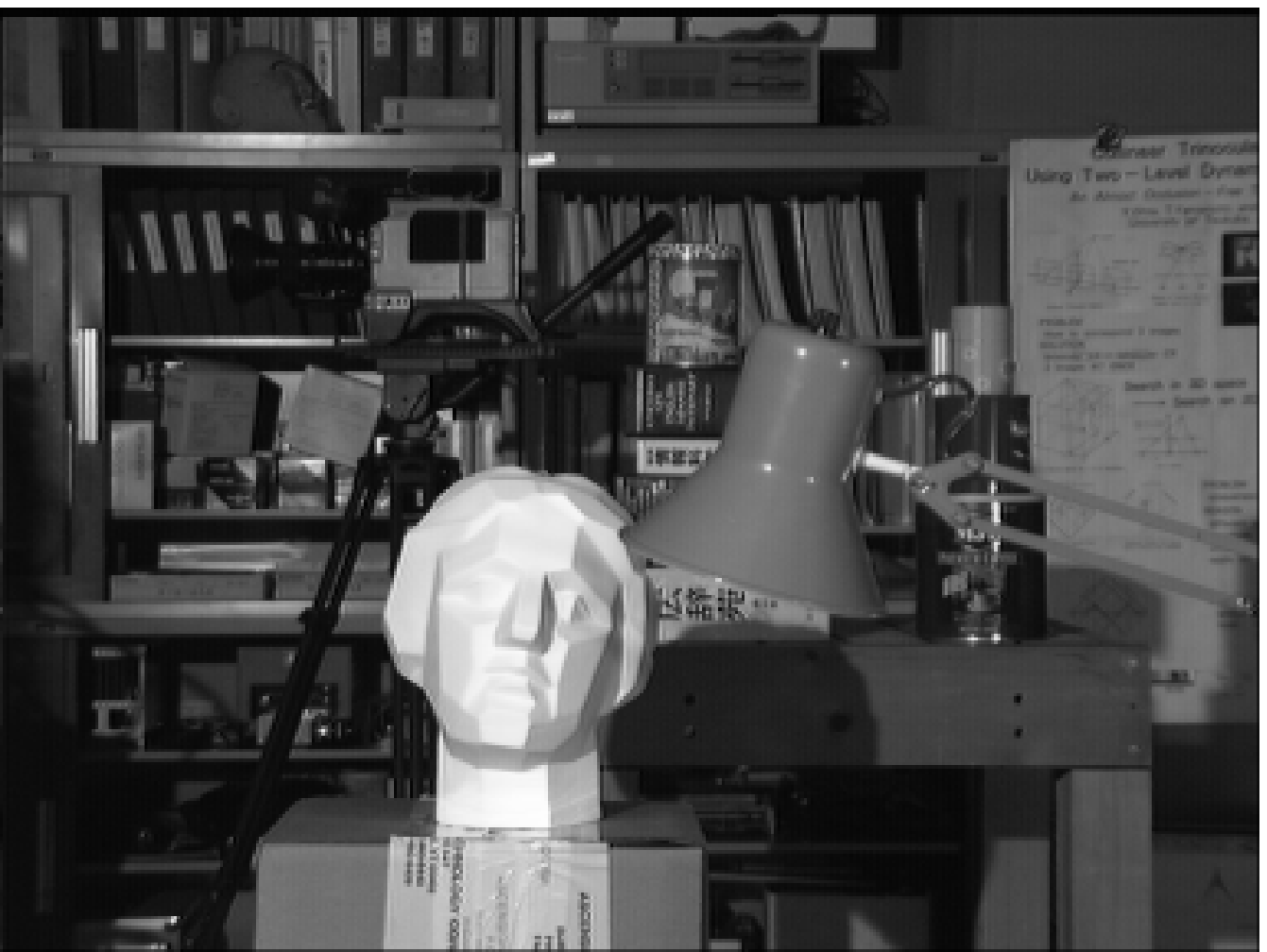} %\hspace{1cm}
    \includegraphics[width = .25\linewidth]{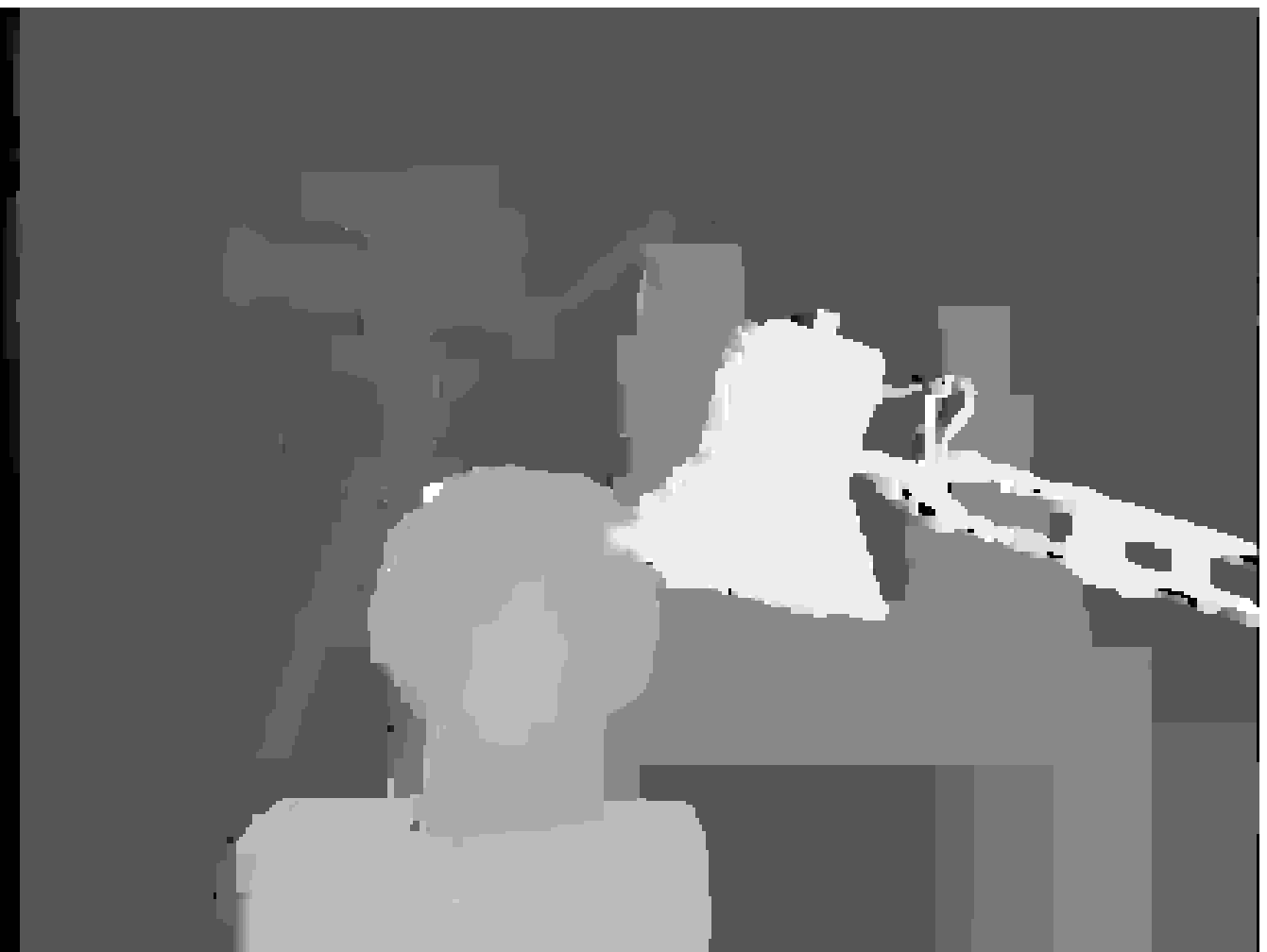} %\hspace{1cm}
    \includegraphics[width = .25\linewidth]{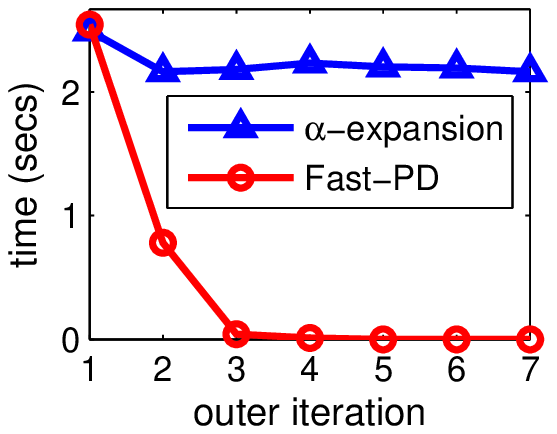}
    }
    \caption{\small FastPD \cite{komodakis_cviu} results for an image denoising (top) and stereo-matching (bottom) problem. The plot in each row  compares the corresponding running time per iteration of the above primal-dual algorithm with the $\alpha$-expansion algorithm, which is a primal-based method (experiments conducted on a 1.6 GHz CPU).
    }
    \label{fig:fastpd}
\end{figure}

A primal-dual based optimization framework has been recently proposed in \cite{DBLP:journals/mia/GlockerKTNP08,ar11} for the problem of deformable registration/fusion, which forms one of the most central and challenging tasks in medical image analysis. This  problem consists of recovering a nonlinear dense deformation field that aligns  two signals that have in general an unknown relationship both in the spatial
and intensity domain. In this framework, towards dimensionality reduction on the variables, the dense registration field is first expressed using a set of control points (registration grid) and an interpolation strategy. Then, the registration cost is expressed using a discrete sum over image costs %(representing an arbitrary similarity measure) 
projected on the control points, and a smoothness term that penalizes local deviations on the deformation field according to a neighborhood system on the grid. One advantage of the resulting optimization framework is that it is able to encode even very complex similarity measures (such as normalized mutual information and Kullback-Leibler divergence) and therefore can be used even when seeking transformations between different modalities (inter-deformable registration). Furthermore, it admits a broad range of regularization terms, and can also be applied to both 2D-2D and 3D-3D registration, as an arbitrary underlying graph structure can be readily employed (see Fig.~\ref{fig:brain} for a result on  3D inter-subject brain registration).

\begin{figure}[t]
    \center
    \renewcommand{\figurename}{Fig.}
    \renewcommand{\captionlabelfont}{\bf}
    \renewcommand{\captionfont}{\small}
    \includegraphics[width=0.75\linewidth]{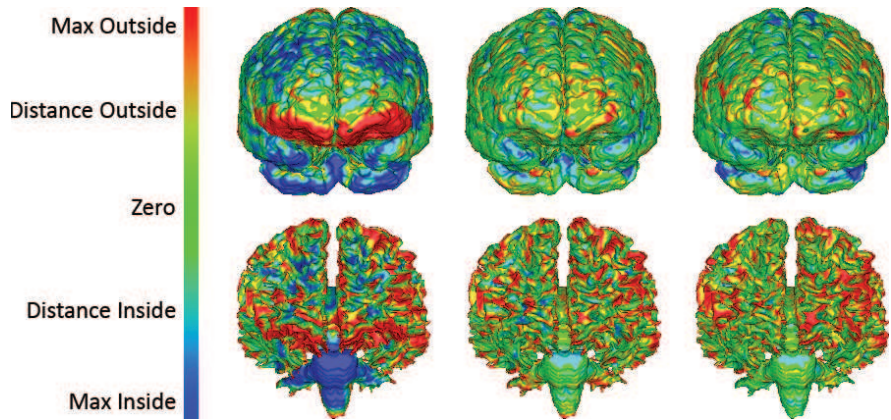}
    \caption{\small Color encoded visualization of the surface distance between warped and expert segmentation after affine (left), FFD-based \cite{Rueckert99nonrigidregistration} (middle), and primal-dual based registration (right) for the Brain 1 data set. The color range is scaled to a maximum and minimum distance of 3 mm. The  average surface distance
(ASD) after registration for the gray matter is 1.66, 1.14, and 1.00 mm for affine,
FFD-based, and primal-dual method, respectively. For the white matter the resulting
ASD is 1.92, 1.31, and 1.06 mm. Note also that the FFD-based method is more than 30 times slower than the primal-dual approach.   }
    \label{fig:brain}
    %\vspace{-10pt}
\end{figure}

Another application of primal-dual methods is in  stereo reconstruction \cite{Pock:2008:CFC:1478172.1478235}, where given as input a pair of left and right images $I_L$, $I_R$ we seek to estimate a function $u: \Omega\rightarrow \Gamma$  representing the depth $u(s)$ at a  point $s$ in the domain $\Omega\subset\mathbb{R}^2$ of the left image (here $\Gamma=[\upsilon_{\rm min}, \upsilon_{\rm max}]$ denotes the allowed depth range). To accomplish this, the following variational problem is proposed in \cite{Pock:2008:CFC:1478172.1478235}:
\begin{equation} 
\minimize{u}\int_{\Omega}f(u(s),s)ds+\int_{\Omega} |\nabla u(s)|ds,
\end{equation}
where $f(u(s),s)$ is a data term favoring different depth values by measuring  the absolute intensity differences  of respective  patches projected in the two input images, and the second term is  a TV regularizer that promotes  spatially smooth depth fields.
The above problem is nonconvex (due to the use of the data term $f$), but it turns out that  there exists an equivalent convex formulation obtained by lifting the original problem to a higher-dimensional space, in which $u$ is represented in terms of its level sets
\begin{equation}
\minimize{\phi\in D} \int_{\Sigma}(|\nabla\phi(s,\upsilon)|+f(s,\upsilon)|\partial_{\upsilon}\phi(s,\upsilon)|)ds d\upsilon.
\end{equation}
In the above formulation,  $\Sigma=\Omega\times \Gamma$, $\phi\colon\Sigma\rightarrow \{0,1\}$ is a binary function such that 
$\phi(s,\upsilon)$ equals $1$ if $u(s)>\upsilon$ and $0$ otherwise, and the feasible set is defined as $D=\left\{\phi\colon\Sigma\rightarrow \{0,1\}\mid (\forall s\in \Omega)\,\phi(s,\upsilon_{\rm min})=1, \phi(s,\upsilon_{\rm max})=0\right\}$. A convex relaxation of the latter problem is obtained by using $D'=\big\{\phi\colon\Sigma\rightarrow [0,1]\mid (\forall s\in \Omega)$\linebreak$\phi(s,\upsilon_{\rm min})=1, \phi(s,\upsilon_{\rm max})=0\big\}$ instead of $D$. A discretized form of the resulting optimization problem can be solved with the algorithms described in Section \ref{se:FBPD}. Fig.~\ref{fig:cremers} shows a sample result of this approach. 

\begin{figure}[t]
    \center
    \renewcommand{\figurename}{Fig.}
    \renewcommand{\captionlabelfont}{\bf}
    \renewcommand{\captionfont}{\small}
    \includegraphics[width=0.55\linewidth]{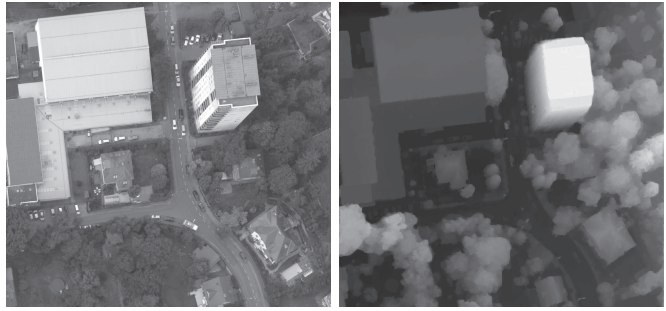}
    \vspace{3pt}\caption{\small Estimated depth map (right) for a large aerial stereo data set of Graz using the primal-dual approach in \cite{Pock:2008:CFC:1478172.1478235}. One of the images of the corresponding stereoscopic pair (of size $1500\times 1400$) is shown on the left.}
    \label{fig:cremers}
    \vspace{-8pt}
\end{figure}

Recently, primal-dual approaches have also been developed for discrete optimization problems  that involve higher-order terms \cite{DBLP:conf/cvpr/KomodakisP09, zabih_14, DBLP:conf/eccv/AroraBKM12}. They have been applied successfully to various tasks, like, for instance, in stereo matching \cite{DBLP:conf/cvpr/KomodakisP09}. In this case, apart from a data term that measures similarity between corresponding pixels in two images, a discontinuity-preserving  smoothness prior of the form $\thh(s_1,s_2,s_3)=\min(|s_1-2s_2+s_3|,\kappa)$ 
with $\kappa \in \RPP$ has been employed as a regularizer that penalizes depth surfaces of high curvature. Indicative stereo matching results from an algorithm based on the dual decomposition principle described in Section~\ref{sec:dd_sec_} are shown in Fig.~\ref{fig:ho}. 

It should be also mentioned  that an advantage of all primal-dual  algorithms (which is especially important  for NP-hard problems) is that they also provide  (for free) per-instance approximation bounds, specifying how far the  cost of an estimated solution can be from the unknown optimal cost. This directly follows from the fact that  these methods are computing both  primal and  dual solutions, which (in the case of a minimization task) provide  respectively upper and lower limits to the true optimum. These approximation bounds are continuously updated throughout an algorithm execution, and thus can be directly used for assessing the performance of a primal-dual method with respect to any particular  problem instance (and without essentially any extra computational cost).  Moreover, often in practice, these sequences converge to a common value, which means that the corresponding estimated solutions are almost optimal (see, e.g., the plots in Fig.~\ref{fig:ho}).   

\begin{figure}[t]
    \center
    \renewcommand{\figurename}{Fig.}
    \renewcommand{\captionlabelfont}{\bf}
    \renewcommand{\captionfont}{\small}
    \subfigure[`Teddy']
    {
    \includegraphics[width = 0.22\linewidth]{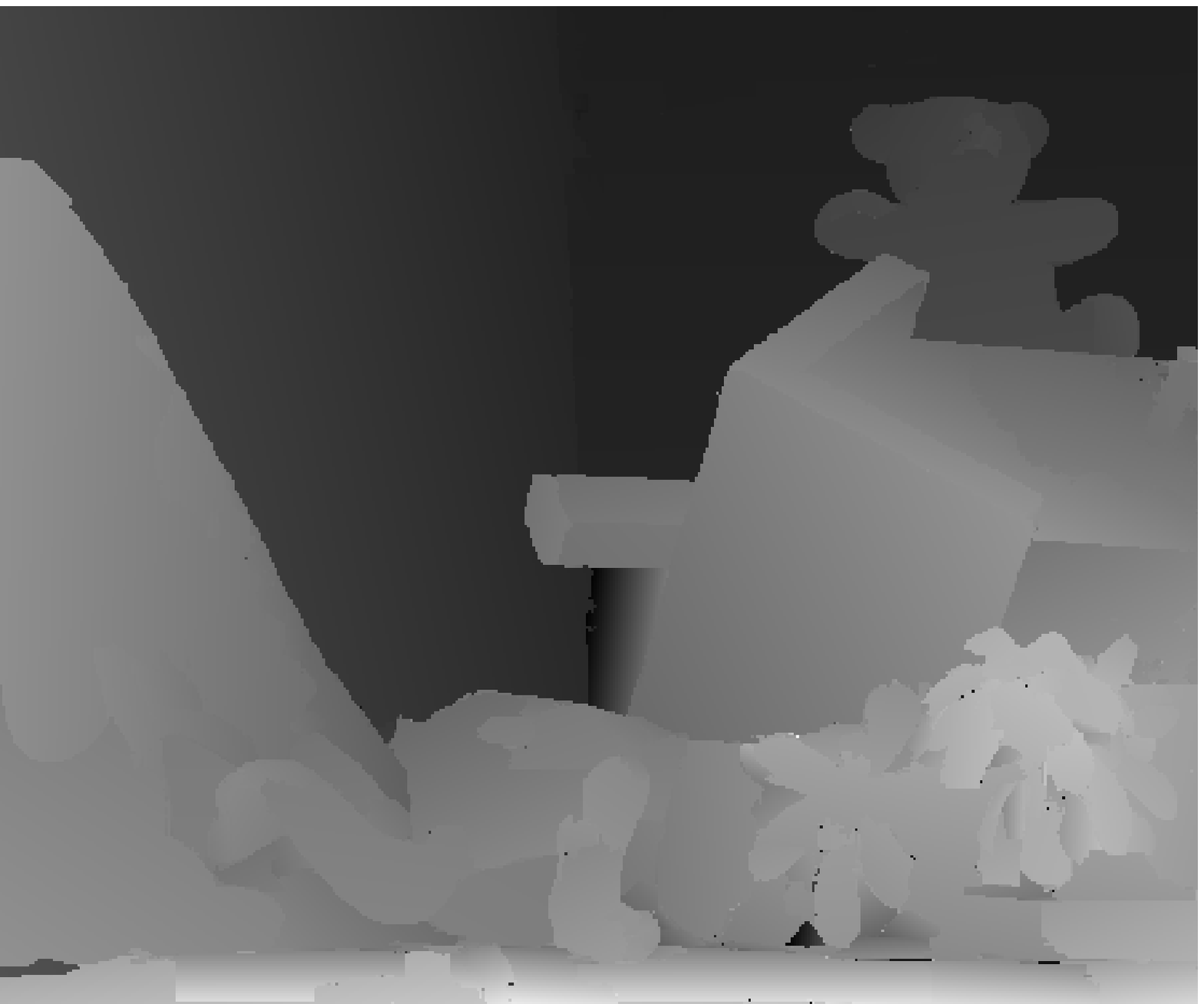}
    \includegraphics[width = 0.22\linewidth]{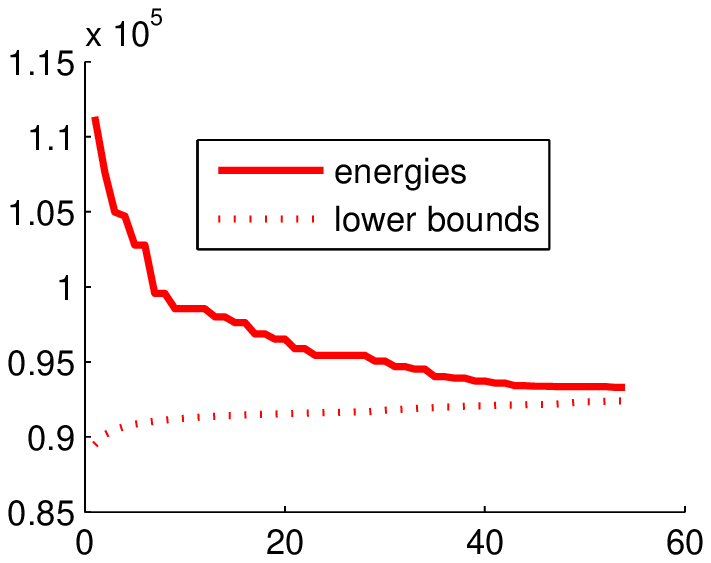}
    }
    \subfigure[`Cones']
    {
    \includegraphics[width = 0.22\linewidth]{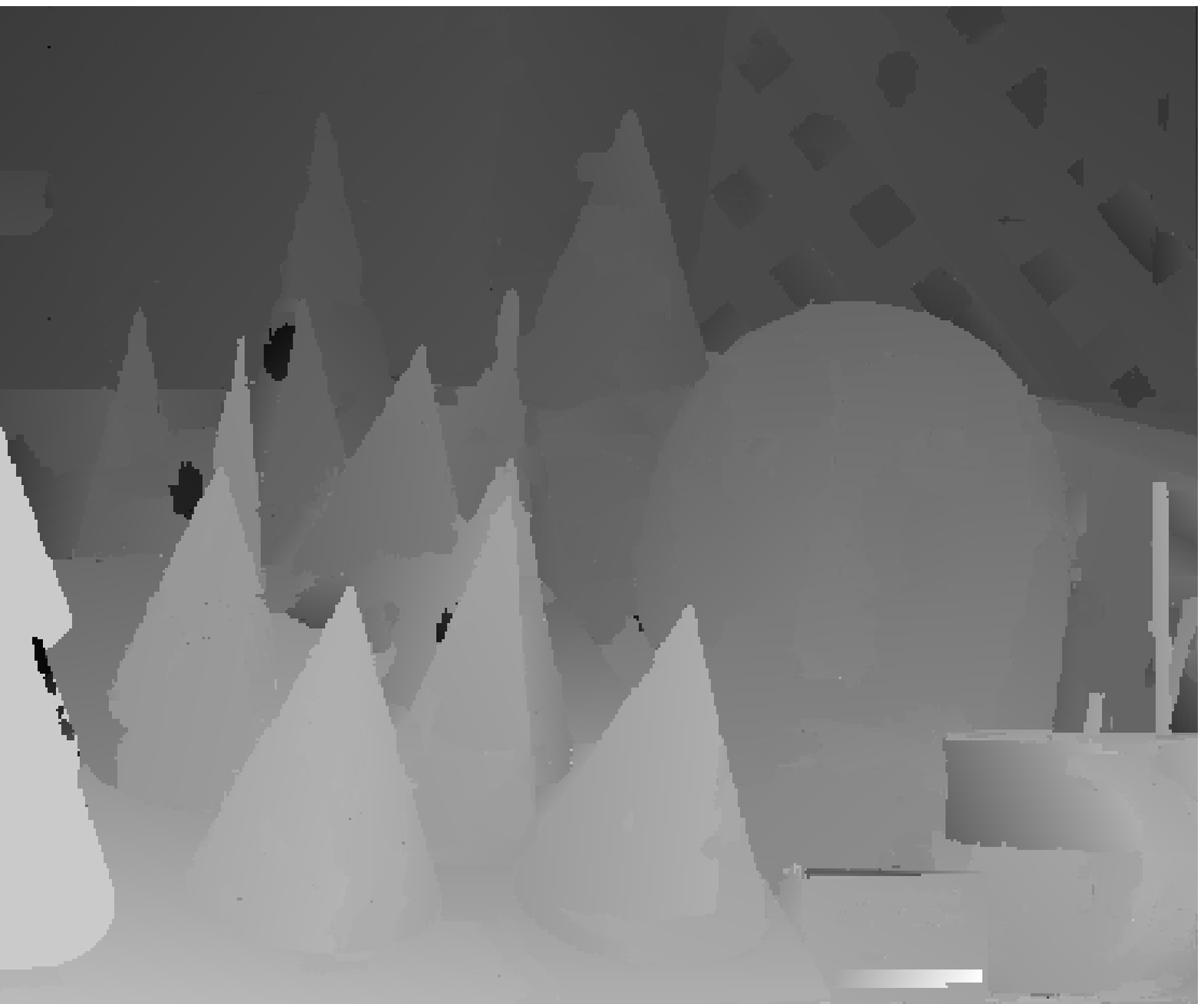}
    \includegraphics[width = 0.22\linewidth]{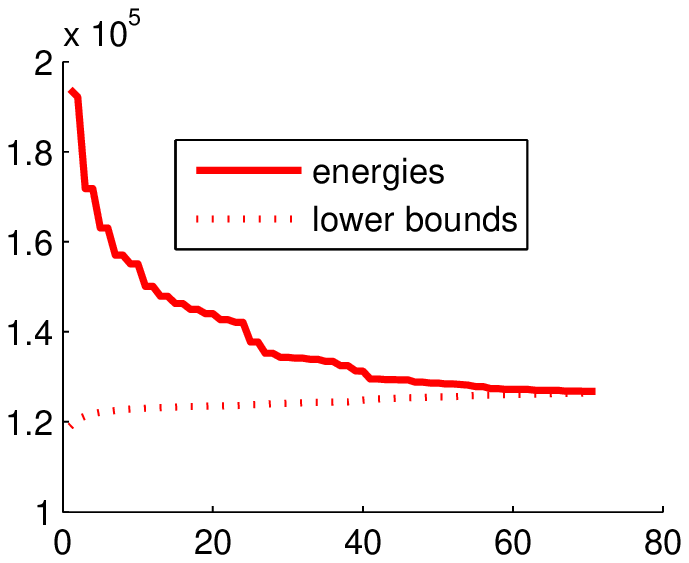}
    }
    \caption{\small Stereo matching results for `Teddy' (a) and 'Cones' (b) when using a higher-order discontinuity preserving smoothness prior. We  show plots for the corresponding sequences of upper and lower bounds generated during the primal-dual method. Notice that these sequences  converge to the same limit, meaning that the estimated solution converges to the optimal value.}
    \label{fig:ho}
    \vspace{-8pt}
\end{figure}

\section{Conclusion}\label{se:conclu}
In this paper, we have reviewed a number of primal-dual optimization methods which can be employed for solving signal and image processing problems.
The links existing between convex approaches and discrete ones were little explored in the literature and one of the contributions of this paper
is to put them in a unifying perspective. 
Although the presented algorithms have been proved to be quite effective in numerous problems, there remains
much room for extending their scope to other application fields, and also for improving them so as to accelerate their convergence. In particular, the parameter choices in these methods may have a strong influence on the convergence speed and
it would be thus interesting to design automatic procedures for setting these parameters.
Various techniques can also be devised for designing faster variants of these methods
(preconditioning, activation of blocks of variables, combination with stochastic strategies, distributed implementations...). 
Another issue to pay attention to
is the robustness to numerical errors although it can be mentioned that most of the existing proximal algorithms are tolerant to summable errors.
Concerning discrete optimization methods, we have shown that the key
to success lies in tight relaxations of combinatorial NP hard problems. Extending these methods to more challenging problems, e.g. those involving higher-order Markov fields or extremely large label sets,
appears to be of main interest in this area. More generally, developing primal-dual strategies that further bridge the gap between continuous and discrete approaches, as well as for solving other kinds of nonconvex optimization problems such as 
those encountered in phase reconstruction or blind deconvolution opens the way to appealing investigations. So, the ground is yours now to play with duality!

\vspace{-10pt}
\bibliographystyle{IEEEbib}
\bibliography{bib/abbr,bib/bibliopd_NK,bib/bibliopd}

\end{document}